\documentstyle[amscd,amssymb,verbatim,12pt]{amsart}
\newcommand{\C}{{\Bbb C}}

\newcommand{\ch}{\mbox{\rm ch}}

\newcommand{\Diff}{\mbox{\rm Diff}}

\newcommand{\End}{\mbox{\rm End}}

\newcommand{\HH}{\mbox{\rm H}}

\newcommand{\Hom}{\mbox{\rm Hom}}

\newcommand{\Image}{\mbox{\rm Im}}

\newcommand{\Ind}{\mbox{\rm Ind}}

\newcommand{\KK}{\mbox{\rm K}}
\newcommand{\Ker}{\mbox{\rm Ker}}

\newcommand{\pt}{\mbox{\rm pt}}

\newcommand{\R}{{\Bbb R}}

\newcommand{\spann}{\mbox{\rm span}}

\newcommand{\tr}{\mbox{\rm tr}}
\newcommand{\Tr}{\mbox{\rm Tr}}

\newcommand{\vol}{\mbox{\rm vol}}
\newcommand{\Z}{{\Bbb Z}}
\theoremstyle{plain}
\newtheorem{definition}{Definition}

\newtheorem{lemma}{Lemma}
\newtheorem{theorem}{Theorem}
\newtheorem{proposition}{Proposition}
\newtheorem{corollary}{Corollary}

\numberwithin{equation}{section}
\renewcommand{\rm}{\normalshape}
\begin{document}
\title{Local Index Theory over \'Etale Groupoids}
\author{Alexander Gorokhovsky}
\address{Department of Mathematics\\
University of Michigan\\
Ann Arbor, MI  48109-1109\\
USA}
\email{gorokhov@@math.lsa.umich.edu}
\author{John Lott}
\email{lott@@math.lsa.umich.edu}
\thanks{The research of the first author was partially supported by a Rackham
grant from the University of Michigan. The research of the second author
was supported by NSF grant DMS-0072154}
\date{}
\maketitle
\begin{abstract}
We give a superconnection proof of Connes' index theorem for
proper cocompact actions of \'etale groupoids. This includes
Connes' general foliation index theorem for foliations with
Hausdorff holonomy groupoid.
\end{abstract}

\section{Introduction}

This paper is concerned with a families index theorem in which the
family of operators is parametrized by 
a noncommutative space which comes from
a smooth Hausdorff \'etale groupoid $G$.  
The relevant index theorem was stated by Connes
in \cite[Section III.7.$\gamma$, Theorem 12]{Connes (1994)}. 
We give a superconnection proof of Connes' theorem.
The desirability of having such a proof was mentioned in 
\cite{Brylinski-Nistor (1994)}.
In the case of a foliation, by taking a complete transversal, one recovers
Connes' general foliation index theorem for a foliation whose
holonomy groupoid is Hausdorff. 
 For the history and significance of 
Connes' foliation index theorem we refer to 
\cite[Sections I.5, II.8-9 and III.6-7]{Connes (1994)},  
along with the references cited therein. 

For concreteness, let us first discuss the case
when $G$ is the cross-product groupoid coming from the action 
of a finitely-generated
discrete group $\Gamma$ on a smooth manifold $B$. In this case the 
geometric setup for the index theorem consists of a manifold
$\widehat{M}$ on which $\Gamma$ acts, and a submersion
$\pi \: : \: \widehat{M} \rightarrow B$ which is $\Gamma$-equivariant.
In addition, we assume that the action of $\Gamma$ on $\widehat{M}$ is
free, properly discontinuous and cocompact. Put $M \: = \: \widehat{M}/\Gamma$,
a compact manifold.  (A relevant example is when $\widehat{M} \: = \:
\R \times S^1$, $B \: = \: S^1$ and $\Gamma \: = \: \Z$, with the action of
$n \in \Z$ on $(r, e^{i \theta}) \in \R \times S^1$ given by
$n \cdot (r, e^{i \theta}) \: = \: \left(r \: + \: n, e^{i (\theta \: + \:
n \alpha)}\right)$ for some $\alpha \in \R$. Then $M \: = \: T^2$.)

There is a quotient map $M \rightarrow B/\Gamma$, which will be the
intuitive setting for our families index theorem. In general $B/\Gamma$ is
highly singular and, following Connes, we will treat it as a ``noncommutative 
space'' ${\frak B}$. 
We will give a superconnection proof of a families index theorem
in such a setting.

Let us consider two special cases.  If $\Gamma \: = \: \{e\}$ then we just
have a submersion $M \rightarrow B$ with $M$ compact. 
In this case, the relevant index theorem
is the Atiyah-Singer families index theorem \cite{Atiyah-Singer (1971)},
for which a superconnection proof was given by Bismut
\cite{Bismut (1985)}. At the 
other extreme, if $B$ is a point then we just have
a covering space $\widehat{M} \rightarrow M$ of a compact manifold $M$. 
In this case the relevant
index theorem is due to Miscenko-Fomenko \cite{Miscenko-Fomenko (1979)}
and Connes-Moscovici \cite{Connes-Moscovici (1990)}, and a superconnection
proof was given by the second author \cite{Lott (1992)}. To some degree,
the present paper combines the superconnection proofs of these two
special cases in order to deal with
general $\Gamma$ and $B$. However, there are some new features, which we
will emphasize. 

One motivation for giving a superconnection proof of
Connes' theorem is that the superconnection formalism gives a somewhat
canonical proof. In particular, it expresses the Chern character of
the index as an explicit differential form. This is one of the reasons 
that the superconnection formalism allows extensions
to the case of manifolds with boundary.  When $\Gamma \: = \: \{e\}$, this
is due to Bismut-Cheeger \cite{Bismut-Cheeger (1990)}. When $B$ is a point,
it is due to the second author \cite{Lott (1992a)} and
Leichtnam-Piazza \cite{Leichtnam-Piazza (1997)}.  Based on the present
paper, it should be possible to extend Connes' index theorem to 
manifolds-with-boundary. In particular, this would give an
index theorem for a foliated manifold-with-boundary if the foliation is 
transverse to the boundary.

We now describe in some detail
Connes' index theorem and the superconnection approach
to its proof.
Before setting up the superconnection formalism, we must first 
describe what we mean by functions and differential forms on the
noncommutative base space ${\frak B}$. There is a clear choice for a class of
``smooth functions'' on ${\frak B}$, namely the algebraic cross-product 
$C^\infty_c(B) \rtimes \Gamma$. 
The choice of ``differential forms''
on ${\frak B}$ is dictated by two facts.
First, if $\Gamma \: = \: \{e\}$, i.e. in the commutative
situation, then we want to recover the smooth compactly-supported
differential forms $\Omega^*_c(B)$. Second, the choice
should extend to general smooth Hausdorff \'etale groupoids $G$. This
dictates that we should take the differential forms to be elements
of $\Omega^*(B, \C \Gamma) \: = \:
\Omega^*_c(B) \: \widehat{\otimes} \: \Omega^*(\C \Gamma)$, where
$\Omega^*(\C \Gamma)$ is the graded differential algebra of noncommutative
forms on $\C \Gamma$ \cite[Section 2.6]{Loday (1998)}, 
and the product in $\Omega^*(B, \C \Gamma)$ takes into
account that $\Gamma$ acts on $\Omega^*_c(B)$. 

Given these choices, we need
to know that the ensuing ``homology'' of ${\frak B}$ is sufficiently rich.
This is shown in the following theorem. Let $GT_{*,\langle e \rangle}$ denote
the complex of
graded traces on $\Omega^*(B, \C \Gamma)$ that are concentrated
at the identity conjugacy class in $\Gamma$. 
(Only these graded traces will be
relevant for the paper.)
We let $C_*(B)$ denote the
currents on $B$ and we let $\overline{\cal C}_*(\Gamma)$ denote a certain 
complex of differential forms on $E \Gamma$,
described in (\ref{1.47.5}). \\ \\
{\bf Theorem 1 \: : \:} The homology of $GT_{*,\langle e \rangle}$ is
isomorphic to the homology of the total complex of the double complex
$\left( (\C \oplus\overline{\cal C}_*(\Gamma)) 
\otimes C_*(B) \right)^\Gamma$. \\

In particular, if $\eta$ is a closed graded $n$-trace on 
$\Omega^*(B, \C \Gamma)$, concentrated
at the identity conjugacy class, then we obtain a corresponding
cohomology class
$\Phi_\eta \in \HH^{n \: + \: \dim(B) \: + \: 2 \Z}_\tau 
\left((E\Gamma \times B)/ \Gamma) \right)$, where the $2 \Z$ denotes
an even-odd grading and the $\tau$ denotes a twisting by the orientation
bundle of $B$. This shows our choice of differential forms gives the
``right'' answer, as $(E\Gamma \times B)/ \Gamma$ is the classifying space
$BG$ for the groupoid $G$; compare the cyclic cohomology calculations in
\cite{Brylinski-Nistor (1994)},
\cite[Section III.2.$\delta$]{Connes (1994)},
\cite{Crainic (1999)},
\cite{Crainic-Moerdijk (2000)}.

For analytic reasons arising from finite propagation speed estimates, 
we introduce a slightly larger space of differential
forms. Let $\parallel \cdot \parallel$ be a word-length metric on $\Gamma$
and put 
\begin{equation}
{\cal B}^\omega \: = \:
\left\{ \sum_{\gamma \in \Gamma} c_\gamma \: \gamma \: : \:
|c_\gamma| \text{ decays faster than any exponential in }\parallel \gamma
\parallel \right\}.
\end{equation}
Put $C^\infty(B, {\cal B}^\omega) \: = \: {\cal B}^\omega \otimes_{\C \Gamma}
(C^\infty_c(B) \rtimes \Gamma)$. It contains
$C^\infty_c(B) \rtimes \Gamma$ as a dense subset. 
There is a corresponding space of
differential forms
$\Omega^*(B, {\cal B}^\omega)$, defined in (\ref{2.30}).

Let us now state Connes' index theorem.
Let $\widehat{M}$ and $M$ be as above. The $\Gamma$-covering $\widehat{M}
\rightarrow M$ is classified by a $\Gamma$-equivariant continuous map
$\mu \: : \: \widehat{M} \rightarrow E\Gamma$, defined up to 
$\Gamma$-homotopy.  Let $\widehat{\nu} \: : \: \widehat{M} \rightarrow
E\Gamma \times B$ be $(\mu, \pi)$ and let $\nu \: : \: M \rightarrow
(E\Gamma \times B)/\Gamma$ be the $\Gamma$-quotient of $\widehat{\nu}$, 
a classifying map for the action of the groupoid $G$ on $\widehat{M}$.

Let $Z$ denote a fiber of
the submersion $\pi \: : \: \widehat{M} \rightarrow B$ and let $TZ$ denote
the vertical tangent bundle, a $\Gamma$-invariant 
vector bundle on $\widehat{M}$.
There is a foliation ${\cal F}$ on $M$ whose leaves are the images
of the fibers $Z$ under $\pi$. The tangent bundle to the foliation is
$T{\cal F} \: = \: (TZ)/\Gamma$, a vector bundle on $M$. Choose a smooth 
$\Gamma$-invariant vertical Riemannian metric $g^{TZ}$ on 
$\widehat{M}$. Assume that $TZ$ has a $\Gamma$-invariant
spin structure, with corresponding
spinor bundle $S^Z$. Let $\widehat{V}$ be a $\Gamma$-invariant vector bundle
on $\widehat{M}$, with $\Gamma$-invariant Hermitian 
connection $\nabla^{\widehat{V}}$. 
Put $V \: = \: \widehat{V}/\Gamma$, a vector bundle on $M$, and
put $\widehat{E} \: = \: S^Z \: \otimes \: \widehat{V}$.

There is an ensuing $\Gamma$-invariant family $D$ of Dirac-type operators
which act fiberwise on $C^\infty_c(\widehat{M}; \widehat{E})$. Equivalently,
the family $D$ is $G$-invariant. Let
$\Ind(D)$ denote the index of this family; we will say more about it later. 
Let $\eta$ be a closed graded trace on
$\Omega^*(B, {\cal B}^\omega)$ which is 
concentrated at the identity conjugacy class. 
In particular, $\eta$ restricts to a closed graded trace on 
$\Omega^*(B, \C \Gamma)$.
In this situation, Connes' index theorem
\cite[Section III.7.$\gamma$, Theorem 12]{Connes (1994)} becomes the statement
that
\begin{equation} \label{0.1}
\langle \ch(\Ind(D)), \eta \rangle \: = \:
\int_M \widehat{A}(T{\cal F}) \: \ch(V) \: \nu^*\Phi_\eta.
\end{equation}
As mentioned before, special cases of Connes' index theorem are the
Atiyah-Singer families index theorem and the covering-space index theorem.

The goal of the paper is to give a superconnection proof of (\ref{0.1}).
A key ingredient will be 
the $\C \Gamma$-vector bundle ${\cal E} \: = \: (\widehat{M} 
\times \C \Gamma)/\Gamma$ on $M$, where $\Gamma$ acts diagonally on
$\widehat{M} \times \C \Gamma$. By construction,
$C^\infty(M; {\cal E}) \: \cong \: C^\infty_c(\widehat{M})$.
The natural flat connection
$\nabla^{1,0}$ on ${\cal E}$ sends $f \in C^\infty_c(\widehat{M})$ to
$df$. An important part of our proof is a certain differentiation
\begin{equation}
\nabla^{0,1} \: : \: C^\infty_c(\widehat{M}) \rightarrow \Omega^1(\C \Gamma)
\otimes_{\C \Gamma} C^\infty_c(\widehat{M})
\end{equation}
in the ``noncommutative''
directions. The explicit formula for $\nabla^{0,1}$ is given in 
Section \ref{Section 2}. 
The sum of $\nabla^{1,0}$ and $\nabla^{0,1}$ is a nonflat
connection
\begin{equation}
\nabla^{can} \: : \: 
C^\infty_c(\widehat{M}) \rightarrow \Omega^1(\widehat{M}, \C \Gamma)
\otimes_{C^\infty_c(\widehat{M}) \rtimes \Gamma} C^\infty_c(\widehat{M}).
\end{equation}

Choose a $\Gamma$-invariant horizontal distribution $T^H \widehat{M}$ on
$\widehat{M}$.
Suppose that $Z$ is even-dimensional.
For $s \: > \: 0$, there is an ensuing Bismut superconnection
$A^{Bismut}_s$ on the submersion $\widehat{M} \rightarrow B$
\cite[Section 10.3]{Berline-Getzler-Vergne (1992)},
\cite[Section IIIa]{Bismut (1985)}.
Our noncommutative superconnection is simply
$A_s \: = \: A^{Bismut}_s \: + \: \nabla^{0,1}$.
Let $\Omega^*(B, {\cal B}^\omega)_{ab}$ be the quotient of 
$\Omega^*(B, {\cal B}^\omega)$ by the closure of its graded commutator.
We show that the Chern character $\Tr_{s,<e>} \left( e^{- \: A_s^2} \right)$
of the superconnection $A_s$
is well-defined in
$\Omega^*(B, {\cal B}^\omega)_{ab}$.  It is a closed form
and its cohomology class
is independent of $s$. The next result gives its $s \rightarrow 0$ limit.

Let $\phi \in C^\infty_c(\widehat{M})$ be such that
$\sum_{\gamma \in \Gamma} \gamma \cdot \phi \: = \: 1$.
Let ${\cal R}$ be the rescaling operator on 
$\Omega^{even}(B, {\cal B}^\omega)_{ab}$ which multiplies an element of
$\Omega^{2k}(B, {\cal B}^\omega)_{ab}$ by 
$(2 \pi i)^{-k}$.\\ \\
{\bf Theorem 2 : } 
\begin{equation}
\lim_{s \rightarrow 0} {\cal R} \: 
\Tr_{s, <e>} \left( e^{- \: A_s^2} \right) \: = \:
\int_Z \phi \: \widehat{A} \left( \nabla^{TZ} \right) \: 
\ch \left( \nabla^{\widehat{V}} \right) \: \ch \left(\nabla^{can} \right)
\: \in \: \Omega^*(B, {\cal B}^\omega)_{ab}.
\end{equation}

The proof of Theorem 2 is by local index theory techniques.

In order to prove (\ref{0.1}), it remains to relate
$\left\langle \Tr_{s, <e>} \left( e^{- \: A_s^2} \right), \eta
\right\rangle$ to
$\left\langle \ch(\Ind(D)), \eta
\right\rangle$. 
Let $C_0(B) \rtimes_r \Gamma$ be the reduced cross-product $C^*$-algebra.
A technical problem is that the Dirac-type operator $D$, when
considered as an operator on a $C_0(B) \rtimes_r \Gamma$-Hilbert 
module, may not have closed range.
This problem also arises in the superconnection proof of the families
index theorem \cite{Bismut (1985)}. Recall that in the
families index theorem, a special case is when 
the kernels and the cokernels of the
fiberwise operators form vector bundles on the base. In this case,
one defines the
analytic index to be the difference of these two vector bundles, as an
element of the $K$-theory of the base. If the
kernels and the cokernels do not form vector bundles then one can deform
the fiberwise operators in order to reduce oneself to the case in which
they do \cite{Atiyah-Singer (1971)}. 

In order to carry out this
deformation argument at the level of superconnection Chern characters
requires a pseudodifferential operator calculus 
\cite[Section 2d]{Bismut (1985)}.
In our context
one can set up such a calculus for a class of operators on 
$C_0(B) \rtimes_r \Gamma$-vector
bundles. However, this would not be enough for our purposes, as
we would need such a calculus for operators on 
$C^\infty(B, {\cal B}^\omega)$-vector bundles.
There seem to be serious problems in constructing such a calculus, as
$C^\infty(B, {\cal B}^\omega)$ is generally not
closed under the holomorphic functional calculus in $C_0(B) \rtimes_r \Gamma$.

To get around this problem, we use a method which seems to be new even
in the case of the families index theorem.  The idea, which is due to
Nistor \cite{Nistor (1997)}, is to define $\Ind(D)$ to be the $\KK$-theory
element represented by the difference between the
index
projection $p$ and a standard projection $p_0$, 
and then relate the Chern character of $[p \: - \: p_0]$ 
to the superconnection
Chern character.  In order to relate the two, Nistor works in a universal
setting and shows that a certain cyclic cohomology group is singly
generated, which implies that the two expressions are related by a
computable constant.  Unfortunately, Nistor's assumptions do not
hold in our setting and his argument does not seem to be adaptable.  
Instead, we give a direct proof relating $[p \: - \: p_0]$
to the superconnection Chern character.  We show that the pairing of
$\eta$ with $[p \: - \: p_0]$ 
can be written as the pairing of $\eta$ with the
Chern character of a certain $\Z_2$-graded
connection $\nabla^\prime$.
We then homotop between the connection $\nabla^\prime$ 
and the superconnection $A_s$. Of course 
one cannot do so in a purely formal way,
as the vector bundles involved are infinite-dimensional.  
However, we show that one can
write things so that one has uniformly smoothing operators inside the
traces during the homotopy, thereby justifying the formal argument.
In this way we prove\\ \\
{\bf Theorem 3 : } For all $s \: > \: 0$,
\begin{equation}
\langle \ch(\Ind(D)), \eta \rangle \: = \:
\langle {\cal R} \: 
\Tr_{s, <e>} \left( e^{- \: A_s^2} \right), \eta \rangle \: \in \: \C.
\end{equation}

Equation (\ref{0.1}) follows from combining Theorems 2 and 3.

Equation (\ref{0.1}) deals with closed graded traces $\eta$ on 
$\Omega^*(B, {\cal B}^\omega)$. On the other hand,
the geometric and
topological consequences of index theory involve the K-theory of
$C_0(B) \rtimes_r \Gamma$.
Equation (\ref{0.1}) is in some ways a stronger result than the analogous 
$\KK_*(C_0(B) \rtimes_r \Gamma)$-valued index theorem, 
as the underlying algebra in
(\ref{0.1}) is
$C^\infty(B, {\cal B}^\omega) \subset C_0(B) \rtimes_r \Gamma$.
However, in order
to obtain geometric consequences from (\ref{0.1}), one must consider
certain ``smooth'' algebras ${\cal A}$ such that $C^\infty(B, {\cal B}^\omega)
\subset {\cal A} \subset C_0(B) \rtimes_r \Gamma$. The requirements
on ${\cal A}$ are that it should be closed under the holomorphic
functional calculus in $C_0(B) \rtimes_r \Gamma$, and that $\eta$ should
extend to a continous cyclic cocycle on ${\cal A}$. The existence of
such subalgebras ${\cal A}$ is discussed in 
\cite[Chapter III]{Connes (1994)} and we do not have anything new to
say about it in this paper.

We give the extensions of 
Theorems \ref{Theorem 2} and \ref{Theorem 3} to the setting of 
smooth Hausdorff \'etale 
groupoids. 
As mentioned above, our proof applies to foliations, to give a 
superconnection proof of Connes' general
foliation index theorem (for foliations with Hausdorff holonomy groupoid).  
There has been
some previous work along these lines.  Theorem 2
was proven by Heitsch, using $A_s^{Bismut}$, in the special case 
when $\eta$ comes from a holonomy-invariant transverse current to the
foliation \cite{Heitsch (1995)}.
A corresponding analog of Theorem 3 was proven by Heitsch and Lazarov
when $\eta$ comes from a 
holonomy-invariant transverse current and under some additional
technical assumptions regarding the spectral densities of the leafwise
Dirac-type operators \cite{Heitsch-Lazarov (1999)}. 
In \cite{Liu-Zhang (1999)}, Liu and Zhang gave
an adiabatic limit proof of a certain case of a vanishing result of Connes 
for foliations with spin leaves of positive scalar curvature.  (Connes'
result is a corollary of his
general foliation index theorem.)  The additional
assumptions in \cite{Liu-Zhang (1999)} were that the foliation is
almost-Riemannian and that the pairing object comes from the
Pontryagin classes of
the normal bundle. The adiabatic limit is closely related to superconnections.

The methods of \cite{Heitsch (1995)},
\cite{Heitsch-Lazarov (1999)} and \cite{Liu-Zhang (1999)}
have the limitation that the dimension of the pairing objects
is at most the codimension of the foliation, as they come from the normal
bundle to the foliation.  Consequently, one misses the
noncommutative features of the foliation, which lead to the phenomenon that
the dimension of the leaf space, treated as a noncommutative space, can be
greater than the codimension of the foliation.  For example, a manifold with 
a codimension-$1$ foliation has a Godbillon-Vey
class which is a three-dimensional cohomology class.
One important aspect of Connes' foliation index theorem is that it allows
a pairing between the foliation index and the Godbillon-Vey class
\cite[Section III.6]{Connes (1994)}. One would miss this pairing
if one treated the foliation in a
``commutative'' way.

Let us also mention the paper \cite{Moriyoshi-Natsume (1996)} which
proves (\ref{0.1}) in the case when $\widehat{M} \: = \: Z \times S^1$,
$B \: = \: S^1$, the action of $\Gamma$ on $S^1$ preserves orientation,
$V$ is trivial and $\eta$ corresponds to the Godbillon-Vey
class.  The method of proof of \cite{Moriyoshi-Natsume (1996)} is
to represent $\Ind(D)$ by means of a ``graph'' projection and then
compute the pairing of this projection with a cyclic cocycle representing
the Godbillon-Vey class.

 Not all foliations have a holonomy groupoid that is Hausdorff.
We expect that our results can be extended to the nonHausdorff case, but
we have not worked this out in detail.  Relevant treatments of the
cyclic cohomology of an \'etale groupoid algebra, in the nonHausdorff case,
are in \cite{Crainic (1999)}, \cite{Crainic-Moerdijk (2000)} and
\cite{Crainic-Moerdijk (2000b)}. 

The results of this paper are valid in the generality of a smooth Hausdorff
\'etale groupoid $G$ with a free proper 
cocompact action on a manifold $P$.  
For simplicity of notation, we first present 
all of the arguments in the case of a cross-product groupoid 
$B \rtimes \Gamma$.
We then explain how to extend the proofs to general $G$.

The structure of the paper is as follows.  In Section \ref{Section 2} 
we prove
Theorem 1. In Section \ref{Section 3} 
we define a class of fiberwise-smoothing operators
on $\widehat{M} \rightarrow B$ and construct their
$\Omega^*(B, {\cal B}^\omega)_{ab}$-valued traces.  In Section 
\ref{Section 4} we
define the superconnection $A_s$ and prove Theorem 2. In Section 
\ref{Section 5} we prove Theorem 3 and
give some consequences.  

In Section 
\ref{Section 6} we give the extension of the preceding results to general
\'etale groupoids.  The extension is not entirely straightforward and the
expressions in Section \ref{Section 6} 
would be unmotivated if it were not for the
results of the preceding sections.
In Subsection \ref{Subsection 6.1}
we show how in the case of a cross-product groupoid,
the expressions of Section \ref{Section 6} 
reduce to expressions of the preceding sections.
The reader may wish to read Subsection 
\ref{Subsection 6.1} simultaneously with the
rest of Section \ref{Section 6}.

More detailed information is given at the beginnings of the sections.

Background information on superconnections and index theory is in 
\cite{Berline-Getzler-Vergne (1992)}. For background results we sometimes
refer to the relevant sections in
\cite{Berline-Getzler-Vergne (1992)}, \cite{Connes (1994)} or
\cite{Loday (1998)}, where references to the original articles can be found.

We thank Victor Nistor for discussions of  
\cite{Brylinski-Nistor (1994)} and \cite{Nistor (1997)}. 
We thank the referee for useful suggestions,
and Eric Leichtnam and Paolo Piazza for corrections to an earlier
version of this paper. 
The second author thanks the Max-Planck-Institut-Bonn and the
Mathematical Sciences Research Institute for their hospitality while
part of this research was performed.

\section{Closed Graded Traces} \label{Section 2}

In this section we compute the homology of $GT_{*,<e>}$. We show in
Proposition \ref{Proposition 1} that it is 
isomorphic to the homology of the total complex of a certain double complex
${\cal T}_{*,*}$.
We then construct a morphism from this double complex to the
double complex $\left( (\C \oplus\overline{\cal C}_*(\Gamma)) 
\otimes C_*(B) \right)^\Gamma$. In order to construct this morphism, we
first describe the connection $\nabla^{can}$ on ${\cal E} \: = \:
(\widehat{M} \times \C \Gamma)/\Gamma$, along with its
Chern character. Taking $\widehat{M}$ to be $E \Gamma$, we obtain a
connection $\nabla^{univ}$, whose Chern character implements the morphism.
We show that the morphism induces an isomorphism between the $E^1$-terms of
the two double complexes.  It follows that it induces an isomorphism
between the homologies of the total complexes.  

 The material in this section is necessarily of a technical
nature.  A  trusting reader may be willing to take the main result of the
section, Corollary \ref{Corollary 1}, on faith.  To read the rest of the
paper, it is also worthwhile to read the Important Digression of the
section, which is labeled as such. 

Let ${\cal B}$ be a unital algebra over $\C$.
Let $\Omega^*({\cal B})$ denote its universal graded differential
algebra (GDA) \cite[Section 2.24]{Karoubi (1987)}. As a vector space,
\begin{equation} \label{1.1}
\Omega^k({\cal B}) \: = \: {\cal B} \: \otimes \: \left( \otimes^k
({\cal B}/\C) \right).
\end{equation}
As a GDA, $\Omega^*({\cal B})$ is generated by ${\cal B}$ and $d {\cal B}$
with the relations
\begin{equation} \label{1.2}
d1 \: = \: 0, \: \: \: d^2 \: = \: 0, \: \: \:
d(\omega_k \omega_l) \: = \: (d\omega_k) \omega_l \: + \: (-1)^k \:
\omega_k \: (d\omega_l)
\end{equation}
for $\omega_k \: \in \Omega^k({\cal B})$, $\omega_l \: \in 
\Omega^l({\cal B})$. 
It will be convenient to write an element $\omega_k$ of
$\Omega^k({\cal B})$ as a finite sum 
$\omega_k \: = \: \sum b_0 \: db_1 \ldots db_k$. 
Let $[\Omega^*({\cal B}), \Omega^*({\cal B})]$ denote the graded
commutator of $\Omega^*({\cal B})$ and
let $\Omega^*({\cal B})_{ab} \: = \:
\frac{\Omega^*({\cal B})}{[\Omega^*({\cal B}), \Omega^*({\cal B})]}$ 
denote the abelianization of $\Omega^*({\cal B})$.
Put $\overline{\Omega}^*({\cal B})_{ab} \: = \: 
\frac{\Omega^*({\cal B})}{\C \: \oplus \:
[\Omega^*({\cal B}), \Omega^*({\cal B})]}$.
Let $HDR^*({\cal B})$ and $\overline{H}DR^*({\cal B})$ denote the
cohomologies of the differential complexes
$\Omega^*({\cal B})_{ab}$ and $\overline{\Omega}^*({\cal B})_{ab}$,
respectively.
Then $HDR^*({\cal B}) \cong \overline{H}DR^*({\cal B})$ if $* \: \ge \: 1$,
and there is a short exact sequence
$0 \rightarrow \C \rightarrow HDR^0({\cal B}) \rightarrow 
\overline{H}DR^0({\cal B}) \rightarrow 0$. Furthermore, Connes and Karoubi
showed that $\overline{H}DR^*({\cal B})$ is expressed in terms of reduced
cyclic and Hochschild homology 
\cite[Theorem 2.6.7]{Loday (1998)} by
\begin{equation} \label{1.3}
\overline{H}DR^*({\cal B}) \: \cong \: \Ker \left(
B \: : \: \overline{HC}_*({\cal B}) \rightarrow
\overline{HH}_{*+1}({\cal B}) \right).
\end{equation}

If $\Gamma$ is a discrete group and ${\cal B}$ is the group algebra
$\C \Gamma$,
i.e. finite sums $\sum_{\gamma \in \Gamma} c_\gamma \gamma$,
let us recall the calculation of $HDR^*(\C \Gamma)$.
It breaks up with respect to the conjugacy classes of $\Gamma$, as do
the Hochschild and cyclic cohomologies of $\C \Gamma$, and we will
only be interested in the component $HDR^*_{<e>}(\C \Gamma)$ corresponding
to the identity conjugacy class.
Let $\HH_*(\Gamma; \C)$ denote the group homology of $\Gamma$ and let
$\overline{\HH}_*(\Gamma; \C)$ denote the reduced group homology,
i.e. $\overline{\HH}_*(\Gamma; \C) \: = \: \HH_*(B\Gamma, *; \C)$.
Then it follows from Burghelea's work that
the reduced Hochschild and cyclic homologies of $\C \Gamma$, when
considered at the identity conjugacy class, are 
\begin{align} \label{1.4}
\overline{HH}_{*,{<e>}}(\C \Gamma) & \cong \overline{\HH}_{*}(\Gamma; \C)
\text{ and }\notag \\
\overline{HC}_{*,{<e>}}(\C \Gamma) & \cong \bigoplus_{i \ge 0} 
\overline{\HH}_{*-2i}(\Gamma; \C),
\end{align}
with the map
$B \: : \: \overline{HC}_{*,{<e>}}(\C \Gamma) \rightarrow
\overline{HH}_{*+1,{<e>}}(\C \Gamma)$ 
vanishing \cite[Section 7.4]{Loday (1998)}.
Hence
\begin{equation} \label{1.5}
HDR^*_{<e>}(\C \Gamma) \cong
\begin{cases}
\HH_0(\Gamma; \C) & \text{ if $* \: = \: 0$,} \\
\bigoplus_{i \ge 0} \overline{\HH}_{*-2i}(\Gamma; \C) & 
\text{ if $* \: > \: 0$.}
\end{cases}
\end{equation}

If ${\cal B}$ is a locally convex topological algebra
then there is a natural
completion of the algebraic GDA $\Omega^*({\cal B})$ to a locally convex GDA
\cite[Section 5.1]{Karoubi (1987)}. For simplicity of notation, 
when the context is clear we will also denote
this completion by $\Omega^*({\cal B})$. We will also denote by
$\Omega^*({\cal B})_{ab}$ the quotient of $\Omega^*({\cal B})$ by the
closure of $[\Omega^*({\cal B}), \Omega^*({\cal B})]$, where the
quotienting by the closure is done in order to obtain a Hausdorff space.
In general,
we take the tensor product of two locally convex topological vector
spaces to be the projective
topological tensor product and 
we let $\widehat{\otimes}$ denote a graded (projective) tensor product.

Let $B$ be a smooth manifold on which a
finitely-generated discrete group $\Gamma$ acts on the right, not necessarily
freely or properly discontinuously.  Given $\gamma \in \Gamma$, let
$R_\gamma \in \Diff(B)$ denote the action of $\gamma$ on $B$.
We let $\Gamma$ act
on $C^\infty_c(B)$ on the left so that
$\gamma \cdot f \: = \: R_\gamma^* f$, i.e. $(\gamma \cdot f)(b) \: = \:
f(b \gamma)$.
Given $\gamma \in \Gamma$,
let $B^\gamma$ denote the subset of $B$ which is pointwise fixed by $\gamma$.
Given $b \in B$, let $\Gamma_b \subset \Gamma$ be the isotropy subgroup 
at $b$ for the action of $\Gamma$ on $B$. 

Let $C^\infty_c(B) \rtimes \Gamma$ denote the cross product
algebra, whose elements are finite sums $\sum_{\gamma \in \Gamma} f_\gamma \: 
\gamma$,
with $f_\gamma \in C^\infty_c(B)$. 
We wish to define an appropriate GDA whose
zeroth-degree component equals $C^\infty_c(B) \rtimes \Gamma$. One's first 
choice might be the universal GDA $\Omega^*(C^\infty_c(B) \rtimes \Gamma)$.
The corresponding de Rham cohomology is expressed in (\ref{1.3}) 
in terms of the
cyclic and Hochschild homologies of $C^\infty_c(B) \rtimes \Gamma$.
Now such homologies have been computed in 
\cite[Section III.2.$\delta$]{Connes (1994)},
\cite{Brylinski-Nistor (1994)},
\cite{Crainic (1999)} and \cite{Crainic-Moerdijk (2000)}
in cases of increasing generality. In particular,
the periodic cyclic cohomology of $C^\infty_c(B) \rtimes \Gamma$ contains
a factor consisting of the
twisted equivariant cohomology of $B$. This implies that
the space of closed graded traces on $\Omega^*(C^\infty_c(B) \rtimes \Gamma)$
is sufficiently rich for our purposes.

On the other hand, the choice of $\Omega^*(C^\infty_c(B) \rtimes \Gamma)$ 
as the GDA is inconvenient from the
point of view of superconnections.  In the case when $\Gamma$ is the
trivial group, and one is dealing with a fiber bundle $M \rightarrow B$, the
usual superconnection formalism uses the graded differential algebra
$\Omega^*_c(B)$. It appears that it would be quite cumbersome to redo
the superconnection proof of the Atiyah-Singer families index theorem using
the noncommutative differential forms
$\Omega^*(C^\infty_c(B))$ instead of $\Omega^*_c(B)$.

For this reason, instead of using $\Omega^*(C^\infty_c(B) \rtimes \Gamma)$
as the GDA, we replace it by an appropriate quotient. By universality,
there will be 
a map from the closed graded traces on the quotient GDA to the closed graded
traces on $\Omega^*(C^\infty_c(B) \rtimes \Gamma)$. We want this map to
have a sufficiently large image. It turns out that an 
appropriate GDA is
$\Omega^*(B, \C \Gamma) \: = \: \Omega^*_c(B) 
\widehat{\otimes} \Omega^*(\C \Gamma)$,
the graded algebraic tensor product over $\C$,
where the multiplication in $\Omega^*(B, \C \Gamma)$ takes into account that
$\C \Gamma$ acts on $\Omega^*_c(B)$. Then $\Omega^*(B, \C \Gamma)$ is a GDA,
with $\Omega^0(B, \C \Gamma) \: = \: C^\infty_c(B) \rtimes \Gamma$.

We wish to compute the homology of the complex of graded traces on
$\Omega^*(B, \C \Gamma)$, or at least the part of the homology which is
concentrated at the identity conjugacy class.
Let $GT_{n,<e>}$ be the graded traces on $\Omega^*(B, \C \Gamma)$ which
are concentrated on the elements $
\sum_{\gamma_0, \ldots, \gamma_n : \gamma_0 \cdot \ldots \cdot
\gamma_n = e}
\omega_{\gamma_0, \ldots, \gamma_n} \gamma_0 \: d\gamma_1 \ldots d\gamma_n$.
Let $d^t : GT_{n,<e>} \rightarrow GT_{n-1,<e>}$ be the boundary operator.
A closed graded trace is an element of $\Ker(d^t)$. 

Consider the space of maps $\tau_k \: : \: \Gamma^{k+1}
\rightarrow \C$.  We now define certain operators that arise in
the computation of the cyclic homology of $\C \Gamma$
\cite[Section 2.21]{Karoubi (1987)}. Namely, put 
\begin{align} \label{1.6}
(t \tau)_k(\gamma_0, \gamma_1, \ldots, \gamma_k) \: & = \:
(-1)^k \: \tau_k(\gamma_1, \ldots, \gamma_k, \gamma_0), \notag \\
(\delta_i \tau)_{k+1}(\gamma_0, \ldots, \gamma_{k+1}) \: & = \:
\tau_k(\gamma_0, \ldots, \gamma_{i-1}, \widehat{\gamma_i}, \gamma_{i+1},
\ldots, \gamma_k), \notag \\
(\sigma_i \tau)_{k-1}(\gamma_0, \ldots, \gamma_{k-1}) \: & = \:
\tau_k(\gamma_0, \ldots, \gamma_i, \gamma_i, \ldots, 
\gamma_{k-1}), \notag \\
(B_0 \tau)_{k-1}(\gamma_0, \ldots, \gamma_{k-1}) \: & = \:
\tau_k(\gamma_{k-1}, \gamma_0, \ldots, \gamma_{k-1}).
\end{align}
Put 
\begin{equation} \label{1.7}
(b \tau)_{k+1} \: = \: \sum_{i=0}^k (-1)^i \: (\delta_i \tau)_{k+1}.
\end{equation}
The operator $B_0 \: b \: + \: b \: B_0$ acts by
\begin{align} \label{1.8}
((B_0 \: b \: + \: b \: B_0) \tau)_k(\gamma_0, \ldots, \gamma_k) \:  = \:
& \tau_k(\gamma_0, \ldots, \gamma_k) \: - \: (-1)^k \: 
\tau_k(\gamma_k, \gamma_0, \ldots, \gamma_{k-1}) \notag \\
& + \: (-1)^k \: \tau_k(\gamma_{k-1}, \gamma_0, \ldots, \gamma_{k-1}).
\end{align}

Let $C_*(B)$ be the currents on $B$. We denote the pairing of
$\omega \in \Omega^*_c(B)$ with $c \in C_*(B)$ by $\langle \omega, c \rangle$.
The de Rham boundary $\partial \: : \: C_*(B) \rightarrow C_{*-1}(B)$ satisfies
$\langle d\omega, c \rangle \: = \: \langle \omega, \partial c \rangle$.
The action of $\Gamma$ on $C_*(B)$ is such that
$\langle \gamma \cdot \omega, \gamma \cdot c \rangle \: = \: 
\langle \omega, c \rangle$.
 
Let $C^k(\Gamma)$ be the vector space of $\C$-valued functions on 
$\Gamma^{k+1}$.
Let ${\cal C}_{k,l}$ be the elements $\tau_k \in C^k(\Gamma) \otimes C_l(B)$
which are
normalized in the sense that $(\sigma_i \tau)_{k-1} \: = \: 0$ for all 
$0 \: \le \: i \: < \: k$,
and which are $\Gamma$-invariant in the sense that
\begin{equation} \label{1.9}
\tau_k(\gamma \gamma_0, \ldots, \gamma \gamma_k) \: = \:
\gamma \cdot \tau_k(\gamma_0, \ldots, \gamma_k).
\end{equation}
Put ${\cal C}_n \: = \: \bigoplus_{k \: + \: l \: = \: n}{\cal C}_{k, l}$.
The operators $\partial$, $b$ and $B_0$ can be defined on ${\cal C}_{*,*}$ 
in the natural way. Let ${\cal T}_{*, *}$ be the double complex given by
\begin{equation} \label{1.10}
{\cal T}_{k,l} \: = \: \Ker(b) \cap \Ker(bB_0 \: + \: B_0 b) \subset 
{\cal C}_{k,l}
\end{equation}
with boundary operators $(-1)^{l} \: B_0$ and $\partial$.

\begin{proposition} \label{Proposition 1}
The vector space $GT_{n,<e>}$ is isomorphic to 
$\bigoplus_{k \: + \: l \: = \: n}{\cal T}_{k, l}$.
Under this isomorphism, the action of $d^t$ on 
$GT_{n,<e>}$ is
equivalent to the action of $\partial \: + \: (-1)^{l} \: B_0$ on
$\bigoplus_{k \: + \: l \: = \: n}{\cal T}_{k, l}$.
\end{proposition}
\begin{pf}
Given $\tau \in {\cal C}_{n}$, write it as
$\tau \: = \: \sum_{k=0}^n \tau_k$, with $\tau_k \in {\cal C}_{k,n-k}$. 
We obtain a linear functional $\Psi_\tau$ on
$\Omega^*(B, \C \Gamma)$, which is concentrated at the identity
conjugacy class, by the
formula
\begin{equation} \label{1.11}
\Psi_\tau(\omega \: \gamma_0 d\gamma_1 \ldots d\gamma_k) \: = \: 
\begin{cases}
\langle \omega, \tau_k(\gamma_0, \gamma_0 \gamma_1, \ldots,
\gamma_0 \ldots \gamma_{k-1}, e) \rangle & \text{ if
$\gamma_0 \gamma_1 \ldots \gamma_k \: = \: e$}, \\
0 & \text{ if $\gamma_0 \gamma_1 \ldots \gamma_k \: \neq \: e$}.
\end{cases}
\end{equation}
Conversely, all linear functionals on 
$\Omega^*(B, \C \Gamma)$, which are concentrated at the identity
conjugacy class, arise in this way.

$\Psi_\tau$ will be a graded trace if and only if it satisfies
\begin{equation} \label{1.12}
\Psi_\tau(\gamma \: \omega \:  \gamma_0 d\gamma_1 d\gamma_2 \ldots d\gamma_k)
 \: = \:
\Psi_\tau(\omega \: \gamma_0 d\gamma_1 d\gamma_2 \ldots d\gamma_k \: \gamma)
\end{equation}
and
\begin{equation} \label{1.13}
\Psi_\tau(d\gamma_k \: \omega \:  \gamma_0 d\gamma_1 d\gamma_2 \ldots 
d\gamma_{k-1}) \: = \:
(-1)^{k-1 + |\omega|} \: 
\Psi_\tau(\omega \: \gamma_0 d\gamma_1 d\gamma_2 \ldots 
d\gamma_k).
\end{equation}

Equation (\ref{1.12}) is equivalent to
\begin{align} \label{1.14}
\Psi_\tau((\gamma \cdot \omega) \: \gamma \gamma_0 \:
d\gamma_1 d\gamma_2 \ldots d\gamma_k) \: = \: 
\Psi_\tau(\omega \: [& 
\gamma_0 d\gamma_1 d\gamma_2 \ldots d(\gamma_k \gamma)
\: - \: \gamma_0 d\gamma_1 d\gamma_2 \ldots d(\gamma_{k-1} \gamma_k) d \gamma)
\notag \\
&
+ \ldots + \: (-1)^k \: \gamma_0 \gamma_1 \: d\gamma_2 \ldots d\gamma_k
d\gamma]).
\end{align}
This is equivalent to
\begin{align} \label{1.15}
& \langle \gamma \cdot \omega, \tau(\gamma \gamma_0,  \:
\gamma \gamma_0 \gamma_1, \ldots, \gamma \gamma_0 \ldots \gamma_{k-1}, e)
\rangle \: = \notag \\
& \langle \omega, \tau(\gamma_0, \gamma_0 \gamma_1, \ldots, 
\gamma_0 \ldots \gamma_{k-1}, e)  
\: - \: \tau(\gamma_0, \gamma_0 \gamma_1, \ldots, \gamma_0 \ldots \gamma_{k-2},
\gamma_0 \ldots \gamma_k, e) \: + \ldots + \notag \\
& \: \: \: \: \: (-1)^k \:
\tau(\gamma_0 \gamma_1, \ldots, \gamma_0 \ldots \gamma_k, e)\rangle,
\end{align}
which is equivalent to
\begin{align} \label{1.16}
& \tau(\gamma_0,  \: \gamma_0 \gamma_1, \ldots, \gamma_0 \ldots \gamma_{k-1}, 
\gamma_0 \ldots \gamma_{k})
\: = \notag \\
& \tau(\gamma_0, \gamma_0 \gamma_1, \ldots, 
\gamma_0  \ldots \gamma_{k-1}, e)  
\: - \: \tau(\gamma_0, \gamma_0 \gamma_1, \ldots, \gamma_0 \ldots \gamma_{k-2},
\gamma_0 \ldots \gamma_k, e) \: + \ldots + \notag \\
& \: \: \: \: \: (-1)^k \:
\tau(\gamma_0 \gamma_1, \ldots, \gamma_0 \ldots \gamma_k, e),
\end{align}
which is equivalent to
\begin{align} \label{1.17}
0 \: = \: & (-1)^k \: \left[ \tau(\gamma_0 \gamma_1, \ldots, \gamma_0 \ldots 
\gamma_k, e) \: + \ldots + (-1)^{k-3} \: 
\tau(\gamma_0, \gamma_0 \gamma_1, \ldots, \gamma_0 \ldots \gamma_{k-2},
\gamma_0 \ldots \gamma_k, e) \: +  \right. \notag \\
& \left.  (-1)^{k-2} \: 
\tau(\gamma_0, \gamma_0 \gamma_1, \ldots, 
\gamma_0 \ldots \gamma_{k-1}, e) \: + \:
(-1)^{k-1} \: \tau(\gamma_0,  \: \gamma_0 \gamma_1, \ldots, \gamma_0 \ldots 
\gamma_{k-1}, \gamma_0 \ldots \gamma_{k}) \right].
\end{align}
After a change of variable and using the $\Gamma$-invariance of $\tau$,
this is equivalent to $b \tau \:= \: 0$.

Next,
\begin{align} \label{1.18}
& (-1)^{|\omega|} \:
\Psi_\tau(d\gamma_k \: \omega \:  \gamma_0 d\gamma_1 d\gamma_2 \ldots 
d\gamma_{k-1}) \:  = \:
\Psi_\tau((\gamma_k \cdot \omega) \: d\gamma_k \:  \gamma_0 d\gamma_1 
d\gamma_2 \ldots d\gamma_{k-1}) \: =  \notag \\
& \Psi_\tau((\gamma_k \cdot \omega) \: d(\gamma_k \gamma_0) d\gamma_1 d\gamma_2
 \ldots d\gamma_{k-1})
\: - \:
\Psi_\tau((\gamma_k \cdot \omega) \: \gamma_k \:  d\gamma_0 d\gamma_1 
d\gamma_2 \ldots d\gamma_{k-1}) \: = \notag \\
& \langle \gamma_k \cdot \omega, \tau_k(e, \gamma_k \gamma_0, 
\gamma_k \gamma_0 \gamma_1, \ldots,
\gamma_k \gamma_0 \ldots \gamma_{k-2}, e) \rangle \notag \\
& \: \: \: \: \: - \:
\langle \gamma_k \cdot \omega, \tau_k(\gamma_k, \gamma_k \gamma_0, 
\gamma_k \gamma_0 \gamma_1, \ldots,
\gamma_k \gamma_0 \ldots \gamma_{k-2}, e) \rangle \: = \: \notag \\
& \langle \omega, \tau_k(\gamma_0 \ldots \gamma_{k-1}, \gamma_0, 
\gamma_0 \gamma_1, \ldots,
\gamma_0 \ldots \gamma_{k-1}) \: - \:
\tau_k(e, \gamma_0, 
\gamma_0 \gamma_1, \ldots,
\gamma_0 \ldots \gamma_{k-1}) \rangle.
\end{align}
Thus (\ref{1.13}) is equivalent to
\begin{align} \label{1.19}
& \tau_k(\gamma_0 \ldots \gamma_{k-1}, \gamma_0, 
\gamma_0 \gamma_1, \ldots,
\gamma_0 \ldots \gamma_{k-1}) \: - \: \tau_k(e, \gamma_0, 
\gamma_0 \gamma_1, \ldots,
\gamma_0 \ldots \gamma_{k-1}) \: = \notag \\
& \: \: \: \: \: (-1)^{k-1} \:
\tau_k(\gamma_0, \gamma_0 \gamma_1, \ldots,
\gamma_0 \ldots \gamma_{k-1}, e). 
\end{align}
By a change of variables, this in turn is equivalent to
\begin{align} \label{1.20}
& \tau_k(g_{k-1}, g_0, 
g_1, \ldots,
g_{k-1}) \: - \: \tau_k(e, g_0,  \ldots,
g_{k-1}) \: = \notag \\
& \: \: \: \: \: (-1)^{k-1} \:
\tau_k(g_0, \ldots,
g_{k-1}, e). 
\end{align}
Using the $\Gamma$-equivariance, this is equivalent to
\begin{align} \label{1.21}
& \tau_k(g_{k-1}, g_0, 
g_1, \ldots,
g_{k-1}) \: - \: \tau_k(g_k, g_0,  \ldots,
g_{k-1}) \: = \notag \\
& \: \: \: \: \: (-1)^{k-1} \:
\tau_k(g_0, \ldots,
g_{k-1}, g_k), 
\end{align}
which, from (\ref{1.8}), is equivalent to 
$(B_0 b \: + \: b B_0) \tau \: = \: 0$.

Finally, given $\tau \in C^k(\Gamma) \otimes C_l(B)$, let
$\omega \: \gamma_0 \: d\gamma_1 \ldots d\gamma_{k^\prime}$ be an element of
$\Omega^{l^\prime}(B) \: \widehat{\otimes} \Omega^{k^\prime}(\C \Gamma)$.
Then
\begin{equation} \label{1.22}
(d^t \Psi_\tau) (\omega \: \gamma_0 \: d\gamma_1 \ldots d\gamma_{k^\prime}) 
\: = \:
\Psi_\tau \left(d \omega \: \gamma_0 \: d\gamma_1 \ldots d\gamma_{k^\prime} 
\: + \:
(-1)^{l^\prime} \: \omega \: d\gamma_0 \ldots d\gamma_{k^\prime} \right).
\end{equation}
If $k^\prime \: = \: k$ and $l^\prime \:= \: l \: - \: 1$ then
\begin{align} \label{1.23} 
(d^t \Psi_\tau) (\omega \: \gamma_0 \: d\gamma_1 \ldots d\gamma_{k^\prime}) 
\: & = \:
\langle d\omega, \tau(\gamma_0, \gamma_0 \gamma_1, \ldots, \gamma_0 \ldots, 
\gamma_{k-1}, e) \rangle
\notag \\
\: & = \:
\langle \omega, \partial
\tau(\gamma_0, \gamma_0 \gamma_1, \ldots, \gamma_0 \ldots, \gamma_{k-1}, e) 
\rangle \notag \\
&  = \: \Psi_{\partial \tau} (\omega \: \gamma_0 \: d\gamma_1 \ldots 
d\gamma_{k^\prime}).
\end{align}
If $k^\prime \: = \: k \: - \: 1$ and $l^\prime \: = \: l$ then
\begin{align} \label{1.24}
(d^t \Psi_\tau) (\omega \: \gamma_0 \: d\gamma_1 \ldots d\gamma_{k^\prime}) 
\: & = \:
(-1)^{l^\prime} \: \langle \omega, \tau(e, \gamma_0, \ldots, 
\gamma_0 \ldots \gamma_{k^\prime-1}, e) \rangle \notag \\ 
& = \: (-1)^{l^\prime} \: \langle \omega, (B_0 \tau)(\gamma_0, \ldots, 
\gamma_0 \ldots \gamma_{k^\prime-1}, e) \rangle \notag \\
& = \:
\Psi_{(-1)^{l} B_0 \tau}(\omega \: \gamma_0 \: d\gamma_1 \ldots 
d\gamma_{k^\prime}).
\end{align}
This proves the proposition.
\end{pf}

From Proposition \ref{Proposition 1}, (\ref{1.7}) and (\ref{1.8}), we see
that any antisymmetric group cocycle for $\Gamma$,
which takes values in the closed currents on $B$,
gives a closed graded trace on $\Omega^*(B, \C \Gamma)$.  
In this way, we have a map
$\HH^k(\Gamma; Z_l(B)) \rightarrow \HH_{k+l}(GT_{*,\langle e \rangle},
d^t)$. In particular, if $k \: = \: 0$ then we obtain closed graded traces
on $\Omega^*(B, \C \Gamma)$
from $\Gamma$-invariant closed currents on $B$. We now use Proposition
\ref{Proposition 1} to describe all of the homology of the complex
$GT_{*,\langle e \rangle}$.

Consider the $E^1$-term of the double complex ${\cal T}_{*, *}$.  That is,
$E^1_{k,l}$ is the $k$-th homology group of the complex ${\cal T}_{*, l}$ with
respect to the differential $(-1)^l \: B_0$.
We first want to compute this homology group. 
To do so, we follow the general method of proof of
\cite[Section III.1.$\beta$, Theorem 22]{Connes (1994)}.
We fix $l$ for the moment.

Let us define operators $b^\prime$, $A$ and $B$
on $C^*(\Gamma) \: \otimes \: C_l(B)$ by the
usual formulas 
\begin{align} \label{1.25}
(b^\prime \tau)_{k+1} \: & = \: \sum_{i=0}^{k-1} (-1)^i \:
(\delta_i \tau)_{k+1}, \notag \\
(A \tau)_k \: & = \: \sum_{i=0}^k (-1)^i \: t^i \tau, \notag \\
B \: & = \: A \: B_0.
\end{align}
(The ``$B$'' in this Connes $B$-operator should not be 
confused with the manifold $B$.)

\begin{lemma} \label{Lemma 1}
Acting on ${\cal C}_{*,l}$,
we have
$\Ker(b) \cap \Ker(B_0) \: = \: \Ker(b) \cap \Ker(1 \: - \: t)$.
\end{lemma}
\begin{pf}
If $\tau \in {\cal C}_{*,l}$ and $(1 \: - \: t) \tau \: = \: 0$ then
$B_0 \tau \: = \: 0$. Thus $\Ker(b) \cap \Ker(1 \: - \: t) \subset
\Ker(b) \cap \Ker(B_0)$. On the other hand, the identity
$b^\prime \: B_0 \: + \: B_0 \: b \: = \: 1 \: - \: t$ shows that
$\Ker(b) \cap \Ker(B_0) \: \subset \: \Ker(b) \cap \Ker(1 \: - \: t)$.
\end{pf}

\begin{lemma} \label{Lemma 2}
An element $\tau \in {\cal T}_{k,l}$ lies in the image of
$B_0 \: : \: {\cal T}_{k+1,l} \rightarrow {\cal T}_{k,l}$ if and only
if $\tau \: = \: B \phi$ for some
$\phi \in \Ker(b) \subset {\cal C}_{k+1, l}$.
\end{lemma}
\begin{pf}
Suppose that $\tau \in {\cal T}_{k,l}$ satisfies
$\tau \: = \: B_0 \phi$ for some $\phi \in 
{\cal T}_{k+1,l}$. Then $B_0 \tau \: = \: 0$. By the previous lemma
$\tau \: = \: t \tau$ and so $\tau \: = \: \frac{1}{k \: + \: 1} \: 
A \tau \: = \:  \frac{1}{k \: + \: 1} \: B \phi$. 

Now suppose that $\tau \in {\cal T}_{k,l}$ satisfies 
$\tau \: = \: B \phi$ for some
$\phi \in \Ker(b) \subset {\cal C}_{k+1, l}$. Put
\begin{equation} \label{1.26}
\Theta \: = \: B_0 \phi \: - \: \frac{1}{k+1} \tau.
\end{equation} 
Note that $B_0 \phi \in \Ker(B_0) \subset {\cal C}_{k, l}$.
Also, as $\tau \: = \: B \phi$, it follows that $\tau \: = \: t \tau$ which,
along with the fact that $\tau \in {\cal T}_{k,l}$, implies that
$\tau \in \Ker(B_0)$. Hence 
$\Theta \in \Ker(B_0) \subset {\cal C}_{k, l}$.
As $A \Theta \: = \: 0$, we can write
$\Theta \: = \: \psi \: - \: t \psi$ where
$\psi \in \Ker(B_0) \subset {\cal C}_{k, l}$ is given by
$\psi \: = \: - \: \frac{1}{k+1} \sum_{i=0}^k i \: t^i \Theta$.
Then $\Theta \: = \: (1 \: - \: t) \: \psi \: = \: 
(b^\prime \: B_0 \: + \: B_0 \: b) \: \psi \: = \: B_0 b \psi$.
Put $\phi^\prime \: = \: \phi \: - \: b\psi$. Then from (\ref{1.26}),
$\tau \: = \: (k \: + \: 1) \: B_0 \phi^\prime$. Furthermore,
$b \phi^\prime \: = \: b (\phi \: - \: b\psi) \: = \: 
b\phi \: = \: 0$ and $(k \: + \: 1) \: b B_0 \phi^\prime \: = \: 
b \tau \: = \: 0$. Hence $\phi^\prime \in {\cal T}_{k+1,l}$.
This proves the lemma.
\end{pf}

Put 
\begin{equation} \label{1.27}
C^k_\lambda(\Gamma; C_l(B)) \: = \: \Ker(B_0) \cap \Ker(B_0 b \: + \:
b B_0) \subset {\cal C}_{k,l}.
\end{equation}
From (\ref{1.6}) and (\ref{1.8}), the elements of $C^k_\lambda(\Gamma; C_l(B))$
can be considered to be cyclic cochains which are reduced if $k \: > \: 0$. 
Put 
\begin{align} \label{1.28}
Z^k_\lambda(\Gamma; C_l(B)) \: & = \: 
\Ker(b \: : \: C^k_\lambda(\Gamma; C_l(B)) \rightarrow
C^{k+1}_\lambda(\Gamma; C_l(B)), \notag \\ 
B^k_\lambda(\Gamma; C_l(B)) \: & = \: 
\Image(b \: : \: C^{k-1}_\lambda(\Gamma; C_l(B)) \rightarrow
C^{k}_\lambda(\Gamma; C_l(B)), \notag \\
H^k_\lambda(\Gamma; C_l(B)) \: &= \: Z^k_\lambda(\Gamma; C_l(B))/
B^k_\lambda(\Gamma; C_l(B)).
\end{align}
We define the Hochschild objects
$CH^k(\Gamma; C_l(B))$, $ZH^k(\Gamma; C_l(B))$, $BH^k(\Gamma; C_l(B))$ and  
$HH^k(\Gamma; C_l(B))$ similarly, but without the cyclic
condition.  The Connes $B$-operator induces a map
$B \: : \: HH^k(\Gamma; C_l(B)) \rightarrow
H^{k-1}_\lambda(\Gamma; C_l(B))$.

We now prove a result which, when $B$ is a point, amounts to the dual
of (\ref{1.3}), when applied to ${\cal B} \: = \: \C \Gamma$ and
considered at the identity conjugacy class.

\begin{lemma} \label{Lemma 3}
There is an isomorphism
\begin{equation} \label{1.29}
\HH^k({\cal T}_{*, l}) \: \cong \: H^{k}_\lambda(\Gamma; C_l(B))/
\Image( B \: : \: HH^{k+1}(\Gamma; C_l(B)) \rightarrow
H^{k}_\lambda(\Gamma; C_l(B)).
\end{equation}
\end{lemma}
\begin{pf}
By Lemma \ref{Lemma 1},
$\Ker(B_0) \subset {\cal T}_{k,l}$ is isomorphic to
$Z^k_\lambda(\Gamma; C_l(B))$. By Lemma \ref{Lemma 2},
\begin{equation} \label{1.30}
\Image \left( B_0 \: : \: 
{\cal T}_{k+1,l} \rightarrow {\cal T}_{k,l} \right) \:
= \: \Image \left( B \: : \: ZH^{k+1}(\Gamma; C_l(B)) \rightarrow 
{\cal T}_{k,l} \right).
\end{equation}
Also,
\begin{equation} \label{1.31}
B^k_\lambda(\Gamma; C_l(B)) \subset
\Image \left( B \: : \: ZH^{k+1}(\Gamma; C_l(B)) \rightarrow 
{\cal T}_{k,l} \right),
\end{equation}
as if $\tau \: = \: b \psi$ with
$\psi \in C^{k-1}_\lambda(\Gamma; C_l(B))$ then 
$\psi \: = \: B \psi^\prime$ for some
$\psi^\prime \in CH^{k}(\Gamma; C_l(B))$
\cite[Section III.1.$\beta$, Corollary 20]{Connes (1994)}, and so
$\tau \: = \: bB \psi^\prime \: = \: - \: Bb\psi^\prime$.
The lemma follows.
\end{pf}

Let ${\HH}^*(\Gamma; C_l(B))$ denote the group cohomology
of $\Gamma$ with coefficients in the $\Gamma$-module $C_l(B)$.
Let $\overline{\HH}^*(\Gamma; C_l(B))$ denote the reduced group cohomology
with coefficients in $C_l(B)$. That is, let $E\Gamma_{(0)}$ be the set of
vertices in $E\Gamma$. Then $\overline{\HH}^*(\Gamma; C_l(B))$ is the
cohomology of the complex
$\left( C^*(E\Gamma, E\Gamma_{(0)}) \otimes_{\Z} C_l(B) \right)^\Gamma$.

We now give a result which, when $B$ is a point, amounts to the dual
of (\ref{1.4}).

\begin{lemma} \label{Lemma 4}
\begin{equation} \label{1.32}
HH^{k}(\Gamma; C_l(B)) \cong
\begin{cases}
\HH^{0}(\Gamma; C_l(B)) & \text{ if $k \: = \: 0$}, \\
\overline{\HH}^{k}(\Gamma; C_l(B)) & \text{ if $k \: > \: 0$}
\end{cases}
\end{equation}
and
\begin{equation} \label{1.33}
H_\lambda^{k}(\Gamma; C_l(B)) \cong
\begin{cases}
\HH^{0}(\Gamma; C_l(B)) & \text{ if $k \: = \: 0$}, \\
\bigoplus_{i \ge 0} \overline{\HH}^{k-2i}(\Gamma; C_l(B)) & \text{ if 
$k \: > \: 0$},
\end{cases}
\end{equation}
with $B \: : \: 
HH^{k+1}(\Gamma; C_l(B)) \rightarrow H_\lambda^{k}(\Gamma; C_l(B))$
vanishing.
\end{lemma}
\begin{pf}
The proof of this follows algebraically from the method of proof of
(\ref{1.4}). That is, we have the same sort of
cyclic structures.  We omit the details.
\end{pf}

Putting together Lemmas \ref{Lemma 3} and \ref{Lemma 4}, we have shown
\begin{proposition} \label{Proposition 2}
The $E^1$-term of ${\cal T}_{*,*}$ is given by
\begin{equation} \label{1.34}
E^1_{k,l} \: \cong \: 
\begin{cases}
\HH^{0}(\Gamma; C_l(B)) & \text{ if $k \: = \: 0$}, \\
\bigoplus_{i \ge 0} \overline{\HH}^{k-2i}(\Gamma; C_l(B)) & 
\text{ if $k \: > \: 0$.}
\end{cases}
\end{equation}
\end{proposition}

Clearly the differential $d^1_{k,l} \: : \: E^1_{k,l} \rightarrow
E^1_{k,l-1}$ is induced from $\partial$.\\ \\
{\bf Important Digression : }
To digress for a moment,
let $\widehat{M}$ be a smooth manifold on which $\Gamma$ acts
freely, properly discontinuously and cocompactly.  
Put $M \: = \: \widehat{M}/\Gamma$, a closed manifold.
We construct a connection
\begin{equation} \label{1.35}
\nabla^{can} : C^\infty_c(\widehat{M}) \rightarrow 
\Omega^1(\widehat{M}, \C \Gamma) 
\otimes_{C^\infty_c(\widehat{M}) \rtimes \Gamma}
C^\infty_c(\widehat{M}).
\end{equation}
Let us note that
\begin{equation} \label{1.36}
\Omega^1(\widehat{M}, \C \Gamma) 
\otimes_{C^\infty_c(\widehat{M}) \rtimes \Gamma}
C^\infty_c(\widehat{M}) \: = \: 
\left( \Omega^1_c(\widehat{M}) \otimes_{C^\infty_c(\widehat{M})}
C^\infty_c(\widehat{M}) \right) \oplus
\left( \Omega^1(\C \Gamma) \otimes_{\C \Gamma}
C^\infty_c(\widehat{M}) \right)
\end{equation}
is isomorphic to 
$\Omega^1_c(\widehat{M}) \oplus
\left( \Omega^1(\C \Gamma) \otimes_{\C \Gamma}
C^\infty_c(\widehat{M}) \right)$.

\begin{lemma} \label{Lemma 5}
\cite[Prop. 9]{Lott (1992)}
Let $h \in C^\infty_c(\widehat{M})$ satisfy 
$\sum_{\gamma \in \Gamma} \gamma \cdot h \: = \: 1$. Define $\nabla^{can}$ by
\begin{equation} \label{1.37}
\nabla^{can} f \: = \: d^{\widehat{M}} f \: \oplus \: 
\sum_{\gamma \in \Gamma} d\gamma \otimes h (\gamma^{-1} \cdot f)
\end{equation}
for $f \in C^\infty_c(\widehat{M})$.
Then $\nabla^{can}$ is a connection on $C^\infty_c(\widehat{M})$.
\end{lemma}

One sees that 
\begin{equation} \label{1.38}
\left( \nabla^{can} \right)^2
\in \Hom_{C^\infty_c(\widehat{M}) \rtimes \Gamma} 
\left(C^\infty_c(\widehat{M}), 
\Omega^2(\widehat{M}, \C \Gamma) 
\otimes_{C^\infty_c(\widehat{M}) \rtimes \Gamma}
C^\infty_c(\widehat{M}) \right)
\end{equation} acts on 
$C^\infty_c(\widehat{M})$ as left multiplication by a $2$-form
$\Theta$ which commutes with
$C^\infty_c(\widehat{M}) \rtimes \Gamma$.
Explicitly, 
\begin{align} \label{1.39}
\Theta \: & = \: \sum_{\gamma \in \Gamma} d^{\widehat{M}} (\gamma \cdot h) 
\: d\gamma \: \gamma^{-1}  \:
- \: \sum_{\gamma, \gamma^\prime \in \Gamma} (\gamma \gamma^\prime \cdot h) \:
(\gamma \cdot h) \: d\gamma \: d\gamma^\prime \: (\gamma \gamma^\prime)^{-1} 
\\
&  = \: - \: \sum_{\gamma \in \Gamma} d^{\widehat{M}} (\gamma^{-1} \cdot h) 
\: \gamma^{-1}  \: d\gamma \:
- \: \sum_{\gamma, \gamma^\prime \in \Gamma} 
((\gamma \gamma^\prime)^{-1} \cdot h) \:
((\gamma^\prime)^{-1} \cdot h) \: (\gamma \gamma^\prime)^{-1}
\: d\gamma \: d\gamma^\prime. \notag
\end{align}
Note that if $\Gamma$ is infinite then $\Theta$ does not lie in 
$\Omega^2(\widehat{M}, \C \Gamma)$, as the sums involved are infinite.
Nevertheless, it sends
$C^\infty_c(\widehat{M})$ to 
$\Omega^2(\widehat{M}, \C \Gamma) 
\otimes_{C^\infty_c(\widehat{M}) \rtimes \Gamma}
C^\infty_c(\widehat{M})$.
Put
\begin{equation} \label{1.40}
\ch(\nabla^{can}) \: = \: e^{- \: 
\frac{ \Theta}{2 \pi i}} \in \End_{\Omega^*(\widehat{M}, \C \Gamma)}
\left( \Omega^*(\widehat{M}, \C \Gamma) 
\otimes_{C^\infty_c(\widehat{M}) \rtimes \Gamma} C^\infty_c(\widehat{M})
\right).
\end{equation}
Then the abelianization of $\ch(\nabla^{can})$ is closed. 
This can be seen by writing
$\Theta \: = \: dA \: - \: A^2$, where
\begin{equation} \label{1.41}
A \: = \: \sum_{\gamma \in \Gamma} d\gamma \: h \: \gamma^{-1} \: = \: - \:  
\sum_{\gamma \in \Gamma} (\gamma^{-1} \cdot h) \: \gamma^{-1} \: d\gamma.
\end{equation}
Then $d\Theta \: = \: - \: [\Theta, A]$ and
$d \ch \left(\nabla^{can} \right) \: = \: - \: \left[
\ch \left(\nabla^{can} \right), A \right]$.
Also, the cohomology class of $\ch(\nabla^{can})$ is
independent of the choice of $h$.
 
In the construction of $\ch \left(\nabla^{can} \right)$, we can allow
$h$ to be a Lipschitz function on $\widehat{M}$ 
(see \cite[Lemma 4]{Lott (1992)}, where $\ch \left(\nabla^{can} \right)$ is
called $\widetilde{\omega}_h$).  Let $E\Gamma$ be the 
bar simplicial complex
for $\Gamma$ (with degenerate simplices collapsed 
\cite[Chapter 1.5, Exercise 3b]{Brown (1982)}). We formally replace 
$\widehat{M}$ by $E \Gamma$. There is a
complex $\Omega^*(E\Gamma)$ of $\C$-valued polynomial forms on $E\Gamma$
defined as in 
\cite[p. 297]{Sullivan (1977)}. Let $j \in C(E\Gamma)$ be the 
barycentric coordinate corresponding to the vertex $e  \in E\Gamma$
\cite[(94)]{Lott (1992)}. Then
$\sum_{\gamma \in \Gamma} \gamma \cdot j \: = \: 1$
\cite[Lemma 5]{Lott (1992)}. (The support of $j$ may not be compact, but
this will not be a problem.) Define $\nabla^{univ}$ as in
(\ref{1.37}), replacing $h$ by $j$.
 
Let $\ch \left( \nabla^{univ} \right) \: = \:  e^{- \: 
\frac{\Theta}{2 \pi i}}$ 
denote the explicit form constructed using $j$. Then
$\ch \left( \nabla^{univ} \right)$ lies in
\begin{equation} \label{1.42}
\prod_{k=0}^\infty 
\bigoplus_{l=0}^\infty
\Omega^{k-l}(E\Gamma) \: \widehat{\otimes} \: \Omega^{k+l}(\C \Gamma).
\end{equation} 
Looking at the formula for $\ch \left( \nabla^{univ} \right)$, 
we see that in fact it lies in
\begin{equation} \label{1.43}
\prod_{k=0}^\infty \left[
\left( Z^{k}(E\Gamma) \: \widehat{\otimes}
\: \Omega^k(\C \Gamma) \right) \oplus 
\bigoplus_{l=1}^\infty \left(
\Omega^{k-l}(E\Gamma) \: 
\widehat{\otimes} \: \Omega^{k+l}(\C \Gamma) \right) \right].
\end{equation} 
Let $\eta \in \bigoplus_{a + b = n} C^a(\Gamma) \otimes C_b(B)$ be a 
graded $n$-trace on $\Omega^*(B, \C \Gamma)$.
We can pair
$\ch \left( \nabla^{univ} \right)$ and $\eta$ with respect to $\C \Gamma$,
to obtain 
\begin{equation} \label{1.44}
\langle \ch \left( \nabla^{univ} \right), \eta \rangle 
\in \bigoplus_{k=0}^\infty \left[
\left( Z^{k}(E\Gamma) \: \otimes \: C_{n-k}(B) \right) \oplus 
\bigoplus_{l=1}^\infty \left(
\Omega^{k-l}(E\Gamma) \: \otimes \: C_{n-k-l}(B) \right) \right].
\end{equation} 
{\bf End of Important Digression}
\begin{lemma} \label{Lemma 6}
The element constructed in (\ref{1.44}) is $\Gamma$-invariant.
\end{lemma}
\begin{pf}
With respect to (\ref{1.42}), let us write
$\ch \left( \nabla^{univ} \right)$ in the form
$\ch \left( \nabla^{univ} \right) \: = \: \sum_i \omega_i \: \widehat{\otimes}
\: \omega_i^\prime$, with
$\omega_i \in \Omega^*(E\Gamma)$ and $\omega_i^\prime \in
\Omega^*(\C \Gamma)$. 
As $\Theta$ commutes with $\Gamma$, it follows that for all 
$\gamma \in \Gamma$,
\begin{equation} \label{1.45}
\sum_i \omega_i \: \widehat{\otimes} \: \omega_i^\prime \: = \:
\sum_i (\gamma \cdot \omega_i) \: \widehat{\otimes} \: \gamma \: 
\omega_i^\prime \:
\gamma^{-1}. 
\end{equation}
Given $\phi \in \Omega^*_c(B)$, we can define the pairing
$\langle \langle \ch \left( \nabla^{univ} \right), \eta \rangle, \phi \rangle
\in \Omega^*(E \Gamma)$.
Then for any $\gamma \in \Gamma$,
\begin{align} \label{1.46}
\langle \langle \ch \left( \nabla^{univ} \right), \eta \rangle, \gamma \cdot
\phi \rangle \: & = \: \sum_i \omega_i \: \eta(\gamma \cdot \phi \:
\widehat{\otimes} \: \omega_i^\prime) \notag \\
& = \: \sum_i \gamma \cdot \omega_i \: \eta(\gamma \cdot \phi \:
\widehat{\otimes} \: \gamma \: \omega_i^\prime \gamma^{-1}) \notag \\
& = \: \sum_i \gamma \cdot \omega_i \: \eta(\gamma \: (\phi \: 
\widehat{\otimes} \: \omega_i^\prime) \: \gamma^{-1}) \notag \\
& = \: \sum_i \gamma \cdot \omega_i \: \eta \left( [\gamma, \: (\phi \: 
\widehat{\otimes} \: \omega_i^\prime) \: \gamma^{-1}] \: + \:
\phi \: \widehat{\otimes} \: \omega^\prime_i \right) \notag \\
& = \: \sum_i \gamma \cdot \omega_i \: \eta \left(
\phi \: \widehat{\otimes} \: \omega^\prime_i \right) \notag \\
& = \: \gamma \cdot 
\langle \langle \ch \left( \nabla^{univ} \right), \eta \rangle, 
\phi \rangle.
\end{align}
This proves the lemma.
\end{pf}

Equivalently, define a complex ${\cal C}_*(\Gamma)$ by
\begin{equation} \label{1.47}
{\cal C}_k(\Gamma) \: = \: Z^k(E \Gamma) \: \oplus \: \bigoplus_{l=1}^\infty
\Omega^{k-2l}(E\Gamma),
\end{equation}
with the natural chain map of degree $-1$. Then
\begin{equation} \label{1.48}
\langle \ch \left( \nabla^{univ} \right), \eta \rangle 
\in \bigoplus_{a+b=n}({\cal C}_a(\Gamma) \otimes C_b(B))^\Gamma.
\end{equation}
To take into account the normalization of $\Omega^*(\C \Gamma)$, 
let $E\Gamma_{(0)}$ denote the vertices of $E \Gamma$ and put
\begin{equation} \label{1.47.5}
\overline{\cal C}_k(\Gamma) \: = \: \Ker \left(
\left( Z^k(E \Gamma) \: \oplus \: \bigoplus_{l=1}^\infty
\Omega^{k-2l}(E\Gamma) \right)
 \rightarrow
\left( Z^k(E \Gamma_{(0)}) \: \oplus \: \bigoplus_{l=1}^\infty
\Omega^{k-2l}(E\Gamma_{(0)}) \right) \right).
\end{equation}
Consider the complex 
$\C \oplus \overline {\cal C}_*(\Gamma)$, where the factor $\C$ is in
degree zero. Then pairing with $\ch \left( \nabla^{univ} \right)$ gives
 a linear map from 
the graded traces on $\Omega^*(B, \C \Gamma)$ to the total
space of the double complex
$\left( (\C \oplus \overline{\cal C}_*(\Gamma)) \otimes C_*(B) 
\right)^\Gamma$.
\begin{lemma} \label{Lemma 7}
Pairing with $\ch \left( \nabla^{univ} \right)$ gives
a morphism from 
the complex of graded traces on $\Omega^*(B, \C \Gamma)$ to the total
complex of the double complex
$\left( (\C \oplus \overline{\cal C}_*(\Gamma)) \otimes C_*(B) 
\right)^\Gamma$.
\end{lemma}  
\begin{pf}
Given $\phi \in \Omega^*_c(B)$, we have an equality in
$\Omega^*(E \Gamma)$ :
\begin{align} \label{1.49}
& d \langle 
\langle \ch \left( \nabla^{univ} \right), \eta \rangle , \phi \rangle 
\: - \: \langle
\langle \ch \left( \nabla^{univ} \right), \eta \rangle , d\phi \rangle
\: - \: \langle
\langle \ch \left( \nabla^{univ} \right), d^t \eta \rangle , \phi \rangle
\:  = \notag \\
& \langle
\langle d\ch \left( \nabla^{univ} \right), \eta \rangle , \phi \rangle \: 
 = \: \langle
\langle \left[\ch \left( \nabla^{univ} \right), A \right], 
\eta \rangle , \phi \rangle. 
\end{align}
Let us write $\ch \left( \nabla^{univ} \right) \: = \:
\sum_i \omega_i \: \widehat{\otimes} \: \omega^\prime_i$ and
$A \: = \:
\sum_j a_j \: \widehat{\otimes} \: a^\prime_j$, with
$\omega_i, \: a_j \in \Omega^*(E\Gamma)$ and 
$\omega_i^\prime, \: a_j^\prime \in \Omega^*(\C \Gamma)$. 
Note that $\ch \left( \nabla^{univ} \right)$ and $A$ are concentrated
at the identity conjugacy class of $\Gamma$, so we can assume the same
about $\omega_i^\prime$ and $a_j^\prime$. Then
\begin{equation} \label{1.50}
\left[\ch \left( \nabla^{univ} \right), A \right] \: = \:
\sum_{i,j} [\omega_i \: \widehat{\otimes} \: \omega_i^\prime,
a_j \: \widehat{\otimes} \: a_j^\prime] \: = \: \pm \:
\sum_{i,j} \omega_i \: a_j \: \widehat{\otimes} \: [ \omega_i^\prime,
 a_j^\prime].
\end{equation}
Hence
\begin{equation} \label{1.51}
\langle \langle \left[\ch \left( \nabla^{univ} \right), A \right], 
\eta \rangle , \phi \rangle \: = \: \pm \: 
\sum_{i,j} \omega_i \: a_j \: \langle [ \omega_i^\prime,
\phi \: \widehat{\otimes} \: a_j^\prime], \eta \rangle,
\end{equation}
which vanishes as $\eta$ is a graded trace.  The lemma follows.
\end{pf}

We note that by the construction of $\overline{\cal C}_{*}$, 
the $E^1$-term $E^1_{k,l}$ of the double complex 
$\left( (\C \oplus \overline{\cal C}_*(\Gamma)) \otimes C_*(B) 
\right)^\Gamma$, i.e. the $k$-th homology of 
$\left( (\C \oplus \overline{\cal C}_*(\Gamma)) \otimes C_l(B) 
\right)^\Gamma$ with respect to the differential of
$\C \oplus \overline{\cal C}_*(\Gamma)$, is isomorphic to
(\ref{1.34}). 
\begin{lemma} \label{Lemma 8}
Pairing with $\ch \left( \nabla^{univ} \right)$ induces
an isomorphism from the $E^1$-term of the double complex ${\cal T}_{*,*}$
to the $E^1$-term of the double complex
$\left( (\C \oplus \overline{\cal C}_*(\Gamma)) \otimes C_*(B) 
\right)^\Gamma$.
\end{lemma}
\begin{pf}
For simplicity of notation, we only address the case when $k \: > \: 0$.
Consider first the component $\overline{\HH}^k(\Gamma; C_l(B))$ of the
$E^1_{k,l}$-term of the double complex ${\cal T}_{*,*}$. It follows from
\cite[Proposition 13]{Lott (1992)} that pairing with
$\ch \left( \nabla^{univ} \right)$ induces an isomorphism
from this component to the corresponding component of the
$E^1_{k,l}$-term of the double complex 
$\left( (\C \oplus \overline{\cal C}_*(\Gamma)) \otimes C_*(B) 
\right)^\Gamma$. 

Next, we remark that for both double complexes,
there is a (reduced) $S$-operator 
$\overline{S} \: : \: E^1_{k,l} \rightarrow E^1_{k+2,l}$ which sends
$\overline{\HH}^{k-2i}(\Gamma; C_l(B))$ to itself.  In the case of
the double complex ${\cal T}_{*,*}$, this $S$-operator is essentially
the dual of the one that acts on the right-hand-side of (\ref{1.3}). 
(A formula for the reduced $S$-operator, as opposed to the $S$-operator, 
is given in \cite[(55)]{Lott (1992a)}.) In the case of the double complex
$\left( \overline{\cal C}_*(\Gamma) \otimes C_*(B) 
\right)^\Gamma$, the $S$-operator is induced from the natural map
from $\left( \overline{\cal C}_k(\Gamma) \otimes C_l(B) 
\right)^\Gamma$ to
$\left( \overline{\cal C}_{k+2}(\Gamma) \otimes C_l(B) 
\right)^\Gamma$.
 
On the other hand, after adjusting the coefficients of its terms,
$\ch \left( \nabla^{univ} \right)$ is $\overline{S}$-invariant
\cite[Proposition 28]{Lott (1992a)}. Then using the known isomorphism
on the $\overline{\HH}^k(\Gamma; C_l(B))$ component of
$E^1_{k,l}$, along with the $\overline{S}$-operator, it follows that
pairing with $\ch \left( \nabla^{univ} \right)$ induces
an isomorphism on all of $E^1_{k,l}$.
\end{pf}

\begin{theorem} \label{Theorem 1}
The homology of $GT_{*,<e>}$, the complex of graded traces on
$\Omega^*(B, \C \Gamma)$ which are concentrated at the identity
conjugacy class, is isomorphic to the homology of the total complex of
the double complex 
$\left( (\C \oplus \overline{\cal C}_*(\Gamma)) \otimes C_*(B) 
\right)^\Gamma$
\end{theorem}
\begin{pf}
By Proposition \ref{Proposition 1}, the homology of $GT_{*,<e>}$ is
isomorphic to the homology of the total complex of the double complex
${\cal T}_{*,*}$.
We have shown that pairing with $\ch \left( \nabla^{univ} \right)$ gives
a morphism of double complexes, which is an isomorphism between 
the $E^1$-terms.
The differentials on the $E^1$-terms are both induced by $\partial$. It
follows that there is an isomorphism between the $E^\infty$-terms.
\end{pf}

Using the isomorphism between the homology of $GT_{*,<e>}$ and the
homology of the double complex ${\cal T}_{*,*}$, we can periodize
with respect to $\overline{S}$ to define the periodic homology
$\HH_*^{per}(GT_{*,<e>})$, with $* \in \Z/2\Z$. 
Let $\HH^*_\tau((E\Gamma \times B)/\Gamma)$ denote
the cohomology of $(E\Gamma \times B)/\Gamma$, twisted by the orientation
bundle of $B$ \cite[Section II.7]{Connes (1994)}.
Let $\overline{\HH}^*_\tau((E\Gamma \times B)/\Gamma)$ denote
the cohomology relative to $(E\Gamma_{(0)} \times B)/\Gamma \: = \: B$.
\begin{corollary} \label{Corollary 1}
There is an isomorphism between
$\HH_*^{per}(GT_{*,<e>})$ and
$\overline{\HH}^{* + \dim(B) + 2\Z}_\tau((E\Gamma \times B)/\Gamma)$.
\end{corollary}
\begin{pf}
As homology commutes with direct limits, Theorem \ref{Theorem 1} implies that
$\HH_*^{per}(GT_{*,<e>})$ is isomorphic to the homology of the
$\Z/2\Z$-graded complex 
\begin{equation}
\Ker \left(
\left( \Omega^*(E\Gamma) \: \otimes \: C_*(B) \right)^\Gamma \rightarrow 
\left( \Omega^*(E\Gamma_{(0)}) \: \otimes \: C_*(B) \right)^\Gamma \right).
\end{equation}
Dualizing with respect to $B$, this is isomorphic to
the $\Z/2\Z$-graded complex 
\begin{equation} \label{1.52}
\Ker \left( \left(
\Omega^*(E\Gamma) \: \otimes \: \Omega^{* \: + \: \dim(B)}_\tau(B) 
\right)^\Gamma \rightarrow 
\left(
\Omega^*(E\Gamma_{(0)}) \: \otimes \: \Omega^{* \: + \: \dim(B)}_\tau(B) 
\right)^\Gamma \right),
\end{equation}
where $\Omega^*_\tau(B)$ 
consists of the differential forms on $B$ with distributional
coefficients and with value in the orientation bundle $o(TB)$. 
The homology of this complex is 
$\overline{\HH}^{* + \dim(B) + 2\Z}_\tau((E\Gamma \times B)/\Gamma)$.
\end{pf}
\noindent
{\bf Remark : } If $\Gamma \: = \: \{e\}$ then the homology of the 
graded traces on $\Omega^*_c(B)$ is the {\it homology} of the currents on  $B$.
If $B$ is a point then the homology of the graded traces on
$\Omega^*(\C \Gamma)$ is essentially the group {\it cohomology} of 
$\Gamma$. In order to put these together into one object,
we have used Poincar\'e duality to convert the homology
of $B$ into the twisted cohomology of $B$. In this way we write the
periodic homology of graded traces on $\Omega^*(B, \C \Gamma)$ in terms
of the twisted cohomology of $(E\Gamma \times B)/\Gamma$. However, this
uniform description only exists after periodizing, because of the
grading reversal in the Poincar\'e duality.  For the 
unperiodized homology of graded traces, we must use the setup of Theorem 
\ref{Theorem 1}. \\

The closed graded traces coming from $(\C \otimes C_*(B))^\Gamma$ are also
relevant; they correspond exactly to the homology of the
$\Gamma$-invariant currents on $B$. In general, forgetting about
reduced cohomology, we have constructed
a map which sends a closed graded $n$-trace $\eta$ on 
$\Omega^*(B, \C \Gamma)$ to an element
$\Phi_\eta \in {\HH}^{n + \dim(B) + 2\Z}_\tau((E\Gamma \times B)/\Gamma)$.

\section{Fiberwise Operators and Traces} \label{Section 3}

In this section we first consider smoothing operators on 
$\widehat{M}$ which act fiberwise, preserve compact support 
and commute with $\Gamma$.
We define a $C^\infty_c(B)$-valued trace
$\Tr_{<e>}$ on such operators.  We then make various extensions of
$\Tr_{<e>}$. First, we extend it to an $\Omega^*(B, \C \Gamma)_{ab}$-valued
trace on form-valued operators.  Next, we extend it to a supertrace on
operators on $\Z_2$-graded vector bundles.  Finally, we extend it to an
$\Omega^*(B, {\cal B}^\omega)_{ab}$-valued trace on smoothing operators
whose Schwartz kernels have sufficiently rapid decay.

Let $\widehat{M}$ be a smooth manifold on which $\Gamma$ acts freely,
properly discontinuously and cocompactly. Put $M \: = \: \widehat{M}/\Gamma$.
Let $B$ be a smooth manifold on which $\Gamma$ acts, not necessarily
freely or properly. Suppose 
that there is a $\Gamma$-invariant submersion $\pi \: : \: \widehat{M}
\rightarrow B$, a fiber of which we denote by $Z$. Then
$M$ is foliated by the images of $Z$ under the map $\widehat{M} \rightarrow M$.
That is, given $b \in B$,
put $Z_b \: = \: \pi^{-1}(b)$. Then the corresponding
leaf of the foliation ${\cal F}$ is $Z_b/\Gamma_b \subset \widehat{M}/\Gamma$.

Let $TZ$ denote the vertical
tangent bundle of $\widehat{M} \rightarrow B$, 
a vector bundle on $\widehat{M}$. Let $g^{TZ}$ be a
$\Gamma$-invariant Euclidean inner product 
on $TZ$. Give $Z_b$ the corresponding
Riemannian metric and induced metric space structure $d$.
As $\Gamma$ acts cocompactly
on $\widehat{M}$, preserving the submersion structure, it follows that
$\{Z_b\}_{b \in B}$ has bounded geometry. That is, 
there is a uniform upper bound on the
absolute values of the sectional curvatures, 
and a uniform lower bound on the injectivity radii.
Let $d\vol_Z$ denote the Riemannian volume forms on the fibers
$\{Z_b\}_{b \in B}$.

An element $K$ of $\End_{C^\infty_c(B) \rtimes \Gamma}
\left(
C^\infty_c(\widehat{M})
\right)$ has a Schwartz kernel $K(z,w)$, with respect to its
fiberwise action, so that we can write
\begin{equation} \label{2.1}
(KF)(z) \: = \: \int_{Z_{\pi(z)}} K(z,w) \: F(w) \: d\vol_{Z_{\pi(z)}}(w).
\end{equation}
for $F \in C^\infty_c(\widehat{M})$.

\begin{definition} $\End^\infty_{C^\infty_c(B) \rtimes \Gamma}
\left( 
C^\infty_c(\widehat{M})
\right)$ is the subalgebra of
$\End_{C^\infty_c(B) \rtimes \Gamma}
\left(
C^\infty_c(\widehat{M})
\right)$
consisting of elements $K$ 
with a smooth integral kernel in $C^\infty(\widehat{M} \times_B \widehat{M})$.
\end{definition}

Note that for each $b \in B$ and each
$w \in Z_b$, the function $K_w(z) \: = \: K(z, w)$ has compact support in $z$.
To simplify notation, if $K \in \End^\infty_{C^\infty_c(B) \rtimes \Gamma}
\left(
C^\infty_c(\widehat{M})
\right)$ then we will write the action of $K$ on 
$C^\infty_c(\widehat{M})$ by
\begin{equation} \label{2.2}
(KF)(z) \: = \: \int_Z K(z,w) \: F(w) \: d\vol_Z(w).
\end{equation}
That is, $\int_Z$ denotes fiberwise integration.
The convolution product on $\End^\infty_{C^\infty_c(B) \rtimes \Gamma}
\left(
C^\infty_c(\widehat{M})
\right)$ is given by
\begin{equation} \label{2.3}
(K K^\prime)(z,w) \: = \: \int_Z K(z,u) \: K(u, w) \: d\vol_Z(u).
\end{equation}
In this way $\End^\infty_{C^\infty_c(B) \rtimes \Gamma}
\left(
C^\infty_c(\widehat{M})
\right)$ is an algebra over $\C$, 
possibly without unit.

Let $\phi \in C^\infty_c(\widehat{M})$ satisfy
$\sum_{\gamma \in \Gamma} \gamma \cdot \phi \: = \: 1$.
Given $K \in \End^\infty_{C^\infty_c(B) \rtimes \Gamma}
\left(
C^\infty_c(\widehat{M})
\right)$ and $b \in B$, put
\begin{equation} \label{2.4}
\Tr(K)(b) \: = \: \sum_{\gamma \in \Gamma_b} \left( \int_{Z_b} 
\phi(w) \: K(w \gamma^{-1}, w) \: 
d\vol_{Z_b}(w) \right) \gamma.
\end{equation}
From the support condition on $K$, $\Tr(K)(b) \in \C \Gamma_b$. 

To express the range of $\Tr$ in a better way, 
let $(B^\gamma)^\C$ be the complex-valued functions on $B^\gamma$.
There is an inclusion $(B^\gamma)^\C \subset B^\C$ coming from extension
by zero.
Then $\bigoplus_{\gamma \in \Gamma} (B^\gamma)^\C \: \gamma$ is
an algebra, as a subalgebra of $B^\C \rtimes \Gamma$. Put
\begin{equation} \label{2.5}
\left( \bigoplus_{\gamma \in \Gamma} (B^\gamma)^\C \: \gamma \right)_{ab} 
\: = \:
\frac{\bigoplus_{\gamma \in \Gamma} (B^\gamma)^\C \:
\gamma}{[\bigoplus_{\gamma \in \Gamma} (B^\gamma)^\C \: \gamma,
\bigoplus_{\gamma \in \Gamma} (B^\gamma)^\C \: \gamma]}.
\end{equation}

Consider $\Tr$ from (\ref{2.4}).

\begin{proposition} \label{Proposition 3}
$\Tr : \End^\infty_{C^\infty_c(B) \rtimes \Gamma}
\left(
C^\infty_c(\widehat{M})
\right) \rightarrow 
\left( \bigoplus_{\gamma \in \Gamma} (B^\gamma)^\C \: \gamma \right)_{ab}$
is a trace.
\end{proposition}
\begin{pf}
Let $\{O_\alpha\}$ be the orbits of $\Gamma$ in $B$ and let
$b_\alpha \in O_\alpha$ be representative elements.
Put $\left( \C \Gamma_{b_\alpha} \right)_{ab} \: = \:
\frac{\C \Gamma_{b_\alpha}}{
[\C \Gamma_{b_\alpha}, \C \Gamma_{b_\alpha}]}$.  Then there
is an isomorphism
\begin{equation} \label{2.6}
I : \left( \bigoplus_{\gamma \in \Gamma} (B^\gamma)^\C \: \gamma \right)_{ab} 
\rightarrow \prod_\alpha \left( \C \Gamma_{b_\alpha} \right)_{ab}.
\end{equation}
Namely, given $\gamma \in \Gamma$ and 
$f \in \left( B^\gamma \right)^\C$,
\begin{equation} \label{2.7}
I([f \: \gamma]) \: = \: \prod_{\alpha} \sum_{b \in b_\alpha \Gamma}
f(b) \: [\gamma_b^\prime \gamma (\gamma_b^\prime)^{-1}],
\end{equation}
where $\gamma_b^\prime \in \Gamma$ is such that 
$b \: = \: b_\alpha \gamma_b^\prime$.
Under this isomorphism, we obtain 
\begin{equation} \label{2.8}
I(\Tr(K)) \: = \: \prod_{\alpha} \sum_{b \in b_\alpha \Gamma}
\sum_{\gamma \in \Gamma_b} \left( \int_{Z_b} 
\phi(w) \: K(w \gamma^{-1}, w) \: 
d\vol_{Z_b}(w) \right) [\gamma_b^\prime \gamma (\gamma_b^\prime)^{-1}].
\end{equation} 
Thus it is enough to show that for each $\alpha$, if we put
\begin{equation} \label{2.9}
I_\alpha(\Tr(K)) \: = \: \sum_{b \in b_\alpha \Gamma}
\sum_{\gamma \in \Gamma_b} \left( \int_{Z_b} 
\phi(w) \: K(w \gamma^{-1}, w) \: 
d\vol_{Z_b}(w) \right) [\gamma_b^\prime \gamma (\gamma_b^\prime)^{-1}]
\end{equation} 
then $I_\alpha(\Tr(K K^\prime)) \: = \: I_\alpha(\Tr(K^\prime K))$.

Let $\{\gamma_\beta\}_{\beta = 1}^\infty$ 
be a sequence of elements of $\Gamma$ such that
$\Gamma_{b_\alpha} \backslash \Gamma \: = \: \{
\Gamma_{b_\alpha} \gamma_\beta \}_{\beta = 1}^\infty$. Then
writing $b \: = \: b_\alpha \: \gamma_\beta$ and
$\gamma \: = \: \gamma_\beta^{-1} \: \gamma_\alpha \: \gamma_\beta$,
\begin{align} \label{2.10}
I_\alpha(\Tr(K)) \: & = \: \sum_{\beta = 1}^\infty
\sum_{\gamma_\alpha \in \Gamma_{b_\alpha}} \left( \int_{Z_{b_\alpha
\gamma_\beta}} 
\phi(w) \: K(w \gamma_\beta^{-1} \gamma_{\alpha}^{-1} \gamma_\beta, w) \: 
d\vol_{Z_{b_\alpha \gamma_\beta}}(w) \right) 
[\gamma_\alpha] \notag \\
& = \: \sum_{\beta = 1}^\infty
\sum_{\gamma_\alpha \in \Gamma_{b_\alpha}} \left( \int_{Z_{b_\alpha
\gamma_\beta}} 
\phi(w) \: K(w \gamma_\beta^{-1} \gamma_\alpha^{-1}, 
w \gamma_{\beta}^{-1}) \: 
d\vol_{Z_{b_\alpha \gamma_\beta}}(w) \right) 
[\gamma_\alpha] \notag \\
& = \: \sum_{\beta = 1}^\infty
\sum_{\gamma_\alpha \in \Gamma_{b_\alpha}} \left( \int_{Z_{b_\alpha}} 
\phi(w \gamma_\beta) \: 
K(w \gamma_\alpha^{-1}, w) \: 
d\vol_{Z_{b_\alpha}}(w) \right) 
[\gamma_\alpha].
\end{align}
Define $\phi_\alpha \in C^\infty_c(Z_{b_\alpha})$ by
$\phi_\alpha(w) \: = \: \sum_{\beta = 1}^\infty \phi(w \gamma_\beta)$.
Then $\sum_{\gamma_\alpha \in \Gamma_{b_\alpha}}
\gamma_\alpha \cdot \phi_\alpha \: = \: 1$ and
\begin{equation} \label{2.11}
I_\alpha(\Tr(K)) \: = \:
\sum_{\gamma_\alpha \in \Gamma_{b_\alpha}} \left( \int_{Z_{b_\alpha}} 
\phi_\alpha(w) \: 
K(w \gamma_\alpha^{-1}, w) \: 
d\vol_{Z_{b_\alpha}}(w) \right) 
[\gamma_\alpha]. 
\end{equation} 
It now follows from \cite[Prop. 7]{Lott (1992)} that $I_\alpha
\circ \Tr$ is a trace.
(The formula in \cite[Prop. 7]{Lott (1992)} is slightly different 
because \cite{Lott (1992)}  considers function spaces to be right 
$\Gamma$-modules instead of left $\Gamma$-modules.) This proves the
proposition.
\end{pf}

One can show that $\Tr$ is independent of the choice of $\phi$.

We can decompose $\Tr(K)$ according to the conjugacy classes
of $\Gamma$. In particular, the component corresponding to the 
conjugacy class of $e \in \Gamma$ is
\begin{equation} \label{2.12}
\Tr_{<e>}(K)(b) \: = \: \left( \int_{Z_b} 
\phi(w) \: K(w, w) \: 
d\vol_{Z_b}(w) \right).
\end{equation}
We see that $\Tr_{<e>}$ is a trace on 
$\End^\infty_{C^\infty_c(B) \rtimes \Gamma}
\left(
C^\infty_c(\widehat{M})
\right)$ 
which takes values in the
co-invariants $\left( C^\infty_c(B) \right)_\Gamma$.

We will need some slight extensions of $\Tr$. First,  
consider the $\Z$-graded algebra
\begin{equation} \label{2.13}
\Hom^\infty_{C^\infty_c(B) \rtimes \Gamma}
\left(
C^\infty_c(\widehat{M}), \:
\Omega^*(B, \C \Gamma)
\otimes_{C^\infty_c(B) \rtimes \Gamma} C^\infty_c(\widehat{M})
\right)
\end{equation}
consisting of elements of
$\Hom_{C^\infty_c(B) \rtimes \Gamma}
\left(
C^\infty_c(\widehat{M}), \:
\Omega^*(B, \C \Gamma)
\otimes_{C^\infty_c(B) \rtimes \Gamma} C^\infty_c(\widehat{M})
\right)$ 
with a smooth integral kernel.
An element $K$ of
$\Hom^\infty_{C^\infty_c(B) \rtimes \Gamma}
\left(
C^\infty_c(\widehat{M}), \:
( \Omega^k(B) \otimes \Omega^l(\C \Gamma) )
\otimes_{C^\infty_c(B) \rtimes \Gamma} C^\infty_c(\widehat{M})
\right)$
can be represented as a finite sum
\begin{equation} \label{2.14}
K \: = \: \sum_{g_1, \ldots, g_l \in \Gamma} dg_1 \ldots dg_l
\: K_{g_1, \ldots, g_l}, 
\end{equation}
where $K_{g_1, \ldots, g_l}$ has a smooth integral kernel
$K_{g_1, \ldots, g_l}(z, w) \in \Lambda^k(T^*_{\pi(z)} B)$.
It acts on
$C^\infty_c(\widehat{M})$ by
\begin{equation} \label{2.15}
(KF)(z) \: = \: \sum_{g_1, \ldots, g_l} 
dg_1 \ldots dg_l \: \left( \int_Z K_{g_1, \ldots, g_l}(z,w)
\: F(w) \: d\vol_Z(w) \right) .
\end{equation}
Then we define $\Tr$ to act on
$\Hom^\infty_{C^\infty_c(B) \rtimes \Gamma}
\left(
C^\infty_c(\widehat{M}), \:
\Omega^*(B, \C \Gamma)
\otimes_{C^\infty_c(B) \rtimes \Gamma} C^\infty_c(\widehat{M})
\right)$
by the formula
\begin{equation} \label{2.16}
\Tr(K)(b) \: = \: \sum_{g_0,g_1, \ldots, g_l \in \Gamma :
g_0 \cdots g_l \in \Gamma_b} 
(dg_1 \ldots dg_l) \: g_0 \: \left( \int_{Z_b} \phi(w) \:
K_{g_1, \ldots, g_l}
(w g_0^{-1}, w)
\: d\vol_{Z_b}(w) \right).
\end{equation}
(Compare \cite[(36)]{Lott (1992)}.)\\ \\
{\bf Example : } Suppose that $\widehat{M} \: = \: \Gamma$ and that $B \: =
\: \pt$.
The action of $\Gamma$ on $C^\infty_c(\widehat{M})$ is given by
$g \cdot \delta_h \: = \: \delta_{h g^{-1}}$. There is an isomorphism
of left $\C \Gamma$-modules $\C \Gamma \rightarrow C^\infty_c(\widehat{M})$
which sends $h$ to $\delta_{h^{-1}}$. 
Consider $h_0 \: dh_1 \in \Omega^1(\C \Gamma)$
and the corresponding element $K \: \in \Hom_{\C \Gamma}(\C \Gamma,
\Omega^*(\C \Gamma) \otimes_{\C \Gamma} \C \Gamma)$  given by
\begin{equation} \label{2.17}
K(h) \: = \: h (h_0 \: dh_1) \: = \: d(h h_0 h_1) \: e \: - \:
d(h h_0) \: h_1.
\end{equation}
Then under the above isomorphism, $K \in
\Hom_{\C \Gamma}(C^\infty_c(\widehat{M}),
\Omega^*(\C \Gamma) \otimes_{\C \Gamma} C^\infty_c(\widehat{M}))$ acts
by 
\begin{equation} \label{2.18}
K(\delta_{h^{-1}}) \: = \: d(h h_0 h_1) \: \delta_e \: - \:
d(h h_0) \: \delta_{h_1^{-1}}.
\end{equation}
Thus 
\begin{equation} \label{2.19}
K_{g_1}(z, w) \: = \: \delta_{w^{-1} h_0 h_1, g_1} \: \delta_{z,e}
 \: - \: \delta_{w^{-1} h_0, g_1} \: \delta_{z, h_1^{-1}}
\end{equation}
and
\begin{equation} \label{2.20}
\sum_{g_0, g_1} dg_1 \: g_0 \: K_{g_1}(w g_0^{-1}, w) \: = \:
d(w^{-1} h_0 h_1) \: w \: - \: d(w^{-1} h_0) \: h_1 w \: = \:
w^{-1} (h_0 \: dh_1) w. 
\end{equation}
If $\phi \: = \: \delta_x$ then we get
\begin{equation} \label{2.21}
\Tr(K)(\pt.) \: = \: x^{-1} (h_0 \: dh_1) x,
\end{equation}
which is equivalent to 
$h_0 \: dh_1$ in $\Omega^*(\C \Gamma)_{ab}$, as it should be.\\
{\bf End of Example} \\
 
As before, we can decompose $\Tr(K)$ according to the conjugacy classes of
$\Gamma$. In this paper we will only be concerned with the component of
$\Tr(K)$ associated to the conjugacy class of $e \in \Gamma$, namely
\begin{equation} \label{2.22}
\Tr_{<e>}(K)(b) \: = \: \sum_{g_0,g_1, \ldots, g_l \in \Gamma : 
g_0 \cdots g_l \: = \: e} (dg_1 \ldots dg_l) \: g_0
\left( \int_{Z_b} \phi(w) \:
K_{g_1, \ldots, g_l}
(w g_0^{-1}, w)
\: d\vol_{Z_b}(w) \right) .
\end{equation}
Then one sees from (\ref{2.22}) that $\Tr_{<e>}$ is a trace on
\begin{equation} \label{2.23}
\Hom^\infty_{C^\infty_c(B) \rtimes \Gamma}
\left(
C^\infty_c(\widehat{M}), \:
\Omega^*(B, \C \Gamma)
\otimes_{C^\infty_c(B) \rtimes \Gamma} C^\infty_c(\widehat{M})
\right)
\end{equation}
which takes values in $\Omega^*(B, \C \Gamma)_{ab}$.

Next, let $\widehat{E}$ be a $\Gamma$-equivariant $\Z_2$-graded  Hermitian
vector bundle on 
$\widehat{M}$. Define 
\begin{equation} \label{2.24}
\End^\infty_{C^\infty_c(B) \rtimes \Gamma}
\left(
C^\infty_c(\widehat{M}; \widehat{E})
\right)
\end{equation}
as before, except with
$K(z, w) \in \Hom(\widehat{E}_w, \widehat{E}_z)$. Also define
\begin{equation} \label{2.25}
\Hom^\infty_{C^\infty_c(B) \rtimes \Gamma}
\left(
C^\infty_c(\widehat{M};
\widehat{E}), \:
\Omega^*(B, \C \Gamma)
\otimes_{C^\infty_c(B) \rtimes \Gamma} C^\infty_c(\widehat{M};
\widehat{E})
\right)
\end{equation}
as before, except with
$K_{g_1, \ldots, g_l}(z, w) \in \Lambda^k(T^*_{\pi(z)} B) \: \otimes \:
\Hom(\widehat{E}_w, \widehat{E}_z)$.
Then we obtain a supertrace
\begin{equation} \label{2.26}
\Tr_{s,<e>} \: : \:
\Hom^\infty_{C^\infty_c(B) \rtimes \Gamma}
\left(
C^\infty_c(\widehat{M};
\widehat{E}), \:
\Omega^*(B, \C \Gamma)
\otimes_{C^\infty_c(B) \rtimes \Gamma} C^\infty_c(\widehat{M};
\widehat{E})
\right)
\rightarrow \Omega^*(B, \C \Gamma)_{ab}
\end{equation} 
by
\begin{equation} \label{2.27}
\Tr_{s,<e>}(K)(b) \: = \: \sum_{g_0,g_1, \ldots, g_l \in \Gamma : 
g_0 \cdots g_l \: = \: e}
(dg_1 \ldots dg_l) \: g_0  
\left( \int_{Z_b} \phi(w) \:
\tr_s (K_{g_1, \ldots, g_l}
(w g_0^{-1}, w))
\: d\vol_{Z_b}(w) \right).
\end{equation}

Finally,
choose a finite generating set for $\Gamma$ and consider the corresponding
right-invariant word metric $\parallel \cdot \parallel$. 
Let ${\cal B}^\omega$ be the formal sums
$\sum_{\gamma \in \Gamma} c_\gamma \: \gamma$ such that
$|c_\gamma|$ decreases faster than any exponential in $\parallel \gamma 
\parallel$
(see \cite[Lemma 2]{Lott (1992)}). Then ${\cal B}^\omega$ is a locally
convex Fr\'echet algebra \cite[Prop. 4]{Lott (1992)}.
The notation ``$\omega$'' comes from the fact that if $\Gamma \: = \: \Z$
then ${\cal B}^\omega$ can be identified with the holomorphic functions on
$\C - 0$.

Put 
\begin{equation} \label{2.28}
C^\infty(B, {\cal B}^\omega) \: = \: {\cal B}^\omega
\otimes_{\C \Gamma} ( C^\infty_c(B) \rtimes \Gamma),
\end{equation}
i.e. the quotient of the locally convex topological vector space
${\cal B}^\omega \otimes ( C^\infty_c(B) \rtimes \Gamma)$ by 
$\overline{\spann 
\{ (a, f) - (\gamma a, \gamma \cdot f) 
\}}$, where $a \in {\cal B}^\omega$, $f \in C^\infty_c(B) \rtimes \Gamma$ and
$\gamma \in \Gamma$.
Then $C^\infty(B, {\cal B}^\omega)$ is a locally convex topological
algebra which has 
$C^\infty_c(B) \: \rtimes \: \Gamma$ as a dense subalgebra.
We can write an element of $C^\infty(B, {\cal B}^\omega)$ as a infinite sum
$\sum_{\gamma \in \Gamma} f_\gamma \: \gamma$, where the functions on $B$
$\{ \gamma^{-1} \cdot f_\gamma \}_{\gamma \in \Gamma}$ all have support
in some compact set $K \subset B$ and have the decay property that for all
$q \in \Z^+$, 
\begin{equation} \label{2.29}
\sup_{b \in K, \: \gamma \in \Gamma} 
\: e^{q \: \parallel \gamma \parallel} \: |f_{\gamma}(b 
\gamma^{-1})| \: < \: \infty,
\end{equation}
along with the analogous statement for the derivatives of $f$.

Similarly, we define a locally convex GDA
$\Omega^*(B, {\cal B}^\omega)$ by saying that an element of type $(k,l)$
is an infinite sum
$\sum_{\gamma_0, \ldots, \gamma_l \in \Gamma} 
\omega_{\gamma_0, \ldots, \gamma_l} \: \gamma_0 \: d\gamma_1 \ldots
d\gamma_l$,
where the $k$-forms on $B$
$\{ (\gamma_0 \ldots \gamma_l)^{-1} \cdot \omega_{\gamma_0, \ldots, \gamma_l}
\}_{\gamma_0, \ldots, \gamma_l \in \Gamma}$ all have support
in some compact set $K \subset B$ and have the decay property that for all
$q \in \Z^+$, 
\begin{equation} \label{2.30}
\sup_{b \in K, \: \gamma_0, \ldots, \gamma_l \in \Gamma} 
\: e^{q \: \parallel \gamma_0 \ldots \gamma_l \parallel} \: 
|\omega_{\gamma_0, \ldots, \gamma_l}(b (\gamma_0 \ldots \gamma_l)^{-1})| \: 
< \: \infty,
\end{equation}
along with the analogous statement for the derivatives.
Then $\Omega^*(B, {\cal B}^\omega)$
has $\Omega^*(B, \C \Gamma)$ as a dense subalgebra.

Put 
\begin{equation} \label{2.31}
C^\infty_{{\cal B}^\omega}(\widehat{M}) \: = \:
{\cal B}^\omega \otimes_{\C \Gamma} 
C^\infty_c(\widehat{M}).
\end{equation}
Then $C^\infty_{{\cal B}^\omega}(\widehat{M})$ is a left 
$C^\infty(B, {\cal B}^\omega)$-module. As in 
\cite[Prop. 5]{Lott (1992)} and using the cocompactness of the
$\Gamma$-action on $\widehat{M}$, the elements of 
$C^\infty_{{\cal B}^\omega}(\widehat{M})$ can be characterized as the
elements $F \in C^\infty(\widehat{M})$ such that 
for any $b \in B$, $m_0 \in Z_b$ and $q \in \Z^+$, 
\begin{equation} \label{2.32}
\sup_{z \in Z_b} \: e^{q \: d(z,m_0)} \: |F(z)| \: < \: \infty,
\end{equation}
along with the analogous property for the covariant derivatives of $f$. 
Let 
$\End^\infty_{C^\infty(B, {\cal B}^\omega)}
\left(
C^\infty_{{\cal B}^\omega}(\widehat{M})
\right)$ be the subalgebra of
$\End_{C^\infty(B, {\cal B}^\omega)} \left( 
C^\infty_{{\cal B}^\omega}(\widehat{M})\right)$ consisting of elements
$K$ with
a smooth integral kernel $K(z,w)$.
Then as in \cite[Prop. 6]{Lott (1992)},
the elements of $\End^\infty_{C^\infty(B, {\cal B}^\omega)}
\left(
C^\infty_{{\cal B}^\omega}(\widehat{M})
\right)$ can be characterized as the 
$\Gamma$-invariant elements $K(z, w) \in C^\infty(
\widehat{M} \times_B \widehat{M})$ such for any $b \in B$ and $q \in \Z^+$, 
\begin{equation} \label{2.33}
\sup_{z,w \in Z_b} \: e^{q \: d(z,w)} \: |K(z,w)| \: < \: \infty,
\end{equation}
along with the analogous property for the covariant derivatives of $K$.
The convolution product in 
$\End^\infty_{C^\infty(B, {\cal B}^\omega)}
\left(
C^\infty_{{\cal B}^\omega}(\widehat{M})
\right)$ is given by
the same expression as (\ref{2.3}), and makes sense because of the bounded
geometry of $\{Z_b\}_{b \in B}$. 
With the natural definition of
\begin{equation} \label{2.34}
\Hom^\infty_{C^\infty(B, {\cal B}^\omega)}
\left(
C^\infty_{{\cal B}^\omega}
(\widehat{M}), \:
\Omega^*(B, {\cal B}^\omega)
\otimes_{C^\infty(B, {\cal B}^\omega)} 
C^\infty_{{\cal B}^\omega}(\widehat{M})
\right),
\end{equation}
an element $K$ can be written as a formal
$\Gamma$-invariant sum (\ref{2.14}). In particular,
the individual terms have 
the decay property that
for any $b \in B$ and $q \in \Z^+$, 
\begin{equation} \label{2.35}
\sup_{w \in Z_b, z \in Z_{b g_l \cdots g_1}} 
e^{q \: d(z,wg_l \cdots g_1)} \: K_{g_1, \ldots, g_l}(z,w) \: < \: \infty.
\end{equation}
The formula (\ref{2.12}) extends to a trace
$\Tr_{<e>} \: : \: 
\End^\infty_{C^\infty(B, {\cal B}^\omega)}
\left(
C^\infty_{{\cal B}^\omega}(\widehat{M})
\right) \rightarrow 
\frac{C^\infty(B, 
{\cal B}^\omega)}{[C^\infty(B, 
{\cal B}^\omega), C^\infty(B, {\cal B}^\omega)]}$
which is concentrated at the identity conjugacy class.
The formula (\ref{2.22}) extends to a trace
\begin{equation} \label{2.37}
\Tr_{<e>} \: : \: 
\Hom^\infty_{C^\infty(B, {\cal B}^\omega)}
\left(
C^\infty_{{\cal B}^\omega}(\widehat{M}), \:
\Omega^*(B, {\cal B}^\omega)
\otimes_{C^\infty(B, {\cal B}^\omega)} 
C^\infty_{{\cal B}^\omega}(\widehat{M})
\right)
\rightarrow \Omega^*(B, {\cal B}^\omega)_{ab}.
\end{equation}

If $\widehat{E}$ is a $\Z_2$-graded 
$\Gamma$-invariant Hermitian vector bundle on $\widehat{M}$,
with an invariant Hermitian connection, then we can define
$\End^\infty_{C^\infty(B, {\cal B}^\omega)}
\left(
C^\infty_{{\cal B}^\omega}(\widehat{M}; \widehat{E})
\right)$
and a supertrace
\begin{equation} \label{2.38}
\Tr_{s,<e>} \: : \:
\Hom^\infty_{C^\infty(B, {\cal B}^\omega)}
\left(
C^\infty_{{\cal B}^\omega}(\widehat{M};
\widehat{E}), \:
\Omega^*(B, {\cal B}^\omega)
\otimes_{C^\infty(B, {\cal B}^\omega)} 
C^\infty_{{\cal B}^\omega}(\widehat{M};
\widehat{E})
\right)
\rightarrow \Omega^*(B, {\cal B}^\omega)_{ab}.
\end{equation}

\section{Superconnections and Small-$s$ Asymptotics} \label{Section 4}

In this section we define the superconnection $A_s$ and compute the
small-$s$ limit of the supertrace of $e^{- \: A_s^2}$, thereby obtaining
the right-hand-side of the index theorem.

Let $\pi : \widehat{M} \rightarrow B$ be a $\Gamma$-invariant submersion
as in the previous section.  We choose
a $\Gamma$-invariant vertical Riemannian metric $g^{TZ}$ on $TZ$ and a 
$\Gamma$-invariant horizontal
distribution $T^H \widehat{M}$ on $\widehat{M}$.

Suppose that $Z$ is even-dimensional. Let $\widehat{E}$ be a
$\Gamma$-invariant Clifford
bundle on $\widehat{M}$ which is equipped 
with a $\Gamma$-invariant connection. For simplicity of notation, 
we asssume that
$\widehat{E} \: = \: S^Z \: \widehat{\otimes} \: \widehat{V}$, 
where $S^Z$ is a vertical spinor bundle
and $\widehat{V}$ 
is an auxiliary vector bundle on $\widehat{M}$. More precisely,
suppose that the vertical tangent
bundle $TZ$ has a spin structure.  Let $S^Z$ be the 
vertical spinor bundle, a $\Gamma$-invariant $\Z_2$-graded Hermitian
vector bundle on $\widehat{M}$. Let
$\widehat{V}$ be another  $\Gamma$-invariant $\Z_2$-graded Hermitian
vector bundle on $\widehat{M}$ which is equipped 
with a $\Gamma$-invariant Hermitian connection $\nabla^{\widehat{V}}$. 
Then we put $\widehat{E} \: = \: S^Z \:
\widehat{\otimes} \: \widehat{V}$. The case of general 
$\Gamma$-invariant Clifford bundles $\widehat{E}$ can
be treated in a way completely analogous to what follows.

Let $Q$ denote the vertical Dirac-type operator acting on
$C^\infty_c(\widehat{M}; \widehat{E})$.
From finite-propagation-speed estimates
as in \cite[Pf. of Prop 8]{Lott (1992)}, 
along with
the bounded geometry
of $\{Z_b\}_{b \in B}$, 
for any $s \: > \: 0$ we have 
\begin{equation} \label{3.1}
e^{- \: s^2 \: Q^2} \in 
\End^\infty_{C^\infty(B, {\cal B}^\omega)} 
\left(
C^\infty_{{\cal B}^\omega}(\widehat{M}; \widehat{E})
\right).
\end{equation}

Let
\begin{equation} \label{3.2}
A^{Bismut}_s \: : \: C^\infty_c(\widehat{M}; \widehat{E}) 
\rightarrow
\Omega_c^*(B) \otimes_{C^\infty_c(B)} C^\infty_c(\widehat{M}; \widehat{E})
\end{equation}
denote the Bismut superconnection 
\cite[Proposition 10.15]{Berline-Getzler-Vergne (1992)}.
In the cited reference it is constructed for fiber bundles with compact
fibers.  However, being a differential operator, it makes sense for any
submersion. 
It is of the form
\begin{equation} \label{3.3}
A^{Bismut}_s \: = \: s \: Q \: + \: \nabla^u \: - \: \frac{1}{4s} \: c(T),
\end{equation}
where $\nabla^u$ is a certain Hermitian connection and
$c(T)$ is Clifford multiplication by the curvature $2$-form $T$ of
the horizontal distribution $T^H \widehat{M}$.
We also denote by
\begin{equation} \label{3.4}
A^{Bismut}_s \: : \: C^\infty_{{\cal B}^\omega}(\widehat{M}; \widehat{E}) 
\rightarrow
({\cal B}^\omega \: \otimes_{\C \Gamma} \: \Omega_c^*(B)) 
\otimes_{({\cal B}^\omega \: \otimes_{\C \Gamma} \: C^\infty_c(B))} 
C^\infty_{{\cal B}^\omega}(\widehat{M}; \widehat{E})
\end{equation}
its extension to 
$C^\infty_{{\cal B}^\omega}(\widehat{M}; \widehat{E})$.
One can use finite-propagation-speed estimates, along with
the bounded geometry
of $\{Z_b\}_{b \in B}$ and the Duhamel expansion as in
\cite[Theorem 9.48]{Berline-Getzler-Vergne (1992)},
to show that we obtain a well-defined element
\begin{equation} \label{3.5}
e^{- \: (A^{Bismut}_s)^2} \in
\Hom^\infty_{({\cal B}^\omega \: \otimes_{\C \Gamma} \: C^\infty_c(B))} 
\left(
C^\infty_{{\cal B}^\omega}(\widehat{M}; \widehat{E}), \:
({\cal B}^\omega \: \otimes_{\C \Gamma} \: \Omega_c^*(B)) 
\otimes_{({\cal B}^\omega \: \otimes_{\C \Gamma} \: C^\infty_c(B))} 
C^\infty_{{\cal B}^\omega}(\widehat{M}; \widehat{E}) \right).
\end{equation}

We now couple $A^{Bismut}_s$ to the connection $\nabla^{can}$ of Section 
\ref{Section 2}
in order to obtain a superconnection
\begin{equation} \label{3.6}
A_s \: : \: C^\infty_{{\cal B}^\omega}(\widehat{M}; \widehat{E}) \rightarrow
\Omega^*(B, {\cal B}^\omega) \otimes_{C^\infty(B, {\cal B}^\omega)}
C^\infty_{{\cal B}^\omega}(\widehat{M}; \widehat{E}).
\end{equation}
Explicitly, 
\begin{equation} \label{3.7}
A_s \: = \: s \: Q \: + \: \nabla^u \: - \: \frac{1}{4s} \: c(T) \: + \:
\sum_{\gamma \in \Gamma} d\gamma \: \otimes \: h \: \gamma^{-1} \cdot
\end{equation}
Let ${\cal R}$ be the rescaling operator on 
$\Omega^{even}(B, {\cal B}^\omega)_{ab}$ which multiplies an element of
$\Omega^{2k}(B, {\cal B}^\omega)_{ab}$ by 
$(2 \pi i)^{-k}$.
Doing a Duhamel expansion around $e^{- \: (A^{Bismut}_s)^2}$ and using the
fact that $h$ has compact support, we can define
\begin{equation} \label{3.8}
e^{- \: A_s^2} \in
\Hom^\infty_{C^\infty(B, {\cal B}^\omega)} 
\left(
C^\infty_{{\cal B}^\omega}(\widehat{M}; \widehat{E}), \:
\Omega^*(B, {\cal B}^\omega) \otimes_{C^\infty(B, {\cal B}^\omega)}
C^\infty_{{\cal B}^\omega}(\widehat{M}; \widehat{E})
\right).
\end{equation}
and hence also define
${\cal R} \: \Tr_{s, <e>} \left( e^{- \: A_s^2} \right) \in 
\Omega^*(B, {\cal B}^\omega)_{ab}$.
From the superconnection formalism 
\cite[Chapter 1.4]{Berline-Getzler-Vergne (1992)},
${\cal R} \: \Tr_{s, <e>} \left( e^{- \: A_s^2} \right)$ is closed and its
cohomology class is independent of $s \: > \: 0$; see
\cite[Theorem 3.1]{Heitsch (1995)} for a detailed proof in the
analogous case of ${\cal R} \: \Tr_s \left(e^{- \: (A^{Bismut}_s)^2} \right)$.

\begin{theorem} \label{Theorem 2}
\begin{equation} \label{3.9}
\lim_{s \rightarrow 0} {\cal R} \: 
\Tr_{s, <e>} \left( e^{- \: A_s^2} \right) \: = \:
\int_Z \phi \: \widehat{A} \left( \nabla^{TZ} \right) \: 
\ch \left( \nabla^{\widehat{V}} \right) \: \ch \left(\nabla^{can} \right)
\: \in \: \Omega^*(B, {\cal B}^\omega)_{ab}.
\end{equation}
\end{theorem}
\begin{pf}
We use a variation of the proof of 
\cite[Theorem 10.23]{Berline-Getzler-Vergne (1992)}. As in
\cite[Theorem 10.23]{Berline-Getzler-Vergne (1992)}, we must first
establish a Lichnerowicz-type formula for $A_s^2$.
Let $\{e_i\}_{i=1}^{\dim(Z)}$ be a local orthonormal basis for
$TZ$ and 
let $\{c^i\}_{i=1}^{\dim(Z)}$ be Clifford algebra generators, with
$(c^i)^2 \: = \: -1$. Let $\{\tau^\alpha\}_{\alpha = 1}^{\dim(B)}$ be
a local basis of $T^*B$ and let 
$E^\alpha$ denote exterior multiplication
by $\tau^\alpha$.
Bismut proved a Lichnerowicz-type formula for $(A^{Bismut}_s)^2$
\cite[Theorem 10.17]{Berline-Getzler-Vergne (1992)}, namely
\begin{align} \label{3.10}
(A^{Bismut}_s)^2 \: = & 
\: s^2 \: {\cal D}^* {\cal D} \: + \: s^2 \: \frac{1}{4} \: r^Z \: + \:
\frac{1}{4} \: \sum_{i,j} F_{ij}(\widehat{V}) \: [c^i, c^j] \: + \:
\sum_{i,\alpha} F_{\alpha i}(\widehat{V}) \: E^\alpha \: c^i \: + \notag \\
& \: \frac{1}{4} \: \sum_{\alpha, \beta} F_{\alpha \beta}(\widehat{V}) \: 
[E^\alpha, E^\beta],
\end{align}
where ${\cal D}$ is a certain vertical differentiation operator and
$r^Z \in C^\infty(\widehat{M})$ is the scalar curvature function
of the fibers.
Then from the formula (\ref{3.7}) for $A_s$, one finds
\begin{align} \label{3.11}
A_s^2 \: = & \: s^2 \: {\cal D}^* {\cal D} \: + \: \frac{1}{4}
\: s^2 \: r^Z \: + \:
\frac{1}{4} \: \sum_{i,j} F_{ij}(\widehat{V}) \: [c^i, c^j] \: + \:
\sum_{\alpha,i} F_{\alpha i}(\widehat{V}) \: E^\alpha \: c^i \: + \notag \\
& \frac{1}{4} \: \sum_{\alpha, \beta} F_{\alpha \beta}(\widehat{V}) \: 
[E^\alpha, E^\beta] \: - \:
s \: \sum_{\gamma \in \Gamma} d\gamma \: 
\left( c(d^{vert} h) \: + \: E(d^{hor} h) \right) \: \gamma^{-1}  \: -
\notag \\
& \: \sum_{\gamma, \gamma^\prime} (\gamma \gamma^\prime \cdot h) \:
(\gamma \cdot h) \: d\gamma \: d\gamma^\prime \: (\gamma \gamma^\prime)^{-1}.
\end{align}

We now perform a Getzler rescaling, as in 
\cite[p. 342]{Berline-Getzler-Vergne (1992)}. Explicitly, we send
$\partial_{x^j} \rightarrow s^{-1} \: \partial_{x^j}$,
$c^j \rightarrow E^j \: - \: I^j \rightarrow
s^{-1} \: E^j \: - \: s \: I^j$ and $\tau^\alpha \rightarrow \tau^\alpha$.
Then following
\cite[Proposition 10.28]{Berline-Getzler-Vergne (1992)},
one finds that in the rescaling limit $A_s^2$ becomes
\begin{equation} \label{3.12}
- \: \sum_{j=1}^{\dim(Z)} \left( \partial_{x^j} \: - \: \frac{1}{8}
\sum_{j=1}^{\dim(Z)} \sum_{a,b=1}^{\dim(\widehat{M})} R^Z_{abij} \: x^j \:
E^a \: E^b \right)^2 \: + \: (\nabla^{\widehat{V}})^2 \: + \: (\nabla^{can})^2.
\end{equation}
The rest of the proof now proceeds as in the proof of
\cite[Theorem 10.21]{Berline-Getzler-Vergne (1992)}; compare 
\cite[Pf. of Prop. 12]{Lott (1992)}.
\end{pf}

Let us note that the right-hand-side of (\ref{3.9}) pairs with closed graded
traces on $\Omega^*(B, \C \Gamma)$, and not just closed graded traces on
$\Omega^*(B, {\cal B}^\omega)$.
In the construction of $\ch \left(\nabla^{can} \right)$, we can allow
$h$ to be a Lipschitz function on $\widehat{M}$ 
(see \cite[Lemma 4]{Lott (1992)}, where $\ch \left(\nabla^{can} \right)$
is called $\widetilde{\omega}_h$). 
Let $c \: : \: \widehat{M} \rightarrow E\Gamma$ be a
$\Gamma$-equivariant
classifying map for the fiber bundle
$\widehat{M} \rightarrow M$. It is defined
up to $\Gamma$-homotopy. 
 As $M$ is compact, we may assume that $c$ is Lipschitz with
respect to a piecewise Euclidean $\Gamma$-invariant metric on the simplicial
complex $E\Gamma$. 
Let $\eta$ be a closed graded $n$-trace on 
$\Omega^*(B, \C \Gamma)$. Then we can describe
$\langle \int_Z \phi(z) \: \widehat{A} \left( \nabla^{TZ} \right) \: 
\ch \left( \nabla^{\widehat{V}} \right) \: \ch \left(\nabla^{can} \right),
\eta \rangle$ as follows.

First, let us dualize (\ref{1.44})
with respect to $B$ to write 
\begin{equation} \label{3.13}
\langle \ch \left( \nabla^{univ} \right), \eta \rangle \in
\left( \bigoplus_{a+b=n} \bigoplus_{l=0}^\infty
\Omega^{a-2l}(E\Gamma) \: \widehat{\otimes} 
\: \Omega^{\dim(B)-b}_\tau(B) \right)^\Gamma,
\end{equation}
where $\Omega^{*}_\tau(B)$ denotes the differential
forms on $B$ with distributional coefficients and with value in the
flat orientation line bundle $o(TB)$.
Passing to a $\Z_2$-grading, we obtain 
\begin{equation} \label{3.14}
\langle \ch \left( \nabla^{univ} \right), \eta \rangle \in
\left( \Omega^{n + \dim(B) + 2 \Z}_\tau(E\Gamma \times B) \right)^\Gamma.
\end{equation}
By construction, the form
that we have obtained is closed, so we have
an element $\Phi_\eta \in 
\HH^{n + \dim(B) + 2 \Z}_\tau((E\Gamma \times B)/\Gamma)$. 

The map $(c, \pi) \: : \: \widehat{M} \rightarrow
E\Gamma \times B$ is $\Gamma$-equivariant and so descends to a
classifying map $\nu \: : \: M \rightarrow (E \Gamma \times B)/\Gamma$.
Let $T{\cal F}$
denote the leafwise tangent bundle on $M$ with respect to the
foliation ${\cal F}$, a vector bundle on $M$. 
Put $V \: = \: \widehat{V}/\Gamma$.
Then we claim that
\begin{equation} \label{3.15}
\langle \int_Z \phi(z) \: \widehat{A} \left( \nabla^{TZ} \right) \: 
\ch \left( \nabla^{\widehat{V}} \right) \: \ch \left(\nabla^{can} \right),
\eta \rangle \: = \:
\int_M \widehat{A}(T{\cal F}) \: \ch(V) \: \nu^* \Phi_\eta \: \in \: \C.
\end{equation}
To see this, take the pairing of
$\ch \left(\nabla^{can} \right)$ and $\eta$ with respect to
$\C \Gamma$,
to get a $\Gamma$-invariant element of
$\Omega^*(\widehat{M}) \: \widehat{\otimes} \: C_*(B)$. 
Dualizing with respect to $B$,
we obtain an element of $\Omega^*(\widehat{M}) \: \widehat{\otimes} \: 
\Omega^*_\tau(B)$.
Applying the product $(\omega_1, \omega_2) \rightarrow
\omega_1 \wedge \pi^*\omega_2$ to this, 
we finally obtain a closed $\Gamma$-invariant 
element of $\Omega^*_\tau(\widehat{M})$, the latter being the
differential forms on $\widehat{M}$ with distributional coefficients
and with value in $\pi^*(o(TB))$. Hence we have an element of
$\HH^*_\tau(M)$ (where the $\tau$ now refers to the flat orientation line
bundle $o(N{\cal F})$), which we denote by
$\langle \ch \left(\nabla^{can} \right), \eta \rangle$.
Let $* ( \widehat{A}(T{\cal F}) \: \ch(V)) \in \HH_*^\tau(M)$
be the Poincar\'e dual of 
$\widehat{A}(T{\cal F}) \: \ch(V) \in \HH^*(M; \R)$.
Then the left-hand-side of (\ref{3.15}) is the pairing between 
$\langle \ch \left(\nabla^{can} \right), \eta \rangle \in \HH^*_\tau(M)$ 
and $* ( \widehat{A}(T{\cal F}) \: \ch(V)) \in \HH_*^\tau(M)$.
Now we may also compute $\langle \ch \left(\nabla^{can} \right), \eta \rangle$
by using the Lipschitz function $c^* j$ instead of $h$ in constructing
$\nabla^{can}$. (We may have to smooth $\omega_1$ before taking the
product.) Then by naturality,
$\langle \ch \left(\nabla^{can} \right), \eta \rangle
\: = \: \nu^* \Phi_\eta$, which proves the claim.

\section{The Index and the Superconnection Chern Character} \label{Section 5}

In this section we prove Theorem \ref{Theorem 3}, relating the
Chern character of $A_s$ to the Chern character of the index. 
We define the index by means of the index projection and show
that its Chern character can be computed by means of a connection $\nabla$.
We then show that the Chern character of the index
can also be written as the supertrace of $e^{- \: (\nabla^\prime)^2}$ for
a certain $\Z_2$-graded connection
$\nabla^\prime$.

We relate the Chern character of the index to the superconnection Chern
character by means of a homotopy from $\nabla^\prime$ to $A_s$.  
This is done in three cases.  In the
first case, that of finitely-generated projective modules, the naive homotopy
argument works.  In the second case, that of the families index theorem, we 
show that smoothing factors in the homotopy allow the naive argument to
be carried through.  In the third case, that of Theorem \ref{Theorem 3},
we again justify the naive homotopy argument. We then give some
geometric consequences of Theorem \ref{Theorem 3}.

Let ${\frak A}$ be an algebra over $\C$ and let
$\Omega^*$ be a GDA equipped with a homomorphism $\rho \: : \: 
{\frak A} \rightarrow \Omega^0$. Let ${\cal E}$ be a left ${\frak A}$-module
and let $\nabla \: : \: {\cal E} \rightarrow \Omega^1 \otimes_{\frak A} 
{\cal E}$ be a connection on ${\cal E}$. 

Let $\widetilde{\Omega}^*$ be a subalgebra of the graded algebra
$\Hom_{\frak A}({\cal E}, \Omega^* \otimes_{\frak A} {\cal E})$. 
Put $\widetilde{\frak A} \: = \: \widetilde{\Omega}^0$.
We assume that
$\widetilde{\Omega}^*$ is closed under $\nabla$ and that the curvature
$\Theta \: = \: \nabla^2 \in
\Hom_{\frak A}({\cal E}, \Omega^2 \otimes_{\frak A} {\cal E})$ 
of the connection lies in $\widetilde{\Omega}^2$. 
Then $\nabla$ extends to a covariant differentiation
$\widetilde{\nabla} \: : \: 
\widetilde{\Omega}^* \rightarrow \widetilde{\Omega}^{*+1}$
on $\widetilde{\Omega}^*$ which satisfies
$\widetilde{\nabla}^2(\widetilde{\omega}) 
\: = \: \Theta \: \widetilde{\omega} \: - \: \widetilde{\omega} \: \Theta$.
Let $\widetilde{\eta} \: : \: 
\widetilde{\Omega}^* \rightarrow \C$ be an even graded trace which
satisfies $\widetilde{\eta}(\widetilde{\nabla} 
\widetilde{\omega}) \: = \: 0$ for all
$\widetilde{\omega} \in \widetilde{\Omega}^*$.

As in \cite[Chapter III.3, Lemma 9]{Connes (1994)}, 
let $X$ be a new formal odd variable of degree $1$ and put
\begin{equation} \label{4.1}
\widetilde{\widetilde{\Omega}}^* \: = \: \widetilde{\Omega} \: \oplus 
X \: \widetilde{\Omega}^* \: \oplus \: \widetilde{\Omega}^* \: X \: \oplus
\: X \: \widetilde{\Omega}^* \: X
\end{equation}
with the new multiplication rules
$(\widetilde{\omega}_1 \: X) \: \widetilde{\omega}_2 \: = \: 
\widetilde{\omega}_1 \: ( X \: \widetilde{\omega}_2) \: = \: 0$
and $(\widetilde{\omega}_1 \: X) \: ( X \: \widetilde{\omega}_2) \: = \: 
\widetilde{\omega}_1 \: \Theta \: 
\widetilde{\omega}_2$.  
Define a graded trace $\widetilde{\widetilde{\eta}}$
on $\widetilde{\widetilde{\Omega}}^*$ by
\begin{equation} \label{4.2}
\widetilde{\widetilde{\eta}}(\widetilde{\omega}_{1} \: + \: X \: 
\widetilde{\omega}_2 \: + \: 
\widetilde{\omega}_3 \: X \: +
\: X \: \widetilde{\omega}_4 \: X) \: = \: 
\widetilde{\eta}(\widetilde{\omega}_1) \: 
+ \: (-1)^{|\widetilde{\omega}_4|} \:
\widetilde{\eta}(\widetilde{\omega}_4). 
\end{equation}
Define a differential $d$ on $\widetilde{\widetilde{\Omega}}^*$ which is
generated by the relations 
\begin{equation} \label{4.3}
d \: \widetilde{\omega} \: = \: \widetilde{\nabla} \: 
\widetilde{\omega} \: + \: 
X \: \widetilde{\omega} \: + \:
(-1)^{\widetilde{\omega}} \: \widetilde{\omega} \: X
\end{equation}
and $d \: X \: = \: 0$. One can check that $d^2 \: = \: 0$ and
$\widetilde{\widetilde{\eta}}(d \widetilde{\widetilde{\omega}}) \: = \: 0$ for 
$\widetilde{\widetilde{\omega}} \in \widetilde{\widetilde{\Omega}}^*$. That is,
$\left( \widetilde{\widetilde{\Omega}}^*, d, 
\widetilde{\widetilde{\eta}} \right)$ defines
a cycle over $\widetilde{\frak A}$ in the sense of 
\cite[Chapter III.1.$\alpha$, Definition 1]{Connes (1994)}.

Suppose that $\widetilde{\frak A}$ is unital.
The cycle structure induces a map from
$\KK_0(\widetilde{\frak A})$ to $\C$ 
\cite[Chapter III.3, Proposition 2]{Connes (1994)}. 
To fix normalizations, let $p \in \widetilde{\frak A}$ be a projection with
corresponding class $[p] \in \KK_0(\widetilde{\frak A})$.
Then the pairing of the Chern character of $[p]$ with $\widetilde{\eta}$ is 
defined to be 
\begin{equation} \label{4.4}
\langle \ch([p]), \widetilde{\eta} \rangle \: = \:
(2 \pi i )^{- \: deg(\widetilde{\eta})/2} \: 
\widetilde{\widetilde{\eta}} \left(
p \: e^{- \: p \: dp \: dp} \: p \right).
\end{equation}
One can check that $p \: dp \: dp \: p \: = \: p \: (\widetilde{\nabla} p) \:
(\widetilde{\nabla} p) \: p \: + \: p \: \Theta \: p$, which in turn equals the
curvature of the connection $p \circ \nabla \circ p$.
Thus 
\begin{equation} \label{4.5}
\langle \ch([p]), \widetilde{\eta} \rangle \: = \: 
(2 \pi i )^{- \: deg(\widetilde{\eta})/2} \: \widetilde{\eta} \left(
p \: e^{- \: (p \circ \nabla \circ p)^2} \: p \right).
\end{equation}
This is consistent with well-known formulae 
if ${\cal E}$ is a finitely-generated projective ${\frak A}$-module, 
but we have not assumed that ${\cal E}$ is finitely-generated projective.
The equation
extends to $p \in M_n(\widetilde{\frak A})$ in the standard way.

We will need an extension of this formula to the nonunital case.
We suppose again that we have the algebra $\widetilde{\Omega}^*$ and
the connection 
$\widetilde{\nabla}$ on it.
In general,
$\widetilde{\nabla}^2$ may not be given in terms of an element of
$\widetilde{\Omega}^2$.
Instead, as in \cite[Section 2]{Nistor (1997)}, we make the weaker
assumption that
$\widetilde{\nabla}^2$ comes from a multiplier
$(l, r)$ of $\widetilde{\Omega}^*$. This means that
$l$ and $r$ are linear maps from $\widetilde{\Omega}^*$ to itself
such that for all $\widetilde{\omega}_1, \widetilde{\omega}_2 \in
\widetilde{\Omega}^*$, we have
$l(\widetilde{\omega}_1 \widetilde{\omega}_2) \: = \: l(\widetilde{\omega}_1) 
\: \widetilde{\omega}_2$,
$r(\widetilde{\omega}_1 \widetilde{\omega}_2) \: = \: \omega_1 \: 
r(\widetilde{\omega}_2)$ and
$\widetilde{\omega}_1 \: l(\widetilde{\omega}_2) \: = \: 
r(\widetilde{\omega}_1) \: \widetilde{\omega}_2$.
Then we assume that $\widetilde{\nabla}^2 (\widetilde{\omega}) \: = \:
l(\widetilde{\omega}) \: - \: r(\widetilde{\omega})$ for some
$(l, r)$ of degree $2$. (If $\widetilde{\frak A}$ is
unital then we recover $\Theta$ by $\Theta \: = \: l(1) \: = \: r(1)$, and
$(l, r)$ are given in terms of ${\Theta}$ by
$l(\widetilde{\omega}) \: = \: \Theta \: \widetilde{\omega}$,
$r(\widetilde{\omega}) \: = \: \widetilde{\omega} \: \Theta$.) With
this understanding,
$p \: dp \: dp \: p \: = \: p \: (\widetilde{\nabla} p) \:
(\widetilde{\nabla} p) \: p \: + \: p \: l(p)$ and equation
(\ref{4.5}) still makes sense for $p \in \widetilde{\frak A}$.

Next, recall that if $\widetilde{\frak A}$ is nonunital and 
$\widetilde{\frak A}^+$ is
the algebra obtained by adding a unit to 
$\widetilde{\frak A}$, with canonical
homomorphism $\pi \: : \: \widetilde{\frak A}^+ \rightarrow \C$, then
$\KK_0(\widetilde{\frak A}) \: = \: \Ker \left( \pi_* \: :\: 
\KK_0(\widetilde{\frak A}^+) \rightarrow \KK_0(\C) \right)$.
Thus an element of $\KK_0(\widetilde{\frak A})$ 
can be represented as $p \: - \: p_0$ with the projections
$p, p_0 \in M_n(\widetilde{\frak A}^+)$ 
satisfying $\pi_*(p) \: = \: \pi_*(p_0) \in
M_n(\C)$.
Then the equation
\begin{equation} \label{4.19}
\langle \ch([p \: - \: p_0]), \widetilde{\eta} \rangle \: = \: 
(2 \pi i )^{- \: deg(\widetilde{\eta})/2} \: \widetilde{\eta} \left(
p \: e^{- \: (p \circ \nabla \circ p)^2} \: p  \: - \: 
p_0 \: e^{- \: (p_0 \circ \nabla \circ p_0)^2} \: p_0 \right).
\end{equation}
gives a well-defined map on $\KK_0(\widetilde{\frak A})$.

\subsection{Finitely-generated projective ${\frak A}$-modules}

Now suppose that ${\cal E} \: = \: {\cal E}^+ \: \oplus \: {\cal E}^-$ is 
$\Z_2$-graded, with
${\cal E}^\pm$ finitely-generated projective ${\frak A}$-modules. 
We assume that $\nabla$ preserves the grading.
Put $\widetilde{\Omega}^* \: = \:
\Hom_{\frak A}({\cal E}, \Omega^* \otimes_{\frak A} {\cal E})$. We assume
that $\widetilde{\frak A} \: = \: \End_{\frak A}({\cal E})$ has a
holomorphic functional calculus.  For example, it suffices that
${\frak A}$ be a complete locally convex topological algebra.
Put
$e \: = \: 
\begin{pmatrix}
1 & 0 \\
0 & 1
\end{pmatrix}
$ and
$v \: = \: 
\begin{pmatrix}
1 & 0 \\
0 & -1
\end{pmatrix}
$.

Given
$D \in \Hom_{\frak A}({\cal E}^+, {\cal E}^-)$ and 
$D^* \in \Hom_{\frak A}({\cal E}^-, {\cal E}^+)$, we assume that the
spectra of $D D^*$ and $D^* D$ are contained in the nonnegative reals.
We construct an index projection following
\cite{Connes-Moscovici (1990)} and
\cite{Moscovici-Wu (1994)}.
Let $u \in C^\infty(\R)$ be an even function such that
$w(x) \: = \: 1 \: - \: x^2 \: u(x)$ is a Schwartz function and
the Fourier transforms of $u$ and $w$ have compact support
\cite[Lemma 2.1]{Moscovici-Wu (1994)}. 
Define $\overline{u} \in C^\infty([0, \infty))$ by $\overline{u}(x)
\: = \: u(x^2)$.
Put ${\cal P} \: = \: \overline{u}(D^* D) D^*$, which we will think of as a
parametrix for $D$, and put
$S_+ \: = \: I \: - \: {\cal P}D$, $S_- \: = \: I \: - \: D{\cal P}$.
Consider the operator 
\begin{equation}
l \: = \: 
\begin{pmatrix}
S_+ & - \: (I + S_+) {\cal P} \\
D & S_-
\end{pmatrix},
\end{equation}
with inverse
\begin{equation}
l^{-1} \: = \: 
\begin{pmatrix}
S_+ & \: {\cal P}(I + S_-) \\
- \: D & S_-
\end{pmatrix}.
\end{equation}
The index projection is defined by  
\begin{equation} \label{4.6}
p \: = \: l \: \frac{e+v}{2} \: 
l^{-1} \: = \: 
\begin{pmatrix}
S_+^2 & S_+ (I + S_+) {\cal P} \\
S_- D & I - S_-^2
\end{pmatrix}.
\end{equation}
Put
\begin{equation} \label{4.8}
p_0 \: = \: \frac{e-v}{2} \: = \: 
\begin{pmatrix}
0 & 0 \\
0 & 1
\end{pmatrix}.
\end{equation}

We note that the supertrace $\Tr_s \: : \: 
\End_{\frak A}({\cal E}) 
\rightarrow {\frak A}/[{\frak A}, {\frak A}]$ is given by
\begin{align} \label{4.9}
\Tr_s \left( M \right) \: & = \:  \Tr \left( \frac{e+v}{2} \:
M \: \frac{e+v}{2} \right) \: - \: \Tr 
\left( \frac{e-v}{2} \: M \: \frac{e-v}{2} \right) \\ 
& = \: \Tr \left( p \: l M l^{-1} \: p \right) \: - \: \Tr \left( p_0
\: M \: p_0 \right). \notag
\end{align}
Let us define a new connection on ${\cal E}$ by
\begin{equation} \label{4.10}
\nabla^\prime \: = \: \left(
\frac{e+v}{2} \: l^{-1} \circ \nabla \circ l \: 
\frac{e+v}{2} \right)\: + \: \left( \frac{e-v}{2} \: \nabla \:
\frac{e-v}{2} \right).
\end{equation}
Then by construction,
\begin{equation} \label{4.11}
l \circ \nabla^\prime \circ l^{-1} \: = \: \left( p \circ \nabla \circ p 
\right) \: + \: \left(
(1-p) \: l \circ \nabla \circ l^{-1} \: (1-p) \right).
\end{equation}
In particular,
\begin{equation} \label{4.12}
\Tr \left( p \: e^{- (p \circ \nabla \circ p)^2} \: p \right) \: = \: 
\Tr \left( p \: l \: e^{- (\nabla^\prime)^2} \: l^{-1} \: p \right).
\end{equation}
Also, from (\ref{4.10}), we have
\begin{equation} \label{4.13}
\Tr \left( p_0 \: e^{- \: \nabla^2} \: p_0 \right) \: = \: 
\Tr \left( p_0 \: e^{- (\nabla^\prime)^2} \: p_0 \right).
\end{equation}
Using (\ref{4.9}), we see that we can write $\ch([p \: - \: p_0])$ 
as the Chern character form of a connection on ${\cal E}$, namely
\begin{equation} \label{4.14}
\ch([p \: - \: p_0]) \: = \:
{\cal R} \: \Tr_s \left( e^{- \: (\nabla^\prime)^2} \right).
\end{equation}

For future use, we note that
$(\nabla^\prime)^+ \: = \: S_+ \nabla^+ S_+ \: + \: 
{\cal P} (I + S_-) \nabla^- D$ and $(\nabla^\prime)^- \: = \: \nabla^-$.

So far in this section we have been working with algebraic tensor products.
If ${\frak A}$ is a complete locally convex topological algebra
then it is straightforward
to extend the statements to the topological setting.

\begin{lemma} \label{Lemma 9}
Suppose that $A$ is a superconnection on ${\cal E}$.
Then for any closed graded trace $\widetilde{\eta}$ on
$\Hom_{\frak A}({\cal E}, \Omega^* \otimes_{\frak A} {\cal E})$,
\begin{equation} \label{4.17}
\langle \ch([p \: - \: p_0]),
\widetilde{\eta} \rangle \: = \:
\widetilde{\eta} \left( \ch \left( A \right) \right).
\end{equation}
\end{lemma}
\begin{pf}
In general, if $\{A(t)\}_{t \in [0,1]}$ is a smooth $1$-parameter family
of superconnections on ${\cal E}$ then
\begin{equation} \label{4.18}
\ch(A(1)) \: - \: \ch(A(0)) \: = \: {\cal R} \: d \: 
\int_{[0,1]} \Tr_s \left( \frac{dA}{dt} \: e^{- \: A(t)^2} \right) \: dt.
\end{equation}
Thus it suffices to construct a smooth path in the space of superconnections
between $\nabla^\prime$ and $A$, for example the linear homotopy
$A(t) \: = \: t \: A \: + \: (1-t) \: \nabla^\prime$.
\end{pf}
\noindent
\subsection{Fiber bundles}
 
Let $\pi : M \rightarrow B$ be a fiber bundle with closed
even-dimensional fiber $Z$.
Endow the fiber bundle
with a vertical Riemannian metric $g^{TZ}$ and a horizontal
distribution $T^HM$. 
Let $E$ be a Hermitian vector bundle on $M$ which is a fiberwise
Clifford bundle, with
compatible connection $\nabla^E$. Let ${\cal E}$ be the smooth sections of
the $\Z_2$-graded
vector bundle $\pi_!(E)$
on $B$, whose fiber over $b$ is $C^\infty(Z_b; E \big|_{Z_b})$.

Take ${\frak A} \: = \: C^\infty(B)$, 
$\Omega^* \: = \: \Omega^*(B)$ and let $\nabla^{\pm}$ be the natural
Hermitian connection on ${\cal E}^{\pm}$ 
\cite[Proposition 10.10]{Berline-Getzler-Vergne (1992)}.
Let $\widetilde{\Omega}^*$ be the subalgebra of
$\Hom_{C^\infty(B)} \left( {\cal E}, \Omega^*(B) \otimes_{C^\infty(B)}
{\cal E} \right)$ consisting of elements with a smooth fiberwise
integral kernel $K(z, w)$.
Put $\widetilde{\frak A} \: = \: \widetilde{\Omega}^0$, a nonunital
algebra if $\dim(Z) \: > \: 0$.
Given a closed current $\eta \in Z_{even}(B; \R)$, let
$\widetilde{\eta}$ be the graded trace on $\widetilde{\Omega}^*$ given by
$\widetilde{\eta}(K) \: = \: \int_{\eta} (\int_Z K(z,z) \: d\vol_Z)$. 

Let $D \: : \: {\cal E}^+ \rightarrow {\cal E}^-$ be the vertical Dirac
operator. Define the index projection $p \in 
\widetilde{\frak A}$ as in (\ref{4.6}). Define $e$ and $v$ as before.
Then the index of $D$ is defined to
be $\Ind(D) \: = \:
[p \: -  \: p_0] \in \KK_0(\widetilde{{\frak A}})$. The Chern character
of $\Ind(D)$ pairs with $\eta$ by (\ref{4.19}).

Let us note that although $p \: 
e^{- \: (p \circ \nabla \circ p)^2} \: p$ and $p_0 \: 
e^{- \: (p_0 \circ \nabla \circ p_0)^2} \: p_0$ 
may not individually lie in $\widetilde{\Omega}^*$,
their difference does.  For example, the component in $\widetilde{\Omega}^0$
is 
\begin{equation} \label{4,20}
\begin{pmatrix}
S_+^2 & S_+ (I + S_+) {\cal P} \\
S_- D & I - S_-^2
\end{pmatrix} \: - \: 
\begin{pmatrix}
0 & 0 \\
0 & 1
\end{pmatrix}
\: = \: 
\begin{pmatrix}
S_+^2 & S_+ (I + S_+) {\cal P} \\
S_- D & - S_-^2
\end{pmatrix}.
\end{equation}
This is related to the fact that the K-theory of a nonunital algebra
is defined in terms of the K-theory of the algebra obtained by adding
a unit.

For $s \: > \: 0$, let $A_s^{Bismut}$ denote the Bismut superconnection on
$\pi_!E$.

\begin{proposition} \label{Proposition 4}
$\langle \ch(\Ind(D)), \widetilde{\eta} \rangle \: = \:
\widetilde{\eta} \left( \ch(A_s^{Bismut}) \right)$.
\end{proposition}
\begin{pf}
We can homotop from $D$ to $s D$ in the definition of the index
projection without changing the $\KK$-theory class and then everywhere
change $s D$ to $D$, to easily reduce to 
the case $s \: = \: 1$. 
Define $\nabla^\prime$ as in (\ref{4.10}). As in the proof of Lemma 
\ref{Lemma 9}, we wish to 
homotop from $\nabla^\prime$ to $A_1^{Bismut}$ and then apply (\ref{4.18}). 
The only issue is see that the formal expressions are
well-defined.

First, for $t \in [0,1]$ put
\begin{equation} \label{4.21}
A(t) \: = \: 
\begin{pmatrix}
(\nabla^\prime)^+ & t \: D^* \\
t \: D &  (\nabla^\prime)^-
\end{pmatrix}.
\end{equation}
Write
\begin{equation}
e^{- \: A(t)^2} \: = \:
\begin{pmatrix}
{e^{- \: A(t)^2}}_{11} & {e^{- \: A(t)^2}}_{12} \\
{e^{- \: A(t)^2}}_{21} & {e^{- \: A(t)^2}}_{22}
\end{pmatrix}.
\end{equation}
Put
\begin{equation} \label{add1}
\ch(A(t)) \: = \: 
{\cal R} \: \left( \Tr \left( S_+ \: {e^{- \: A(t)^2}}_{11} \: S_+ \right)
\: + \: \Tr
\left( D \: {e^{- \: A(t)^2}}_{11} \: {\cal P} \: (I \: + \: S_-) \:  - 
\: {e^{- \: A(t)^2}}_{22} \right) \right).
\end{equation}
Formally, the right-hand-side of (\ref{add1}) equals
${\cal R} \left(  \Tr \left( \: {e^{- \: A(t)^2}}_{11} \right) \: - \:
\Tr \left( \: {e^{- \: A(t)^2}}_{22} \right) \right)$.
To see that the traces in the right-hand-side
of (\ref{add1}) make sense, let us compute $A(t)^2$.
In an ungraded notation, we have
\begin{equation} \label{add2}
A(t)^2 \: = \: 
\begin{pmatrix}
((\nabla^\prime)^+)^2 \: + \: t^2 \: D^* D & 
t \: [(\nabla^\prime)^-, D^*] \: + \: t ((\nabla^\prime)^+ \: - \:
(\nabla^\prime)^-) D^* \\
t \: [(\nabla^\prime)^+, D] \: - \: t ((\nabla^\prime)^+ \: - \:
(\nabla^\prime)^-)D &  ((\nabla^\prime)^-)^2 \: + \: t^2 DD^*
\end{pmatrix}.
\end{equation}
The term in the lower left-hand corner of (\ref{add2}) is
\begin{align}
t \: [(\nabla^\prime)^+, D] \: - \: t ((\nabla^\prime)^+ \: - \:
(\nabla^\prime)^-)D \: & = \: 
t \: [(\nabla^\prime)^-, D] \: - \: t D 
((\nabla^\prime)^+ \: - \:
(\nabla^\prime)^-) \\
& = \: 
t \: [\nabla^-, D] \: - \: t D 
( S_+ \nabla^+ S_+ \: + \: {\cal P} (I \: + \: S_-) \nabla^- D \: - \:
\nabla^-) \notag \\
&  = \: t \left( S_-^2 \nabla^- D \: - \: D S_+ \nabla^+ S_+ \right). \notag
\end{align}
Then modulo uniformly
smoothing operators,
\begin{equation} \label{4.22}
A(t)^2 \equiv 
\begin{pmatrix}
((\nabla^\prime)^+)^2 \: + \: t^2 \: D^* D  &  
t \: [(\nabla^\prime)^-, D^*] \: + \: t ((\nabla^\prime)^+ \: - \:
(\nabla^\prime)^-) D^*  \\
0 &  D \:
\left[  ((\nabla^\prime)^+)^2 \: + \: t^2 \: D^* D) \right] 
\: {\cal P}
\end{pmatrix}
\end{equation}
and
\begin{equation} \label{add3}
e^{- \: A(t)^2} \equiv 
\begin{pmatrix}
e^{- \left( ((\nabla^\prime)^+)^2 \: + \: t^2 \: D^* D \right)} &  
{\cal Z}  \\
0 &  
D \:
e^{- \left( ((\nabla^\prime)^+)^2 \: + \: t^2 \: D^* D \right)} 
\: {\cal P}
\end{pmatrix},
\end{equation}
where
\begin{align} \label{add4}
{\cal Z} \: = \: - \: \int_0^1 
e^{- \: u \: \left( ((\nabla^\prime)^+)^2 \: + \: t^2 \: D^* D \right)} 
& \left( t \: [(\nabla^\prime)^-, D^*] \: + \: t ((\nabla^\prime)^+ \: - \:
(\nabla^\prime)^-) D^* \right) \\
& e^{- \: (1-u) \:
\left( (\nabla^-)^2 \: + \: t^2 \: D D^* \right)} \: du. \notag
\end{align}
It follows that the right-hand-side of (\ref{add1}) is well-defined.

We claim that 
$\langle \ch(\Ind(D)), \widetilde{\eta} \rangle \: = \:
\widetilde{\eta} \left( \ch(A(0)) \right)$. To see this,
we have
\begin{align}
p \: e^{- \: (p \circ \nabla \circ p)^2} \: p \: & = \:
l \:
\begin{pmatrix}
e^{- \: ((\nabla^\prime)^+)^2} & 0 \\
0 & 0
\end{pmatrix}
\: l^{-1} \\
& = \:
\begin{pmatrix}
S_+ \: e^{- \: ((\nabla^\prime)^+)^2} \: S_+ &
S_+ \: e^{- \: ((\nabla^\prime)^+)^2} \: {\cal P} \: (I \: + \: S_-) \\
D \: e^{- \: ((\nabla^\prime)^+)^2} \: S_+ &
D \: e^{- \: ((\nabla^\prime)^+)^2} \: {\cal P} \: (I \: + \: S_-) 
\end{pmatrix}. \notag
\end{align}
Then
\begin{align}
\Tr \left( p \: e^{- \: (p \circ \nabla \circ p)^2} \: p \: - \: 
p_0 \: e^{- \: (p_0 \circ \nabla \circ p_0)^2} \: p_0 \right) \: = \:
& \Tr \left( S_+ \: e^{- \: ((\nabla^\prime)^+)^2} \: S_+ \right) \: +\\
& 
\Tr \left( D \: e^{- \: ((\nabla^\prime)^+)^2} \: {\cal P} \: (I \: + \: S_-)
\: - \: e^{- \: (\nabla^-)^2} \right), \notag
\end{align}
from which the claim follows.

Let us note that the terms being traced in (\ref{add1}) are
in fact uniformly smoothing with respect to $t$, due to factors of the
form $S_\pm$.

We now wish to write the analog of equation 
(\ref{4.18}). Although $\ch(A(t))$ is well-defined, it is not clear
in the present setting that
the integrand in (\ref{4.18}) is integrable for small $t$.
Let us first do a formal
calculation. With respect to (\ref{4.21}), (\ref{add3}) and (\ref{add4}), 
we have
\begin{equation}
\frac{dA}{dt} \: = \: \begin{pmatrix}
0 & D^* \\
D & 0
\end{pmatrix}
\end{equation}
and
\begin{align}
& \Tr_s \left( \frac{dA}{dt} \: 
\begin{pmatrix}
e^{- \left( ((\nabla^\prime)^+)^2 \: + \: t^2 \: D^* D \right)} &  
{\cal Z}  \\
0 &  
D \:
e^{- \left( ((\nabla^\prime)^+)^2 \: + \: t^2 \: D^* D \right)} 
\: {\cal P}
\end{pmatrix} \right) \: = \:
- \: \Tr \left( D \: {\cal Z} \right) \: = \\
& t \: \Tr \left( D \: 
\int_0^1 
e^{- \: u \: \left( ((\nabla^\prime)^+)^2 \: + \: t^2 \: D^* D \right)} 
\: \left( [(\nabla^\prime)^-, D^*] \: + \: ((\nabla^\prime)^+ \: - \:
(\nabla^\prime)^-) D^* \right)
 \right. \notag \\
& \left.  
\: \: \: \: \: \: \: \: \: \: \: \:
\: \: \: \: \: \: \: \: \: \: \: \: e^{- \: (1-u) \:
\left( (\nabla^-)^2 \: + \: t^2 \: D D^* \right)} \: du \right). \notag
\end{align}
Modulo smoothing operators,
\begin{align}
& D \: 
\int_0^1 
e^{- \: u \: \left( ((\nabla^\prime)^+)^2 \: + \: t^2 \: D^* D \right)} 
\: \left( [(\nabla^\prime)^-, D^*] \: + \: ((\nabla^\prime)^+ \: - \:
(\nabla^\prime)^-) D^* \right) \: e^{- \: (1-u) \:
\left( (\nabla^-)^2 \: + \: t^2 \: D D^* \right)} \: du \: \equiv \: \\
& \int_0^1 
e^{- \: u \: \left( (\nabla^-)^2 \: + \: t^2 \: D D^* \right)} 
\: D \: \left( [(\nabla^\prime)^-, D^*] \: + \: ((\nabla^\prime)^+ \: - \:
(\nabla^\prime)^-) D^* \right) \: e^{- \: (1-u) \:
\left( (\nabla^-)^2 \: + \: t^2 \: D D^* \right)} \: du \: \equiv \: \notag 
\\  
& \int_0^1 
e^{- \: u \: \left( (\nabla^-)^2 \: + \: t^2 \: D D^* \right)} 
\: [\nabla^-, D D^*] \: e^{- \: (1-u) \:
\left( (\nabla^-)^2 \: + \: t^2 \: D D^* \right)} \: du. \notag
\end{align}
Then
\begin{align}
& \Tr \left( 
\int_0^1 
e^{- \: u \: \left( (\nabla^-)^2 \: + \: t^2 \: D D^* \right)} 
\: [\nabla^-, D D^*] \: e^{- \: (1-u) \:
\left( (\nabla^-)^2 \: + \: t^2 \: D D^* \right)} \: du \right) \: = \: \\
& \Tr \left( 
[\nabla^-, D D^*] \: e^{- \: 
\left( (\nabla^-)^2 \: + \: t^2 \: D D^* \right)} \right) \: = \:
- \:  t^{-2} \: d \: \Tr \left( 
e^{- \: \left( (\nabla^-)^2 \: + \: t^2 \: D D^* \right)} \right) \notag
\end{align}
The upshot is that we can write
\begin{equation} \label{add5}
\ch(A(1)) \: - \: \ch(A(0)) \: = \: {\cal R} \: d \: 
\int_{[0,1]} \left(
\Tr_s \left( \frac{dA}{dt} \: e^{- \: A(t)^2} \right)
\: + \: t^{-1} \: d \: \Tr \left( 
e^{- \: \left( (\nabla^-)^2 \: + \: t^2 \: D D^* \right)} \right) \right)
\: dt,
\end{equation}
where the integrand in the right-hand-side of (\ref{add5}), after the terms
are appropriately grouped, is 
the trace of a smoothing operator that is continuous in $t$.
Thus
$\widetilde{\eta} \left( \ch(A(0)) \right) \: = \:
\widetilde{\eta} \left( \ch(A(1)) \right)$.

Next, we perform a linear homotopy from $A(1)$ to $A_1^{Bismut}$.
As the $0$-th order part of the superconnection always equals
$\begin{pmatrix}
0 &  D^*  \\
D &  0
\end{pmatrix}$
during this homotopy, it is easy to justify the formal superconnection 
argument, using (\ref{4.18}), that
$\widetilde{\eta} \left( \ch(A(1)) \right) \: = \:
\widetilde{\eta} \left( \ch(A_1^{Bismut}) \right)$. This proves the
proposition.
\end{pf}
\subsection{Fiber Bundles over Cross-Product Groupoids}

Using the notation of Section \ref{Section 4}, put
${\frak A} \: = \: C^\infty(B, {\cal B}^\omega)$, 
$\Omega^* \: = \: \Omega^\infty(B, {\cal B}^\omega)$,
${\cal E} \: = \: C^\infty_{{\cal B}^\omega}(\widehat{M};
\widehat{E})$,
$\widetilde{\frak A} \: = \: 
\End^\infty_{C^\infty(B, {\cal B}^\omega)}
\left(
C^\infty_{{\cal B}^\omega}(\widehat{M};
\widehat{E}) \right)$
and $$\widetilde{\Omega}^* \: = \: 
\Hom^\infty_{C^\infty(B, {\cal B}^\omega)}
\left(
C^\infty_{{\cal B}^\omega}(\widehat{M};
\widehat{E}), \:
\Omega^*(B, {\cal B}^\omega)
\otimes_{C^\infty(B, {\cal B}^\omega)} 
C^\infty_{{\cal B}^\omega}(\widehat{M};
\widehat{E})
\right).$$ Define the index projection as in (\ref{4.6}).
Let $\eta$ be a closed graded trace on $\Omega^*(B, {\cal B}^\omega)$.
Then we can go through the same steps as in the proof of
Proposition \ref{Proposition 4} to conclude 
\begin{theorem} \label{Theorem 3}
\begin{equation} \label{4.30}
\langle \ch(\Ind(D)), \widetilde{\eta} \rangle \: = \:
\langle {\cal R} \: 
\Tr_{s, <e>} \left( e^{- \: A_s^2} \right), \eta \rangle.
\end{equation}
\end{theorem}
We note that we use finite-propagation-speed estimates
in order to know that we can
carry out the arguments in $\widetilde{\Omega}^*$. That is,
we use the fact that if $f \in \Sigma^{(\infty)}(\R)$ is a function whose 
Fourier transform
$\widehat{f}(k)$ has exponential decay in $|k|$ then the Schwartz kernel 
$f(sQ)(z, w)$ has exponential decay in $d(z, w)$. In order to obtain
uniform decay estimates in the analog of the 
last step of Proposition \ref{Proposition 4}, as in 
\cite[Theorem 9.48]{Berline-Getzler-Vergne (1992)} we use the
fact that in the $n$-th order term of the Duhamel expansion of 
$e^{- \: A_s^2}$, there is always a factor of the form
$e^{- \: r \: s^2 \: Q^2}$ with $r \: \ge \: \frac{1}{n+1}$.

Putting together Theorems \ref{Theorem 2} and \ref{Theorem 3}, we obtain
\begin{theorem} \label{Theorem 4}
\begin{equation} \label{4.31}
\langle \ch(\Ind(D)), \widetilde{\eta} \rangle \: = \:
\int_M \widehat{A}(T{\cal F}) \: \ch(V) \: \nu^* \Phi_\eta.
\end{equation}
\end{theorem}

If $Z$ is instead odd-dimensional then one can prove Theorem
\ref{Theorem 4} by a standard trick involving taking the product with
a circle.

\begin{corollary} \label{Corollary 2}
Let ${\cal A}$ be a subalgebra of the reduced cross-product $C^*$-algebra
$C_0(B) \rtimes_r \Gamma$ which is stable under the holomorphic functional
calculus in $C_0(B) \rtimes_r \Gamma$ and which contains
$C^\infty(B, {\cal B}^\omega)$. Let $\eta$ be a closed graded trace on
$\Omega^*(B, \C \Gamma)$ which extends to give a cyclic cocycle on
${\cal A}$. Suppose that $TZ$ is spin and that $g^{TZ}$ has fiberwise
positive scalar curvature.
Then $\int_M \widehat{A}(T{\cal F}) \: \nu^* \Phi_\eta \: = \: 0$.
\end{corollary}
\begin{pf}
Let $D$ be the pure Dirac operator. 
As $\eta$ is a closed graded trace on $\Omega^*(B, \C \Gamma)$, it gives
rise to a cyclic cocycle on $C^\infty_c(B) \rtimes \Gamma$ through its
character \cite[Section III.1.$\alpha$]{Connes (1994)}. By assumption,
this has an extension $\eta^\prime$ to ${\cal A}$. Now we have
$\Ind(D) \in \KK_*(C_0(B) \rtimes_r \Gamma)
\: \cong \: \KK_*({\cal A})$. Then (\ref{4.31}) becomes
\begin{equation}
\langle \Ind(D), \eta^\prime \rangle \: = \:
\int_M \widehat{A}(T{\cal F}) \: \nu^* \Phi_\eta.
\end{equation}
However, by the Lichnerowicz argument, $\Ind(D) \: = \: 0$. The
corollary follows.
\end{pf}

Suppose that $B \: = \: S^1$, with $\Gamma$ acting by orientation-preserving
diffeomorphisms.  There is
a left action of $\Gamma$ on $\Omega^1(B)$. 
Let $v \in \Omega^1(B)$ be a volume form.
Define a closed graded trace on $\Omega^*(B, \C \Gamma)$ by
\begin{equation} \label{4.32}
\eta(f \: g_0 \: dg_1 \: dg_2) \: = \: \int_B
f \: \left( \ln \frac{v}{g_0 g_1 \cdot v} \: d\ln 
\frac{g_0 g_1 \cdot v}{g_0 \cdot v} \: - \: 
\ln \frac{g_0 g_1 \cdot v}{g_0 \cdot v} \: d\ln 
\frac{v}{g_0 g_1 \cdot v} \right).
\end{equation}
Then $\nu^* \Phi_\eta$ is proportionate to the Godbillon-Vey class
${\rm GV} \in \HH^3(M; \R)$ 
\cite[Chapter III.6.$\beta$, Theorem 17]{Connes (1994)}.
Furthermore, the hypotheses of Corollary \ref{Corollary 2} are satisfied
\cite[Chapter III.7.$\beta$]{Connes (1994)}.

\begin{corollary} \label{Corollary 3}
Suppose that $B \: = \: S^1$, $\Gamma$ acts on $B$ by orientation-preserving
diffeomorphisms, 
$TZ$ is spin and $g^{TZ}$ has fiberwise
positive scalar curvature. Then
$\int_M \widehat{A}(T{\cal F}) \: {\rm GV} \: = \: 0$.
\end{corollary}

In general, if $\dim(B) \: = \: q$ and the action of $\Gamma$ on $B$ is
orientation-preserving then one can write down a closed
graded trace $\eta$ on $\Omega^*(B, \C \Gamma)$ so that
$\nu^* \Phi_\eta$ is proportionate to the Godbillon-Vey class
${\rm GV} \in \HH^{2q+1}(M; \R)$, and the above results extend.

\section{\'Etale Groupoids} \label{Section 6}

In this section we generalize the results of the previous sections
from cross-product groupoids to general smooth Hausdorff \'etale
groupoids.  In Subsection \ref{Subsection 6.1} we explain in detail how,
in the case of cross-product groupoids, the expressions of this section
specialize to the expressions of the previous sections.

We follow the groupoid conventions 
of \cite[Sections II.5 and III.2.$\delta$]{Connes (1994)}.
Let $G$ be a smooth Hausdorff groupoid, with units $G^{(0)}$. We suppose that
$G$ is \'etale, i.e. that the range map $r : G \rightarrow G^{(0)}$ and 
the source map $s : G \rightarrow G^{(0)}$
are local diffeomorphisms.
To construct the product of $\gamma_0, \gamma_1 \in G$, we must have
$s(\gamma_0) \: = \: r(\gamma_1)$. Then $r(\gamma_0 \gamma_1) \: = \:
r(\gamma_0)$ and $s(\gamma_0 \gamma_1) \: = \: s(\gamma_1)$.
Given $x \in G^{(0)}$, put $G^x \: = \: r^{-1}(x)$, $G_x \: = \: 
s^{-1}(x)$ and $G_x^x \: = \: G^x \cap G_x$.

Given
$f_0, f_1 \in C^\infty_c(G)$, the convolution product is
\begin{equation} \label{5.1} 
(f_0 \: f_1)(\gamma) \: = \: \sum_{\gamma_0  \gamma_1 \: = \: \gamma}
f_0(\gamma_0) \: f_1(\gamma_1).
\end{equation} 
The sum in (\ref{5.1}) is finite.

We write $G^{(n)}$ for the $n$-chains of composable elements of
$G$, i.e.
\begin{equation} \label{5.2}
G^{(n)} \: = \: \{ (\gamma_1, \ldots, \gamma_n) \in G^{n} \: : \:
s(\gamma_1) \: = \: r(\gamma_2), \cdots, s(\gamma_{n-1}) \: = \:
r(\gamma_n)\}.
\end{equation}
As $G$ is \'etale, $G^{(n)}$ is a manifold of the same dimension as
$G$. As in \cite[Section III.2.$\delta$]{Connes (1994)}, we define a 
double complex
by letting $\Omega^{m,n}_c(G)$ be the quotient of
$\Omega^m_c \left( G^{(n+1)} \right)$ by the forms which are supported
on $\{(\gamma_0, \ldots, \gamma_n) \: : \: \gamma_j \text{ is a unit for some }
j \: > \:  0 \}$.  The product of 
$\omega_1 \in \Omega_c^{*,n_1}(G)$ and $\omega_2 \in \Omega_c^{*,n_2}(G)$ is
given by
\begin{align} \label{5.3} 
(\omega_1 \omega_2)(\gamma_0, \ldots, \gamma_{n_1 + n_2}) \: = \:
& \sum_{\gamma \gamma^\prime \: = \: \gamma_{n_1}} 
\omega_1(\gamma_0, \ldots, \gamma_{n_1 - 1}, \gamma) \: \wedge \:
\omega_2(\gamma^\prime, \gamma_{n_1 + 1}, \ldots, \gamma_{n_1 + n_2}) \: - \:
\\
& \: \left( (-1)^{n_1 - 1} \:
\sum_{\gamma \gamma^\prime \: = \: \gamma_0}
\omega_1(\gamma, \gamma^\prime,
\gamma_{1}, \ldots, \gamma_{n_1 - 1}) \: + \: \right. \notag \\
& 
(-1)^{n_1 - 2} \:
\sum_{\gamma \gamma^\prime \: = \: \gamma_1}
\omega_1(\gamma_0, \gamma, \gamma^\prime,
\gamma_{2}, \ldots, \gamma_{n_1 - 1})
 \: + \ldots + \: \notag \\ 
& \left. \sum_{\gamma \gamma^\prime \: = \: \gamma_{n_1 - 1}}
\omega_1(\gamma_0, \ldots, \gamma_{n_1 - 2}, \gamma, \gamma^\prime) \right)
\wedge 
\omega_2(\gamma_{n_1}, \ldots, \gamma_{n_1 + n_2}). \notag
\end{align}
In forming the wedge product in (\ref{5.3}), the maps $r$ and $s$ are used
to identify cotangent spaces.
The first differential $d_1$ on $\Omega_c^{*,*}(G)$ 
is the de Rham differential.
To define the second differential $d_2$, 
let $\chi_{G^{(0)}} \in C^\infty(G)$ be
the characteristic function for the units.  Then
\begin{equation} \label{5.4}
(d_2 \omega)(\gamma_0, \ldots, \gamma_{n + 1}) \: = \:
\chi_{G^{(0)}}(\gamma_0) \: \omega(\gamma_1, \ldots, \gamma_n).
\end{equation}
We let $\Omega_c^*(G)$ denote the GDA formed by the total complex of
$\Omega_c^{*,*}(G)$.

Let $P$ be a smooth $G$-manifold \cite[Section II.10.$\alpha$,
Definition 1]{Connes (1994)}. That is, first of all, there is a submersion
$\pi \: : \: P \rightarrow G^{(0)}$. Given $x \in G^{(0)}$, we write
$Z_x \: = \: \pi^{-1}(x)$. Putting 
\begin{equation} \label{5.5}
P \times_{G^{(0)}} G \: = \: \{ (p, \gamma) \in P \times G \: : \:
p \: \in \: Z_{r(\gamma)} \},
\end{equation}
we must also have a map $P \times_{G^{(0)}} G \rightarrow P$, denoted
$(p, \gamma) \rightarrow p \gamma$, such that
$p  \gamma \: \in \: Z_{s(\gamma)}$ and  
$(p \gamma_1) \gamma_2 \: = \:
p (\gamma_1 \gamma_2)$ for all $(\gamma_1, \gamma_2) \in G^{(2)}$.
It follows that for each $\gamma \in G$, the map $p \rightarrow p  \gamma$
gives a diffeomorphism from $Z_{r(\gamma)}$ to 
$Z_{s(\gamma)}$.
The groupoid ${\cal G} \: = \: P \rtimes G$ has underlying space
$P \times_{G^{(0)}} G$, units ${\cal G}^{(0)} \: = \: P$ and
maps $r(p, \gamma) \: = \: p$ and $s(p, \gamma) \: = \: p \gamma$.

We assume that $P$ is a proper $G$-manifold, 
i.e. that the map $P \times_{G^{(0)}} G
\rightarrow P \times P$ given by $(p, \gamma) \rightarrow (p, p \gamma)$
is proper. Then ${\cal G} \: = \: 
P \rtimes G$ is a proper groupoid, i.e.
the map ${\cal G} \rightarrow {\cal G}^{(0)} \times {\cal G}^{(0)}$ given by
$\gamma \rightarrow (r(\gamma), s(\gamma))$ is proper 
\cite[Section II.10.$\alpha$, Definition 2]{Connes (1994)}. We also assume
that $G$ acts cocompactly on $P$, i.e. that the quotient of $P$ by the
equivalence relation ($p \sim p^\prime$ if $p = p^\prime \gamma$
for some $\gamma \in G$) is compact.  Equivalently,
${\cal G} \: = \: P \rtimes G$ is a cocompact groupoid, i.e. the quotient of 
${\cal G}^{(0)}$ by the
equivalence relation ($x \sim x^\prime$ if $(x , x^\prime) \: = \:
(r(\gamma), s(\gamma))$ for some $\gamma \in {\cal G}$) is compact.
 Finally, we assume that $G$ acts freely on $P$, i.e. that the
preimage of the diagonal in $P \times P$ under the map
$P \times_{G^{(0)}} G
\rightarrow P \times P$ equals
$P \times_{G^{(0)}} G^{(0)}$. Equivalently, 
${\cal G} \: = \: P \rtimes G$ is a free groupoid, i.e.
the preimage of the diagonal in ${\cal G}^{(0)} \times {\cal G}^{(0)}$ under 
the map
$(r, s) \: : \: {\cal G} \rightarrow {\cal G}^{(0)} \times {\cal G}^{(0)}$
equals ${\cal G}^{(0)}$. 

Now let ${\cal G}$ be any free proper cocompact \'etale groupoid.
The product $C^\infty_c({\cal G}) \times 
C^\infty_c \left( {\cal G}^{(0)} \right) \rightarrow
C^\infty_c \left( {\cal G}^{(0)} \right)$ is given explicitly by
\begin{equation} \label{5.6}
(fF)(x) \: = \: \sum_{\gamma \in {\cal G}^x} f(\gamma) \: F(s(\gamma)).
\end{equation} 
We wish to define a connection 
\begin{equation} \label{5.7}
\nabla^{can} \: : \: C^\infty_c \left( {\cal G}^{(0)} \right) \rightarrow
\Omega_c^1({\cal G}) \otimes_{C^\infty_c ( {\cal G} )}
C^\infty_c \left( {\cal G}^{(0)} \right).
\end{equation}
To do so, we use isomorphisms
$\Omega_c^{1,0}({\cal G}) \otimes_{C^\infty_c ( {\cal G} )}
C^\infty_c \left( {\cal G}^{(0)} \right) \equiv 
\Omega^1_c \left( {\cal G}^{(0)} \right)$ and
$\Omega_c^{0,1}({\cal G}) \otimes_{C^\infty_c ( {\cal G} )}
C^\infty_c \left( {\cal G}^{(0)} \right) \equiv
C^\infty_c({\cal G})/C^\infty_c({\cal G}^{(0)})$. 
The latter isomorphism is realized by saying that
the image of $\omega \otimes F$
in $C^\infty_c({\cal G})/C^\infty_c({\cal G}^{(0)})$ is given by 
\begin{equation} \label{5.8}
(\omega F)(\gamma_0) \: = \:
\sum_{\gamma \gamma^\prime \: = \: \gamma_0} \omega(\gamma, \gamma^\prime) \:
F(s(\gamma_0)) \: - \: \sum_{\gamma \in {\cal G}^{s(\gamma_0)}} 
\omega(\gamma_0, \gamma) \: F(s(\gamma))
\end{equation}
for $\gamma_0 \notin {\cal G}^{(0)}$.
Then with this isomorphism, the multiplication
\begin{equation} \label{5.9}
C^\infty_c ( {\cal G} ) \times
\left( \Omega_c^{0,1}({\cal G}) \otimes_{C^\infty_c ( {\cal G} )}
C^\infty_c \left( {\cal G}^{(0)} \right) \right) \rightarrow
\Omega_c^{0,1}({\cal G}) \otimes_{C^\infty_c ( {\cal G} )}
C^\infty_c \left( {\cal G}^{(0)} \right),
\end{equation}
i.e. the multiplication
\begin{equation} \label{5.10}
C^\infty_c ( {\cal G} ) \times
\frac{C^\infty_c({\cal G})}{C^\infty_c({\cal G}^{(0)})}
\rightarrow
\frac{C^\infty_c({\cal G})}{C^\infty_c({\cal G}^{(0)})},
\end{equation}
is given by
\begin{equation} \label{5.11}
(f{\cal F})(\gamma_0) \: = \: \sum_{\gamma \gamma^\prime \: = \: \gamma_0}
f(\gamma) \: {\cal F}(\gamma^\prime) \: - \:
f(\gamma_0) \: \sum_{\gamma \in {\cal G}^{s(\gamma_0)}} {\cal F}(\gamma)
\end{equation}
for $\gamma_0 \notin {\cal G}^{(0)}$.

More generally, there is an isomorphism between
$\Omega_c^{m,n}({\cal G}) \otimes_{C^\infty_c({\cal G})} C^\infty_c
\left( {\cal G}^{(0)} \right)$ and the quotient of
$\Omega_c^m \left( {\cal G}^{(n)} \right)$ by the forms which are supported
on $$\{(\gamma_0, \ldots, \gamma_{n-1}) 
\: : \: \gamma_j \text{ is a unit for some }
j \: \ge \:  0 \},$$ under which the image of $\omega \otimes F$ is given by
\begin{align} \label{5.12}
(\omega F)(\gamma_0, \ldots, \gamma_{n-1}) \:  = \: &
\sum_{\gamma \gamma^\prime \: = \: \gamma_0} 
\omega(\gamma, \gamma^\prime, \gamma_1, \ldots, \gamma_{n-1}) \: 
F(s(\gamma_{n-1})) 
\: - \: \\ 
& \sum_{\gamma \gamma^\prime \: = \: \gamma_1} 
\omega(\gamma_0, 
\gamma, \gamma^\prime, \gamma_2, 
\ldots, \gamma_{n-1}) \: F(s(\gamma_{n-1})) \: 
+ \ldots + \notag \\
& (-1)^{n-1} \: \sum_{\gamma \gamma^\prime \: = \: \gamma_{n-1}} 
\omega(\gamma_0, \ldots, \gamma_{n-2}, \gamma, \gamma^\prime) \: 
F(s(\gamma_{n-1}) \: + \notag \\
& (-1)^n \: \sum_{\gamma \in {\cal G}^{s(\gamma_{n-1})} } 
\omega(\gamma_0, \ldots, \gamma_{n-1}, \gamma) \: F(s(\gamma)) \notag
\end{align}
for $\gamma_0, \ldots, \gamma_{n-1} \notin {\cal G}^{(0)}$.

Now let $h \in C^\infty_c \left( {\cal G}^{(0)} \right)$ satisfy
\begin{equation} \label{5.13}
\sum_{\gamma \in {\cal G}^x} h(s(\gamma)) \: = \: 1
\end{equation}
for all $x \in {\cal G}^{(0)}$.
Then there is a connection 
\begin{equation} \label{5.14}
\nabla^{can} \: = \: \nabla^{1,0} \: \oplus \: \nabla^{0,1}
\end{equation}
on $C^\infty_c \left( {\cal G}^{(0)} \right)$
where $\nabla^{1,0} (F) \in \Omega^1_c \left( {\cal G}^{(0)} \right)$ 
is the de Rham differential of $F \in 
C^\infty_c \left( {\cal G}^{(0)} \right)$ and $\nabla^{0,1} (F) \in
\Omega_c^{0,1}({\cal G}) \otimes_{C^\infty_c ( {\cal G} )}
C^\infty_c \left( {\cal G}^{(0)} \right) \equiv
C^\infty_c({\cal G})/C^\infty_c({\cal G}^{(0)})$ is given by
\begin{equation} \label{5.15}
(\nabla^{0,1} (F))(\gamma_0) \: = \: F(r(\gamma_0)) \: h(s(\gamma_0))
\end{equation}
for $\gamma_0 \notin {\cal G}^{(0)}$.

One sees that 
$\left( \nabla^{can} \right)^2 \in 
\Hom_{C^\infty_c({\cal G})}
\left( C^\infty_c \left( {\cal G}^{(0)} \right), 
\Omega_c^2({\cal G}) \otimes_{C^\infty_c ( {\cal G} )}
C^\infty_c \left( {\cal G}^{(0)} \right) \right)$ acts on 
$C^\infty_c \left( {\cal G}^{(0)} \right)$ as left multiplication by a
$2$-form $\Theta$ which
commutes with $C^\infty_c({\cal G})$.
Explicitly, $\Theta \: = \: \Theta^{1,1} \: + \: \Theta^{0,2}$ where
\begin{equation} \label{5.16}
\Theta^{1,1}(\gamma_0, \gamma_1) \: = \: - \:
\chi_{{\cal G}^{(0)}}(\gamma_0 \gamma_1) \: d^{de \: Rham}h(s(\gamma_0))
\end{equation}
and
\begin{equation} \label{5.17}
\Theta^{0,2}(\gamma_0, \gamma_1, \gamma_2) \: = \: - \:
\chi_{{\cal G}^{(0)}}(\gamma_0 \gamma_1 \gamma_2) 
\: h(s(\gamma_0)) \: h(s(\gamma_1))
\end{equation}
for $\gamma_1, \gamma_2 \notin {\cal G}^{(0)}$.
Put
\begin{equation} \label{5.18}
\ch(\nabla^{can}) \: = \: e^{- \: 
\frac{ \Theta}{2 \pi i}} \in 
\End_{\Omega^*_c({\cal G})} \left( \Omega^*_c({\cal G})
\otimes_{C^\infty_c({\cal G})} C^\infty_c({\cal G}^{(0)}) \right).
\end{equation}
Then the abelianization of $\ch(\nabla^{can})$ is closed
and its cohomology class is
independent of the choice of $h$.

Now suppose that $G$ acts freely, properly and cocompactly on $P$. Give $P$ a
$G$-invariant fiberwise Riemannian metric.
An element $K$ of $\End_{C^\infty_c(G)}(C^\infty_c(P))$ has a Schwartz
kernel $K(p|p^\prime)$ with respect to its fiberwise action, so that
we can write
\begin{equation} \label{5.19}
(KF)(p) \: = \: \int_{Z_{\pi(p)}} K(p|p^\prime) \: F(p^\prime) \:
d\vol_{Z_{\pi(p)}}.
\end{equation}
Let $\End^\infty_{C^\infty_c(G)}(C^\infty_c(P))$ denote the subalgebra
of $\End_{C^\infty_c(G)}(C^\infty_c(P))$ consisting of elements with a
smooth integral kernel. 

Let $\phi \in C^\infty_c (P)$ satisfy
\begin{equation} \label{5.20}
\sum_{\gamma \in G^{\pi(p)}} \phi(p  \gamma) \: = \: 1
\end{equation}
for all $p \in P$. Define a trace $
\Tr_{<e>}$ on $\End^\infty_{C^\infty_c(G)}(C^\infty_c(P))$ by
\begin{equation} \label{5.22}
\Tr_{<e>}(K)(\gamma_0) \: = \: \int_{Z_{\gamma_0}} \phi(p) \:
K(p|p) \:  d\vol_{Z_{\gamma_0}},
\end{equation}
for $\gamma_0 \in G^{(0)}$. Then $\Tr_{<e>}$ takes value in
$\frac{C^\infty_c(G)}{[C^\infty_c(G), C^\infty_c(G)]}$ and is
concentrated at the units.

Put
\begin{equation} \label{5.23}
G^{(n)} \times_s P \: = \: \{ (\gamma_0, \ldots, \gamma_{n-1}, p) \in
G^{(n)} \times P \: : \: p \in Z_{s(\gamma_{n-1})} \}.
\end{equation}
There is an isomorphism, as in (\ref{5.12}), between
$\Omega_c^{m,n}(G) \otimes_{C^\infty_c(G)} C^\infty_c(P)$ and
the quotient of $\Omega^m_c(G^{(n)} \times_s P)$ by the forms which are
supported on 
\begin{equation} \label{5.24}
\{(\gamma_0, \ldots, \gamma_{n-1}, p) 
\: : \: \gamma_j \text{ is a unit for some }
j \: \ge \:  0 \}.
\end{equation}
Consider the $\Z$-graded algebra
\begin{equation} \label{5.25}
\Hom^\infty_{C^\infty_c(G)} \left( C^\infty_c(P), \Omega_c^*(G) 
\otimes_{C^\infty_c(G)} C^\infty_c(P) \right)
\end{equation}
consisting of elements $K$ of 
$\Hom_{C^\infty_c(G)} \left( C^\infty_c(P), \Omega_c^*(G) 
\otimes_{C^\infty_c(G)} C^\infty_c(P) \right)$ with a smooth integral kernel.
Using the above isomorphism, the kernel of an element $K$
can be written in the form
\begin{equation} \label{5.26}
K(\gamma_0, \ldots, \gamma_{n-1}, p | p^\prime) \in
\Lambda^m(T^*_{r(\gamma_0)} G^{(0)})
\end{equation}
where $(\gamma_0, \ldots, \gamma_{n-1}, p) \in G^{(n)} \times_s P$,
$\gamma_0, \ldots, \gamma_{n-1} \notin G^{(0)}$, 
and $p^\prime \in Z_{r(\gamma_{0})}$.
The action
of $K$ on $C^\infty_c(P)$ is given by
\begin{equation} \label{5.27}
(KF)(\gamma_0, \ldots, \gamma_{n-1}, p) \: = \: 
\int_{Z_{r(\gamma_0)}} K(\gamma_0, \ldots, \gamma_{n-1}, p | p^\prime) \: 
F(p^\prime) \:
d\vol_{Z_{r(\gamma_0)}}
\end{equation}
for $\gamma_0, \ldots, \gamma_{n-1} \notin G^{(0)}$. Then there is a trace
\begin{equation} \label{5.28}
\Tr_{<e>} \: : \: \Hom^\infty_{C^\infty_c(G)}
\left( C^\infty_c(P), \Omega_c^*(G) 
\otimes_{C^\infty_c(G)} C^\infty_c(P) \right) \rightarrow
\Omega_c^*(G)_{ab}
\end{equation}
given by
\begin{align} \label{5.29} 
\Tr_{<e>}(K)(\gamma_0, \ldots, \gamma_n) \: = \: & 
\chi_{G^{(0)}}(\gamma_0 \cdots \gamma_n) \: \int_{Z_{r(\gamma_0)}} \phi(p) 
\\
& \left[ 
\chi_{G^{(0)}}(\gamma_0) \: \sum_{\gamma \gamma^\prime \: = \: \gamma_n}
\: 
K(\gamma_1, \ldots, \gamma_{n-1}, \gamma, p ( \gamma^\prime)^{-1} | p) 
\: - \:
\right. \notag \\
&  \chi_{G^{(0)}}(\gamma_0) \: \sum_{\gamma \gamma^\prime \: = \: 
\gamma_{n-1}} \:
K(\gamma_1, \ldots, \gamma_{n-2}, \gamma, \gamma^\prime, p ( \gamma_n)^{-1} | 
p) \: + \: \notag \\
& \chi_{G^{(0)}}(\gamma_0) \: \sum_{\gamma \gamma^\prime \: = \: 
\gamma_{n-2}} \:
K(\gamma_1, \ldots, \gamma_{n-3}, \gamma, \gamma^\prime, \gamma_{n-1}, 
p ( \gamma_n)^{-1} | 
p) \: + \notag \\
& \ldots \: + \notag \\
& 
(-1)^{n-1} \chi_{G^{(0)}}(\gamma_0) \: \sum_{\gamma \gamma^\prime \: = \: 
\gamma_{1}} \:
K(\gamma, \gamma^\prime, \gamma_2, \ldots,  \gamma_{n-1}, 
p ( \gamma_n)^{-1} | 
p) \: + \notag \\
& \left. (-1)^n \: 
K(\gamma_0, \ldots,  \gamma_{n-1}, 
p ( \gamma_n)^{-1} | 
p)  \right] \: d\vol_{Z_{r(\gamma_0)}}. \notag
\end{align}

Let $L$ be a topological space which is the total space of a submersion
$\sigma : L \rightarrow G^{(0)}$. We suppose that
each fiber $L_x \: = \: \sigma^{-1}(x)$ is
a complete length space with metric $d_x$. 
We also assume that $G$ acts isometrically, properly
and cocompactly on $L$. Let $i : G \rightarrow L$ be a $G$-equivariant 
map, not necessarily continuous. That is, for each $x \in G^{(0)}$, $i$ sends
$G_x$ to $L_x$ and for each $\gamma \in G$, the composite map
$G_{r(\gamma)} \stackrel{\cdot \gamma}{\longrightarrow} G_{s(\gamma)}
\stackrel{i}{\longrightarrow} L_{s(\gamma)}$ equals the composite map
$G_{r(\gamma)} \stackrel{i}{\longrightarrow} L_{r(\gamma)}
\stackrel{\cdot \gamma}{\longrightarrow} L_{s(\gamma)}$.
We assume that $i$ is proper in the sense that the preimage of a
compact set has compact closure. Note
that $i$ gives a possibly-discontinuous section of $\sigma$. 
We assume in addition that for any compact subset $K$ of $G^{(0)}$,
$i(K)$ has compact closure. 
 (The existence of $L$, $\sigma$ and $i$ is a definite
assumption.) 
Define a ``length function'' on $G$ by 
\begin{equation} \label{5.30}
l(\gamma) \: = \: 
d_{s(\gamma)}( i(s(\gamma)) , i(\gamma)),
\end{equation}
where we think of $\gamma$ and $s(\gamma) \in
G^{(0)}$ as living in $G_{s(\gamma)}$. Then $l(\gamma_0 \gamma_1) \: \le \:
l(\gamma_0) \: + \: l(\gamma_1)$. Furthermore, for each $x \in G^{(0)}$,
the restriction of $l$ to $G_x$ is proper.

Let $C^\infty_\omega(G)$ be the set of $f \in C^\infty(G)$ such that\\
1. $s^* f$ has support in some compact subset $K$ of $G^{(0)}$ and\\
2. For all $q \in \Z^+$,
\begin{equation} \label{5.31}
\sup_{x \in K} \sup_{\gamma \in G_x} 
\: e^{q l(\gamma)} |f(\gamma)| \: < \: \infty,
\end{equation}
along with the analogous property for derivatives.
Then $C^\infty_\omega(G)$ is an algebra with the same formal multiplication
as in (\ref{5.1}).  It is independent
of the choices of $L$ and $i$; 
compare \cite[Proposition 3]{Lott (1992)}. 
We define $\Omega^*_\omega(G)$ similarly. That is, first define
$s_n$ on $G^{(n+1)}$ by $s_n(\gamma_0, \ldots, \gamma_n) \: = \:
s(\gamma_n)$. Let
$\widetilde{\Omega}^{m,n}_\omega(G)$
be the elements $\omega$ of 
$\Omega^m(G^{n+1})$ such \\
1. $s_n^* \omega$ has support in some compact subset $K$ of $G^{(0)}$ and\\
2. For all $q \in \Z^+$,
\begin{equation} \label{5.32}
\sup_{x \in K} \sup_{(\gamma_0, \ldots, \gamma_n) \in s_n^{-1}(x)} 
\: e^{q l(\gamma_0 \ldots \gamma_n)} |\omega(\gamma_0, \ldots, \gamma_n)| 
\: < \: \infty,
\end{equation}
along with the analogous property for derivatives.  
Let $\Omega^*_\omega(G)$ be the quotient of 
$\widetilde{\Omega}^{m.n}_\omega(G)$ by the forms which are supported on
$\{(\gamma_0, \ldots, \gamma_n) \: : \: \gamma_j \text{ is a unit for some }
j \: > \: 0\}$.
Then $\Omega^*_\omega(G)$ is a GDA, with the same formal multiplication
as in (\ref{5.3}).

Suppose now that $G$ acts freely, properly and cocompactly on $P$ as before.
Put 
\begin{equation} \label{5.33}
C^\infty_\omega(P) \: = \: C^\infty_\omega(G) \otimes_{C^\infty_c(G)}
C^\infty_c(P).
\end{equation} Using the cocompactness of the $G$-action on $P$, the
elements of $C^\infty_\omega(P)$ can be characterized as elements
$F \in C^\infty(P)$ such that for any $x \in G^{(0)}$, $p \in Z_x$ and
$q \in \Z^+$, we have
\begin{equation} \label{5.34}
\sup_{z \in Z_x} \: e^{q \: d(z, x)} \: |F(z)| \: < \: \infty,
\end{equation}
along with the analogous property for the covariant derivatives of $F$.
Let 
$\End^\infty_{C^\infty_\omega(G)}
\left(
C^\infty_{\omega}(P)
\right)$ be the subalgebra of
$\End_{C^\infty_\omega(G)} \left( 
C^\infty_{\omega}(P) \right)$ consisting of elements
$K$ with
a smooth integral kernel $K(z,w)$.
Then the elements of $\End^\infty_{C^\infty_\omega(G)}
\left(
C^\infty_{\omega}(P)
\right)$ can be characterized as the $G$-invariant
elements $K(z| w) \in C^\infty(
P \times_{G^{(0)}} P)$ such for any $x \in G^{(0)}$ and $q \in \Z^+$, 
\begin{equation} \label{5.35}
\sup_{z,w \in Z_x} e^{q \: d(z,w)} \: |K(z|w)| \: < \: \infty,
\end{equation}
along with the analogous property for the covariant derivatives of $K$.
With the natural definition of $\Hom^\infty_{C^\infty_\omega(G)}
\left(
C^\infty_{\omega}
(P), \:
\Omega^*_\omega(G)
\otimes_{C^\infty_\omega(G)} 
C^\infty_{\omega}(P)
\right)$,
an element $K$ has a kernel as in (\ref{5.26}).
The formula (\ref{5.22}) extends to a trace
$\Tr_{<e>} \: : \: 
\End^\infty_{C^\infty_\omega(G)} \rightarrow
\frac{C^\infty_\omega(G)}{[C^\infty_\omega(G), C^\infty_\omega(G)]}$.
The formula (\ref{5.29}) extends to a trace
\begin{equation} \label{5.37}
\Tr_{<e>} \: : \: 
\Hom^\infty_{C^\infty_\omega(G)}
\left(
C^\infty_{\omega}
(P), \:
\Omega^*_\omega(G)
\otimes_{C^\infty_\omega(G)} 
C^\infty_{\omega}(P)
\right)
\rightarrow \Omega^*_\omega(G)_{ab}.
\end{equation}

If ${E}$ is a $\Z_2$-graded 
$G$-invariant Hermitian vector bundle on $P$,
with an invariant Hermitian connection, then we can define
$C^\infty_{\omega}(P; {E})$
and a supertrace
\begin{equation} \label{5.38}
\Tr_{s,<e>} \: : \:
\Hom^\infty_{C^\infty_\omega(G)}
\left(
C^\infty_{\omega}(P;
{E}), \:
\Omega^*_\omega(G)
\otimes_{C^\infty_\omega(G)} 
C^\infty_{\omega}(P;
{E})
\right)
\rightarrow \Omega^*_\omega(G)_{ab}.
\end{equation}

We now choose
a $G$-invariant vertical Riemannian metric $g^{TZ}$ on the
submersion $\pi \: : \: P \rightarrow G^{(0)}$ and a 
$G$-invariant horizontal
distribution $T^H P$.
Suppose that $Z$ is even-dimensional. Let $\widehat{E}$ be a
$\Gamma$-invariant Clifford
bundle on $P$ which is equipped 
with a $G$-invariant connection. For simplicity of notation, 
we asssume that
$\widehat{E} \: = \: S^Z \: \widehat{\otimes} \: \widehat{V}$, 
where $S^Z$ is a vertical spinor bundle
and $\widehat{V}$ 
is an auxiliary vector bundle on $P$. More precisely,
suppose that the vertical tangent
bundle $TZ$ has a spin structure.  Let $S^Z$ be the 
vertical spinor bundle, a $G$-invariant $\Z_2$-graded Hermitian
vector bundle on $P$. Let
$\widehat{V}$ be another $G$-invariant $\Z_2$-graded Hermitian
vector bundle on $P$ which is equipped 
with a $G$-invariant Hermitian connection. Then we put 
$\widehat{E} \: = \: S^Z \:
\widehat{\otimes} \: \widehat{V}$. The case of general 
$G$-invariant Clifford bundles $\widehat{E}$ can
be treated in a way completely analogous to what follows.

Let $Q$ denote the vertical Dirac-type operator acting on
$C^\infty_c(P; \widehat{E})$.
From finite-propagation-speed estimates
as in \cite[Pf. of Prop 8]{Lott (1992)}, 
along with
the bounded geometry
of $\{Z_x\}_{x \in G^{(0)}}$, 
for any $s \: > \: 0$ we have 
\begin{equation} \label{5.39}
e^{- \: s^2 \: Q^2} \in 
\End^\infty_{C^\infty_\omega(G)} 
\left(
C^\infty_{\omega}(P; \widehat{E})
\right).
\end{equation}

Let
\begin{equation} \label{5.40}
A^{Bismut}_s \: : \: C^\infty_c(P; \widehat{E}) \rightarrow
\Omega^*_c(G^{(0)}) \otimes_{C^\infty_c(G^{(0)})}
C^\infty_c(P; \widehat{E}) 
\end{equation}
denote the Bismut superconnection on the submersion
$\pi \: : \: P \rightarrow G^{(0)}$
\cite[Proposition 10.15]{Berline-Getzler-Vergne (1992)}. 
It is of the form
\begin{equation} \label{5.41}
A^{Bismut}_s \: = \: s \: Q \: + \: \nabla^u \: - \: \frac{1}{4s} \: c(T),
\end{equation}
where $\nabla^u$ is a certain Hermitian connection and
$c(T)$ is Clifford multiplication by the curvature $2$-form $T$ of
the horizontal distribution $T^H P$.
We also denote by
\begin{equation} \label{5.42}
A^{Bismut}_s \: : \: C^\infty_{\omega}(P; \widehat{E}) 
\rightarrow
\Omega^*_\omega(G) \otimes_{C^\infty_\omega(G)}
C^\infty_\omega(P; \widehat{E}) 
\end{equation}
its extension to 
$C^\infty_{\omega}(P; \widehat{E})$.
One can use finite-propagation-speed estimates, along with
the bounded geometry
of $\{Z_x\}_{x \in G^{(0)}}$ and the Duhamel expansion as in
\cite[Theorem 9.48]{Berline-Getzler-Vergne (1992)},
to show that we obtain a well-defined element
$e^{- \: (A^{Bismut}_s)^2}$.

We now couple $A^{Bismut}_s$ to the connection $\nabla^{can}$
in order to obtain a superconnection
\begin{equation} \label{5.44}
A_s \: : \: C^\infty_{\omega}(P; \widehat{E}) \rightarrow
\Omega^*_\omega(G) \otimes_{C^\infty_\omega(G)}
C^\infty_{\omega}(P; \widehat{E}).
\end{equation}
Let ${\cal R}$ be the rescaling operator on 
$\Omega^{even}_\omega(G)_{ab}$ which multiplies an element of
$\Omega^{2k}_\omega(G)_{ab}$ by 
$(2 \pi i)^{-k}$.
Doing a Duhamel expansion around $e^{- \: (A^{Bismut}_s)^2}$ and using the
fact that $h$ has compact support, we can define
\begin{equation} \label{5.45}
e^{- \: A_s^2} \in
\Hom^\infty_{C^\infty_\omega(G)} 
\left(
C^\infty_{\omega}(P; \widehat{E}), \:
\Omega^*_\omega(G) \otimes_{C^\infty_\omega(G)}
C^\infty_{\omega}(P; \widehat{E})
\right).
\end{equation}
and hence also define
${\cal R} \: \Tr_{s, <e>} \left( e^{- \: A_s^2} \right) \in 
\Omega^*_\omega(G)_{ab}$.
From the superconnection formalism 
\cite[Chapter 1.4]{Berline-Getzler-Vergne (1992)},
${\cal R} \: \Tr_{s, <e>} \left( e^{- \: A_s^2} \right)$ is closed and its
cohomology class is independent of $s \: > \: 0$; see
\cite[Theorem 3.1]{Heitsch (1995)} for a detailed proof in the
analogous case of ${\cal R} \: \Tr_s \left(e^{- \: (A^{Bismut}_s)^2} \right)$.

The proof of the next theorem is analogous to that of Theorem \ref{Theorem 2}.

\begin{theorem} \label{Theorem 5}
\begin{equation} \label{5.46}
\lim_{s \rightarrow 0} {\cal R} \: 
\Tr_{s, <e>} \left( e^{- \: A_s^2} \right) \: = \:
\int_Z \phi(z) \: \widehat{A} \left( \nabla^{TZ} \right) \: 
\ch \left( \nabla^{\widehat{V}} \right) \: \ch \left(\nabla^{can} \right)
\: \in \: \Omega^*_\omega(G)_{ab}.
\end{equation}
\end{theorem}

Let us note that the right-hand-side of (\ref{5.46}) pairs with closed graded
traces on $\Omega_c^*(G)$, and not just closed graded traces on
$\Omega^*_\omega(G)$. Let $\eta$ be a closed graded trace on
$\Omega_c^*(G)$.


Let $EG$ denote the bar construction of a universal space on which
$G$ acts freely. That is, $EG$ is the geometric realization of a
simplicial manifold given by $E_nG \: = \: G^{(n+1)}$, with
face maps
\begin{equation} \label{5.47}
d_i(\gamma_0, \ldots, \gamma_n) \: = \: 
\begin{cases}
(\gamma_1, \ldots, \gamma_n) &\text{ if } i = 0,\\
(\gamma_0, \ldots, \gamma_{i-1} \gamma_i, \ldots, \gamma_n) &
\text{ if } 1 \: \le \:  i \: \le \: n
\end{cases}
\end{equation}
and degeneracy maps
\begin{equation} \label{5.48}
s_i(\gamma_0, \ldots, \gamma_n) \: = \: 
(\gamma_0, \ldots, \gamma_{i}, 1, \gamma_{i+1}, \ldots \gamma_n),
\: \: \: \: \: 0 \: \le \: i \: \le \: n.
\end{equation}
Here $1$ denotes a unit.
The action of $G$ on $EG$ is induced from the action on $E_nG$ given
by $(\gamma_0, \ldots, \gamma_n) \: \gamma \: = \:
(\gamma_0, \ldots, \gamma_n \gamma)$. Let $BG$ be the quotient space.
Let $\pi \: : \: EG \rightarrow G^{(0)}$ be the map induced from
the maps $E_nG \rightarrow G^{(0)}$ given by 
$(\gamma_0, \ldots, \gamma_n) \rightarrow s(\gamma_n)$.
Let $J \in C(EG)$ be the ``barycentric coordinate'' corresponding
to the units $G^{(0)} \: \subset \: E_0G$. 
 That is, for each $x \in G^{(0)}$, $\pi^{-1}(x)$ is a simplicial
complex and $J \big|_{\pi^{-1}(x)}$ is the function on
$\pi^{-1}(x)$ defined as in \cite[(94)]{Lott (1992)}, with respect to the 
vertex $x$ instead of the vertex $e$.  Then for all $p \in EG$,
$\sum_{\gamma \in G^{\pi(p)}} J(p  \gamma) \: = \: 1$. 
Let $\nabla^{univ}$ be the connection constructed as in (\ref{5.14}), using
$J$ in place of $h$. Then pairing $\ch \left( \nabla^{univ} \right)$ with
$\eta$, we construct an element
$\Phi_\eta \in \HH^*_\tau(BG)$. 

Put $M \: = \: P /G$, a compact manifold. 
It inherits a foliation ${\cal F}$ from the submersion
$\pi \: : \: P \rightarrow G^{(0)}$.
Let $\nu \: : \: M \rightarrow BG$ be the classifying map for the
$G$-action on 
$P$. Put $V \: = \: \widehat{V}/G$, a vector bundle on $M$.
By naturality, 
\begin{equation} \label{5.49}
\langle \int_Z \phi(z) \: \widehat{A} \left( \nabla^{TZ} \right) \: 
\ch \left( \nabla^{\widehat{V}} \right) \: \ch \left(\nabla^{can} \right),
\eta \rangle \: = \:
\int_M \widehat{A}(T{\cal F}) \: \ch(V) \: \nu^* \Phi_\eta.
\end{equation}

As in the proof of Theorem \ref{Theorem 4}, we obtain
\begin{theorem} \label{Theorem 6}
Let $\eta$ be a closed graded trace on $\Omega^*_\omega(G)$. Then
\begin{equation} \label{5.50}
\langle \ch(\Ind(D)), \widetilde{\eta} \rangle \: = \:
\int_M \widehat{A}(T{\cal F}) \: \ch(V) \: \nu^* \Phi_\eta.
\end{equation}
\end{theorem}
\noindent
{\bf Remark : } Theorem \ref{Theorem 6} also follows from
\cite[Section III.7.$\gamma$, Theorem 12]{Connes (1994)}.

\begin{corollary} \label{Corollary 4}
Let $M^n$ be a compact manifold with a codimension-$q$ foliation ${\cal F}$.
Let $V$ be a vector bundle on $M$ and let
$D$ be a leafwise Dirac-type operator coupled to $V$.
Let $\HH^*(\Tr {\cal F})$ 
denote the Haefliger cohomology of $(M, {\cal F})$ \cite{Haefliger (1980)}. 
Recall that there is a linear map $\int_{\cal F} \: : \: \HH^*(M) \rightarrow
\HH^{*-n+q}(\Tr {\cal F})$. 
Let $\eta$ be a holonomy-invariant closed  transverse current.
Then 
\begin{equation} \label{5.51}
\langle \ch(\Ind(D)), \widetilde{\eta} \rangle \: = \:
\langle \int_{\cal F} \widehat{A}(T{\cal F}) \: \ch(V), \eta \rangle.
\end{equation}
\end{corollary}
\begin{pf}
 Let ${\cal H}$ be the holonomy groupoid of ${\cal F}$, with
source and range maps $s_{\cal H}, r_{\cal H} : \: {\cal H} \rightarrow M$ 
\cite[Section II.8.$\alpha$]{Connes (1994)}. 
Let $T$ be a complete transversal for ${\cal F}$. That is, $T$ is
a $q$-dimensional submanifold of $M$, not necessarily connected, which
is transverse to ${\cal F}$ and has the property that every leaf of
$(M, {\cal F})$ intersects $T$. Put $G \: = \: r_{\cal H}^{-1}(T) \: \cap \:
s_{\cal H}^{-1}(T)$, the reduced holonomy groupoid.
That is, an element of $G$ is an equivalence class of smooth leafwise paths in
$M$ from $T$ to $T$, where two paths are equivalent if they have the
same endpoints and the same holonomy. The units are $G^{(0)} \: = \: T$.

Put $P \: = \: s_{\cal H}^{-1}(T)$. 
Define $\pi \: : \: P \rightarrow G^{(0)}$ to be the
restriction of $s_{\cal H}$ to $P$. Then for $x \in T$, 
$\pi^{-1}(x)$ 
is the holonomy cover of the leaf through $x$, which we give the induced
Riemannian metric. One can see that $G$ acts freely, properly and
cocompactly on $P$. 

Put $L \: = \: P$, $\sigma \: = \: \pi$ and
let $i \: : \: G \rightarrow L$ be the
inclusion from $r_{\cal H}^{-1}(T) \: \cap \:
s_{\cal H}^{-1}(T)$ to $s_{\cal H}^{-1}(T)$. It is easy to check that
$(L, \sigma, i)$ satisfies the requirements to define 
$\Omega^*_\omega(G)$. 
Then $\eta$ defines a closed graded
trace on $\Omega^*_\omega(G)$. The right-hand-side of (\ref{5.50}) becomes
the right-hand-side of (\ref{5.51}). 
\end{pf}
{\bf Remark : } In order to prove
Corollary \ref{Corollary 4}, we do not have to assume that the holonomy
groupoid is Hausdorff. This is because the pairing with
the transverse current $\eta$ amounts to an integration over
$G^{(0)} \: = \: T$. Because of this we are effectively dealing with
forms of the type $\Omega^{*,0}_\omega(G)$, and so the Hausdorffness of $G$ 
does not play a role.  \\ \\
{\bf Remark : }
To see the relationship between Corollary \ref{Corollary 4}
and Connes' index theorem for
a foliation with a holonomy-invariant transverse measure $\mu$
\cite[Section I.5.$\gamma$, Theorem 7]{Connes (1994)}, let
$RS \in \HH_{n-q}(M; \R)$ denote the Ruelle-Sullivan current associated to
$\mu$
\cite[Section I.5.$\beta$]{Connes (1994)}. Then
$\langle \int_{\cal F} \widehat{A}(T{\cal F}) \: \ch(V), \mu \rangle
\: = \: \langle \widehat{A}(T{\cal F}) \: \ch(V), RS \rangle$. \\ \\ 
{\bf Remark : } In some cases of foliations, 
a heat equation proof of Corollary
\ref{Corollary 4}, using the Bismut superconnection,
 was given in \cite{Heitsch-Lazarov (1999)}.

\begin{corollary} \label{Corollary 5}
Let ${\cal A}$ be a subalgebra of the reduced groupoid $C^*$-algebra
$C^*_r(G)$ which is stable under the holomorphic functional
calculus in $C^*_r(G)$ and which contains
$C^\infty_\omega(G)$. Let $\eta$ be a closed graded trace on
$\Omega^*_c(G)$ which extends to give a cyclic cocycle on
${\cal A}$. Suppose that $TZ$ is spin and that $g^{TZ}$ has fiberwise
positive scalar curvature.
Then $\int_M \widehat{A}(T{\cal F}) \: \nu^* \Phi_\eta \: = \: 0$.
\end{corollary}

Suppose that $\dim(G^{(0)}) \: = \: 1$, with $G$ acting on 
$G^{(0)}$ so as to preserve orientation.
Let $v \in \Omega^1(G^{(0)})$ be a volume form.
With a hopefully-clear notation,
define a closed graded trace on $\Omega^*_c(G)$, concentrated on
$\Omega^{0,2}_c(G)$, by
\begin{equation} \label{5.52}
\eta(\omega) \: = \: \int_{\gamma_0 \gamma_1 \gamma_2 \in G^{(0)}}
\omega(\gamma_0, \gamma_1, \gamma_2) \: 
\left( \ln \frac{v}{\gamma_0 \gamma_1 \cdot v} \: d\ln 
\frac{\gamma_0 \gamma_1 \cdot v}{\gamma_0 \cdot v} \: - \: 
\ln \frac{\gamma_0 \gamma_1 \cdot v}{\gamma_0 \cdot v} \: d\ln 
\frac{v}{\gamma_0 \gamma_1 \cdot v} \right).
\end{equation}
Then $\nu^* \Phi_\eta$ is proportionate to the Godbillon-Vey class
${\rm GV} \in \HH^3(M; \R)$ 
\cite[Chapter III.6.$\beta$, Theorem 17]{Connes (1994)}.
Furthermore, the hypotheses of Corollary \ref{Corollary 5} are satisfied
\cite[Chapter III.7.$\beta$]{Connes (1994)}.

\begin{corollary} \label{Corollary 6}
Suppose that $\dim(G^{(0)}) \: = \: 1$, $G$ acts on $G^{(0)}$ so as
to preserve orientation,
$TZ$ is spin and $g^{TZ}$ has fiberwise
positive scalar curvature.
Then $\int_M \widehat{A}(T{\cal F}) \: {\rm GV} \: = \: 0$.
\end{corollary}

In general, if $\dim(G^{(0)}) \: = \: q$ and the action of $G$ on 
$G^{(0)}$ is
orientation-preserving then one can write down a closed
graded trace $\eta$ on $\Omega^*_c(G)$ so that
$\nu^* \Phi_\eta$ is proportionate to the Godbillon-Vey class
${\rm GV} \in \HH^{2q+1}(M; \R)$, and the above results extend.

\begin{corollary} \label{Corollary 7}
Let $M$ be a compact manifold with a codimension-$q$ foliation ${\cal F}$.
We assume that the foliation is tranversally orientable and that
$T{\cal F}$ is spin. We also assume that
the holonomy groupoid of the foliation is Hausdorff. Let $g^{T{\cal F}}$ be
a leafwise metric on $(M, {\cal F})$. If $g^{T{\cal F}}$ has
positive scalar curvature on the leaves then
$\int_M \widehat{A}(T{\cal F}) \: GV \: = \: 0$.
\end{corollary}
\begin{pf}
Let $T$ be a complete transversal for ${\cal F}$.
Let $G$ be the reduced holonomy groupoid.
Then the corollary is an
application of Theorem \ref{Theorem 6}. 
\end{pf}
{\bf Remark : }  Corollary \ref{Corollary 7} also follows
from \cite[Section III.7.$\beta$, Corollary 10]{Connes (1994)}.

\subsection{Translation} \label{Subsection 6.1}

In this subsection we show how the results of Section \ref{Section 6} 
specialize
to those of Section \ref{Section 3}, in the case when the groupoid comes from
the action of $\Gamma$ on $B$. We use the notation of Section \ref{Section 3}.

We put $G \: = \: B \times \Gamma$, with $G^{(0)} \: \cong \: B$,
$r(b, \gamma) \: = \: b$ and $s(b, \gamma) \: = \: b \gamma$.
A form 
\begin{equation} \label{5.53}
\sum_{g_0, \ldots, g_n} \eta_{g_0,\ldots,g_n} \: g_0 \:
dg_1 \ldots dg_n \in
\Omega^*(B, \C \Gamma)
\end{equation}
gets translated to the form $\omega \in \Omega^*(G)$ given by
\begin{equation} \label{5.54}
\omega((b_0, g_0), \ldots, (b_n, g_n)) \: = \:
\eta_{g_0,\ldots,g_n}(b_0).
\end{equation}
Then the product (\ref{5.3}) is equivalent to the calculation
\begin{align} \label{5.55}
& \left( \left( \sum_{g_0, \ldots, g_n} \eta_{g_0,\ldots,g_n} \: g_0 \:
dg_1 \ldots dg_n \right) \cdot
\left( \sum_{g^\prime_0, \ldots, g^\prime_{n^\prime}} 
\eta^\prime_{g^\prime_0,\ldots,g^\prime_{n^\prime}} \: g^\prime_0 \:
dg^\prime_1 \ldots dg^\prime_{n^\prime} \right) \right)(b) \: = \notag \\
& \sum_{g_0, \ldots, g_n, g^\prime_0, \ldots, g^\prime_{n^\prime}}
\eta_{g_0,\ldots,g_n}(b) \: 
\eta^\prime_{g^\prime_0,\ldots,g^\prime_n}(b  g^\prime_0 \ldots 
g^\prime_{n^\prime}) \:
g_0 \: dg_1 \ldots dg_n \: 
g^\prime_0 \:
dg^\prime_1 \ldots dg^\prime_{n^\prime} \: = \notag \\
& \sum_{g_0, \ldots, g_n, g^\prime_0, \ldots, g^\prime_{n^\prime}}
\eta_{g_0,\ldots,g_n}(b) \: 
\eta^\prime_{g^\prime_0,\ldots,g^\prime_n}(b g^\prime_0 \ldots 
g^\prime_{n^\prime}) \: \cdot \notag \\
& \left[ g_0 \: dg_1 \ldots d(g_n
g^\prime_0) \:
dg^\prime_1 \ldots dg^\prime_{n^\prime} \: + \ldots + \:  
(-1)^n \: g_0 g_1 \: dg_2 \ldots dg_n \: 
dg^\prime_0 \ldots dg^\prime_{n^\prime}
\right].
\end{align}
The differential $d$, given by
\begin{align} \label{5.56}
d \left( \sum_{g_0, \ldots, g_n} \eta_{g_0,\ldots,g_n} \: g_0 \:
dg_1 \ldots dg_n \right) \: = \: &
\sum_{g_0, \ldots, g_n} \left( d^{de \: Rham} \: \eta_{g_0,\ldots,g_n} \right)
\: g_0 \: dg_1 \ldots dg_n \: + \notag \\
& (-1)^{|\eta|} \:
\sum_{g_0, \ldots, g_n} \eta_{g_0,\ldots,g_n} \: 1 \: dg_0 \:
dg_1 \ldots dg_n
\end{align}
becomes the sum of $d_1$ and the differential $d_2$ of (\ref{5.4}).

Take $P \: = \: \widehat{M}$. Then ${\cal G} \: = \: \widehat{M} \times 
\Gamma$. The product (\ref{5.6}) becomes
\begin{equation} \label{5.57}
((\sum_g f_g \: g) F)(p) \: = \: \sum_g \: f_g(p) \: F(p g).
\end{equation}
We illustrate the right-hand-side of (\ref{5.57}) by the diagram
$p \stackrel{g}{\longleftarrow} pg$.
Equation (\ref{5.8}) is the translation of
\begin{align} \label{5.58}
(\sum_{g_0,g_1} f_{g_0,g_1} \: g_0 \: dg_1) \: F \: & = \:
\left( \sum_{g_0,g_1} \: g_0 \: dg_1 \: (g_0 g_1)^{-1} \cdot 
f_{g_0,g_1} \right)
\: F \notag \\
& = \: 
\left( \sum_{g_0,g_1} \: [d(g_0 g_1) \: - \: dg_0 \: g_1] \: 
(g_0 g_1)^{-1} \cdot f_{g_0,g_1} \right) \cdot F \notag \\
& = \:
\sum_{g_0,g_1} \: d(g_0 g_1) \: \left( (g_0 g_1)^{-1} \cdot f_{g_0,g_1}
\right) \: F
- \sum_{g_0,g_1} dg_0 
\left( g_0^{-1} \cdot f_{g_0,g_1} \right) \: (g_1 \cdot F) \notag \\
& = \:
\sum_{g_0 g_1 \: = \: g} \: dg \: \left( g^{-1} \cdot f_{g_0,g_1}
\right) \: F
- \sum_{g, g^\prime} dg 
\left( g^{-1} \cdot f_{g,g^\prime} \right) \: (g^\prime \cdot F),
\end{align}
or
\begin{equation} \label{5.59}
\left( (\sum_{g_0,g_1} f_{g_0,g_1} \: g_0 \: dg_1) \: F \right)(p) \: = \:
\sum_{g_0 g_1 \: = \: g} \: dg \: f_{g_0,g_1}(p g^{-1})
\: F(p)
- \sum_{g, g^\prime} dg 
f_{g,g^\prime}(pg^{-1}) \: F(p g^\prime).
\end{equation}
We illustrate the right-hand-side of (\ref{5.59}) by the diagrams
$pg^{-1} \stackrel{g_0}{\longleftarrow} p g_1^{-1} 
\stackrel{g_1}{\longleftarrow}
p$ and 
$pg^{-1} \stackrel{g}{\longleftarrow} p \stackrel{g^\prime}{\longleftarrow}
pg^\prime$.
Equation (\ref{5.11}) is the translation of
\begin{align} \label{5.60}
(\sum_{g} f_{g} \: g) \: (\sum_{g^\prime} dg^\prime 
{\cal F}_{g^\prime}) \: 
& = \:
\sum_{g,g^\prime} f_g \: g \: dg^\prime \: {\cal F}_{g^\prime} \notag \\
& = \: \sum_{g,g^\prime} g \: dg^\prime  \: ((g g^\prime)^{-1} \cdot f_g) 
\: {\cal F}_{g^\prime} \notag \\
& = \: \sum_{g,g^\prime} [d(g g^\prime) \: - \: dg \: 
g^\prime ]  \: ((g g^\prime)^{-1} \cdot f_g) 
\: {\cal F}_{g^\prime} \notag \\
& = \: \sum_{g,g^\prime} d(g g^\prime)
\: ((g g^\prime)^{-1} \cdot f_g) 
\: {\cal F}_{g^\prime} \: - \: \sum_{g,g^\prime}
dg
\: (g^{-1} \cdot f_g) 
\: (g^\prime \cdot {\cal F}_{g^\prime})
\end{align}
or
\begin{align} \label{5.61}
\left( (\sum_{g} f_{g} \: g) \: (\sum_{g^\prime} dg^\prime \: 
{\cal F}_{g^\prime}) \right)(p) \: = \: &
\sum_{gg^\prime \: = \: g_0} dg_0
\: f_g(p g_0^{-1}) 
\: {\cal F}_{g^\prime}(p) \: - \notag \\
& \sum_{g} dg
\: \cdot f_g(p g^{-1}) 
\: \sum_{g^\prime} {\cal F}_{g^\prime}(p g^\prime).
\end{align}
We illustrate the right-hand-side of (\ref{5.61}) by the diagrams
$pg_0^{-1} \stackrel{g}{\longleftarrow} p (g^\prime)^{-1}
\stackrel{g^\prime}{\longleftarrow} p$ and 
$pg^{-1} \stackrel{g}{\longleftarrow} p \stackrel{g^\prime}{\longleftarrow}
pg^\prime$.
Equation (\ref{5.12}) is the translation of
\begin{align} \label{5.62}
(\eta_{g_0,\ldots,g_n} \: g_0 \: dg_1 \ldots dg_n) \: F \: & = \: 
\left[ d(g_0 g_1) \ldots dg_n \: - \: dg_0 \: d(g_1 g_2) \ldots dg_n \: + \: 
\ldots \: + \: \right. \notag \\
& \: \: \: \: \: \left. (-1)^{n-1} \: dg_0 \ldots d(g_{n-1} g_n) \: + \:
(-1)^n \: dg_0 \ldots dg_{n-1} \: g_n \right] \notag \\
& \: \: \: \: \: 
\left((g_0 \ldots g_n)^{-1} \cdot \eta_{g_0,\ldots,g_n} \right) 
\: F \notag \\
& = \:
\left[ d(g_0 g_1) \ldots dg_n \: - \: dg_0 \: d(g_1 g_2) \ldots dg_n \: + \:
\ldots \: + \: \right. \notag \\
& \left. \: \: \: \: \: (-1)^{n-1} \: dg_0 \ldots d(g_{n-1} g_n)
 \right] \: 
\left((g_0 \ldots g_n)^{-1} \cdot \eta_{g_0,\ldots,g_n} \right) \: F \: + 
\notag \\
& \: \: \: \: \:
(-1)^n \: dg_0 \ldots dg_{n-1} \: 
\left( (g_0 \ldots g_{n-1})^{-1} \cdot \eta_{g_0,\ldots,g_n} \right)
(g_n \cdot F).
\end{align}
Equation (\ref{5.15}) is the translation of
\begin{equation} \label{5.63}
\left( \nabla^{(0,1)} F \right)(p) \: = \:
\left( \sum_g dg \: h \: (g^{-1} \cdot F) \right)(p) \: = \:
\sum_g dg \: h(p) \: F(p g^{-1}). 
\end{equation}
We illustrate the right-hand-side of (\ref{5.63}) by the diagram
$pg^{-1} \stackrel{g}{\longleftarrow} p$.
Equations (\ref{5.16}) and (\ref{5.17}) are the translations of (\ref{1.39}).
Equation (\ref{5.22}) is the translation of (\ref{2.12}).
Equation (\ref{5.29}) is the translation of 
\begin{align} \label{5.64}
\Tr_{<e>}(K)(b) \: & = \: \sum_{g_0,g_1, \ldots, g_l \in \Gamma :
g_0 \cdots g_l \: = \: e} 
(dg_1 \ldots dg_l) \: g_0 \: \left( \int_{Z_b} \phi(w) \:
K_{g_1, \ldots, g_l}
(w g_0^{-1}, w)
\: d\vol_{Z_b}(w) \right) \notag \\
& = \: \sum_{g_0,g_1, \ldots, g_l \in \Gamma :
g_0 \cdots g_l \in \Gamma_b} 
[dg_1 \ldots d(g_l g_0) \: - \: 
dg_1 \ldots d(g_{l-1} g_l) dg_0 \: + \ldots + \notag \\
& \: \: \: \: \: \: \: \: \: \: \: \: \: \: \:
(-1)^{l-1} \: d(g_1 g_2) \ldots dg_l \: dg_0 \: + \: (-1)^l \:
g_1 \: dg_2 \ldots dg_0 ]  \notag \\
& \: \: \: \: \: \: \left( \int_{Z_b} \phi(w) \:
K_{g_1, \ldots, g_l}
(w g_0^{-1}, w)
\: d\vol_{Z_b}(w) \right).
\end{align}

Choose a finite generating set for $\Gamma$. Let ${\cal C}$ be
the corresponding Cayley graph, on which $\Gamma$ acts on the right
by isometries.  If $B$ is compact, put $L \: = \: B \times {\cal C}$. Let
$i \: : \: G \rightarrow L$ be the natural inclusion
$B \times \Gamma \rightarrow B \times {\cal C}$. 
 (In this case, the requirements on
$L$ and $i$ are satisfied because $\Gamma$ is 
finitely-generated.)  Then
$C^\infty_\omega(G) \: = \: C^\infty(B, {\cal B}^\omega)$,
$\Omega^*_\omega(G) \: = \: \Omega^*(B, {\cal B}^\omega)$ and 
$C^\infty_\omega(P) \: = \: C^\infty_{{\cal B}^\omega}(\widehat{M})$.

We have $EG \: = \: E\Gamma \times B$. If $p_0 \: : \:
E\Gamma \times B \rightarrow E\Gamma$ is the projection map then
$J \: = \: p_0^* (j)$.

\end{document}